\documentclass{compositio}
\usepackage{amsmath}
\title[A two-variable refinement of Stark's conjecture]{A
two-variable refinement of the Stark conjecture in the function field     
case}
\author{Greg W. Anderson}
\email{gwanders@math.umn.edu}
\address{School of Mathematics, University of Minnesota, MN 55455, U.S.A.}
\classification{11R58, 11G45}
\keywords{Stark conjecture, Tate's thesis, Stirling formula, soliton} 
\date{May 26, 2005}
\DeclareMathOperator{\Aut}{{\mathrm{Aut}}}
\newcommand{\MMM}{{\mathcal{M}}}
\newcommand{\Pbold}{{\mathbf{P}}}
\newcommand{\HHH}{{\mathcal{H}}}
\newcommand{\ordbold}{{\mathbf{ord}}}

\newcommand{\gggg}{{\mathfrak{g}}}
\newcommand{\II}{{\mathcal{I}}}
\newcommand{\JJ}{{\mathcal{J}}}

\newcommand{\weight}{{\mathrm{wt}}}
\newcommand{\ab}{{\mathrm{ab}}}
\newcommand{\card}{{\#}}
\newcommand{\UUU}{{\mathcal{U}}}
\DeclareMathOperator{\St}{{\mathrm{St}}}
\DeclareMathOperator{\lead}{{\mathrm{lead}}}
\DeclareMathOperator{\NNN}{{\mathcal{N}}}
\DeclareMathOperator{\KKK}{{\mathcal{K}}}
\DeclareMathOperator{\LLL}{{\mathcal{L}}}
\DeclareMathOperator{\FFF}{{\mathcal{F}}}
\DeclareMathOperator{\ord}{{\mathrm{ord}}}
\DeclareMathOperator{\tr}{{\mathrm{tr}}}
\DeclareMathOperator{\Res}{{\mathrm{Res}}}
\DeclareMathOperator{\Gal}{{\mathrm{Gal}}}

\newcommand{\sch}{{\mathrm{Sch}}}
\newcommand{\iso}{{\stackrel{\sim}{\rightarrow}}}
\newcommand{\kk}{{\bf k}}
\newcommand{\KK}{{\bf K}}
\newcommand{\adeles}{{\mathbb{A}}}
\newcommand{\FF}{{\mathbb{F}}}

\newcommand{\OO}{{\mathcal{O}}}
\newcommand{\norm}[1]{{\left\|#1\right\|}}
\newcommand{\ZZ}{{\mathbb Z}}
\newcommand{\ee}{{\bf e}}
\newcommand{\CC}{{\mathbb C}}

{\theoremstyle{plain}
\newtheorem{Conjecture}[subsection]{Conjecture}
\newtheorem{Lemma}[subsection]{Lemma}
\newtheorem{Theorem}[subsection]{Theorem}
\newtheorem{Proposition}[subsection]{Proposition}
\newtheorem{Corollary}[subsection]{Corollary}
}
{\theoremstyle{remark}

}
\newcommand{\one}{{\mathbf{1}}}
\newcommand{\perf}{{\mathrm{perf}}}
\DeclareMathOperator{\Moore}{{\mathrm{Moore}}}
\DeclareMathOperator{\Ore}{{\mathrm{Ore}}}
\newcommand{\toy}{{\mathrm{toy}}}
\newcommand{\RR}{{\mathbb{R}}}
\newcommand{\sep}{{\mathrm{sep}}}
\DeclareMathOperator{\Hom}{{\mathrm{Hom}}}

\newcommand{\abold}{{\mathbf{a}}}
\newcommand{\CCC}{{\mathcal{C}}}
\DeclareMathOperator{\Frob}{{\mathrm{Frob}}}
\DeclareMathOperator{\graph}{{\mathrm{graph}}}
\newcommand{\VVV}{{\mathcal{V}}}
\begin{document}
\begin{abstract}
We propose a conjecture refining the Stark
conjecture St($K/k$,$S$) in the function field case. 
Of course St($K/k$,$S$) in this case is a theorem due to Deligne
and independently to Hayes. 
The novel
feature of our conjecture is that it involves {\em two}- rather than
one-variable algebraic functions over finite fields.  Because the
conjecture is formulated in the language and framework of Tate's thesis
we have powerful standard techniques at our disposal to
study it.
 We build a case for our conjecture by (i)  proving a parallel
result in the framework of adelic harmonic analysis which we dub the
{\em adelic Stirling formula}, (ii) proving the conjecture in the genus
zero case and (iii) explaining in detail how to deduce St($K/k$,$S$)
from our conjecture. In the genus zero case the class of two-variable
algebraic functions thus arising includes all the {\em solitons} over a genus zero global field previously studied and applied by the author,
collaborators and others.  Ultimately the inspiration for this paper
comes from striking examples given some years ago by R.\ Coleman and D.\
Thakur. 
\end{abstract}
\maketitle 
\tableofcontents
\section{Introduction}
We propose here a conjecture (Conjecture~\ref{TheConjecture} below)
refining the  Stark conjecture $\St(K/k,S)$ (Tate's formulation
\cite{Tate}) in the function field case, and we verify our conjecture
in the genus zero case. Of course
$\St(K/k,S)$ in the function field case is a theorem due to Deligne \cite{Tate}
and independently to Hayes \cite{Hayes}.  The main novelty of our conjecture
is that it involves  {\em two-} rather than one-variable algebraic
functions.
 The conjecture is formulated in the powerful language and framework of Tate's
thesis. Ultimately this
paper is inspired by remarkable examples of Coleman
\cite{Coleman} and Thakur \cite[\S9.3]{Thakur}.
To motivate the paper we devote
most of this introduction to the discussion of an example 
similar to the ones originally considered by Coleman and by Thakur. 
 We
first discuss the example in ``raw form'' and then we explain how to
look at it from the adelic point of view cultivated in this paper.

Consider the smooth projective model $C/\FF_q$
of the affine plane curve
$$Y^q-Y=X^{q-1}$$
over $\FF_q$. On the surface $C\times C$ (product
formed over
$\FF_q$) consider the meromorphic function
$$\varphi=Y\otimes 1-1\otimes Y-\frac{X\otimes 1}{1\otimes
X}.$$
Following the lead of Coleman, it can easily be verified
that
\begin{equation}\label{equation:ColemanDivisor}
\begin{array}{c}
\mbox{divisor of $\varphi$}\\\\
\equiv \;\mbox{transpose of}\\\\
\graph((X,Y)\mapsto
(X^q,Y^q))\\\\
\displaystyle+\sum_{\alpha\in \FF_q^\times}
\graph((X,Y)\mapsto (\alpha X,\alpha+Y)),\\\\
\mbox{modulo horizontal and vertical divisors.}
\end{array}
\end{equation}
By the theory of correspondences on curves it follows that
the analogously formed sum of endomorphisms of the Jacobian of $C$
vanishes. Thus  we reprove following Coleman---by a strikingly
elementary argument---that the ``universal Gauss
sum''
 $$-\sum_{\alpha\in \FF_q^\times}((X,Y)\mapsto (\alpha X, \alpha+Y))\in
\ZZ[\Aut(C/\FF_q)]
$$
acts on the Jacobian of
$C$ in the same way as does the $q^{th}$ power Frobenius endomorphism.

Now view $C$  as a covering of the $T$-line over $\FF_q$ by setting
$$-T=Y^q-Y=X^{q-1}.$$
Still following Coleman, notice that by specializing the
second variable to $\FF_{q^d}$ in the example
$\varphi$,  we obtain a ``Stark-like'' function
in the extension
$\FF_{q^d}(X,Y)/\FF_{q}(T)$.   We say ``Stark-like'' because  the
functions arising this way do not have precisely the right properties
to qualify as Stark units in the sense of
$\St(K/k,S)$; nonetheless, these functions
would appear to be promising ``raw material'' for building Stark units
since their divisors and Galois properties are evident. What is yet more intriguing is that this approach gives insight into the variation of the Stark unit as we vary the choice of completely split place.
It is Coleman's brilliant idea to link two-variable algebraic
functions to Stark's conjecture in this way. 

We turn  to follow Thakur's lead. We have
\begin{equation}\label{equation:ThakurInterpolation}
\left.\begin{array}{l}
\\
\varphi\\
\\
\end{array}\right|\begin{array}{l}
X\otimes 1\mapsto \alpha X\\
Y\otimes 1\mapsto \beta+Y\\
1\otimes X\mapsto X^{q^{N+1}}\\
1\otimes Y\mapsto Y^{q^{N+1}}
\end{array}=
\displaystyle\prod_{\begin{subarray}{c}
a\in \FF_q[T]\\
\deg a<N
\end{subarray}}
\frac{a+\beta T^N+T^{N+1}+\alpha/T}{a+T^N},
\end{equation}
for all $\alpha\in \FF_q^\times$, $\beta\in \FF_q$ and 
integers $N\geq 0$, as can be verified by manipulating Moore
determinants. By definition, the left side of the equation is
the pull-back of $\varphi$ via the map
$$((X,Y)\mapsto((\alpha X,\beta+Y),(X^{q^{N+1}},Y^{q^{N+1}}))):C\rightarrow C\times C.$$
The system of identities (\ref{equation:ThakurInterpolation}) uniquely
characterizes our example
$\varphi$, and illustrates an extremely valuable heuristic which we
credit to Thakur, namely:  ``Frobenius interpolation formulas''  with
``$\Gamma$-partial-products''  on the right side automatically have
``functions of Coleman type'' on the left side.

We pause  to remark that a fairly extensive body of theory and applications of phenomena  of
Coleman/Thakur type has already been developed. For example, see
\cite{AndersonStickelberger},
\cite{AndersonAHarmonic},  \cite{ABP},
\cite{Sinha}, \cite{SinhaDR}; the applications to transcendence theory have been particularly notable. For a survey with a wealth of examples, see  
\cite[Chaps.\ 4 and 8]{ThakurBook}.  
While the main emphasis in the theory to this point has been on the case of the basefield $\FF_q(T)$, significant efforts to generalize the base have been made (see \cite{AndersonAHarmonic} and \cite{ThakurBook}). But a general framework for understanding and classifying all these examples has been lacking.
The main point of this paper is to rebuild the theory on the foundation of Tate's thesis so that 
a clear (albeit conjectural) picture of the situation over a general basefield emerges.

We return to the analysis of our example $\varphi$. We explain how
to look at the example from the point of view peculiar to this paper. 
Put
$$k=\FF_q(T),\;\;\;K=\FF_q(X,Y),\;\;\;G=\Gal(K/k).$$
The group $G$ is a copy of $\FF_q^\times \times \FF_q$ and in
particular is abelian. Notice
that
$\varphi$ is  invariant under the diagonal action
$$(x\otimes y)^\sigma=x^\sigma\otimes y^\sigma$$ of $G$ on
$K\otimes_{\FF_q}K$. 

The example $\varphi$ lives in a setting with somewhat subtle
commutative algebra properties. Notice that (i)
the ring extension 
$$K\otimes_{\FF_q}K/k\otimes_{\FF_q}k$$
is a finite etale extension
of Dedekind domains, (ii) the diagonal prime
$$\Delta_K=\ker(x\otimes y\mapsto
xy):K\otimes_{\FF_q}K\rightarrow K
$$ lies above the diagonal prime
$$\Delta=\Delta_k=\ker(x\otimes y\mapsto xy):k\otimes_{\FF_q}k\rightarrow
k
$$ and (iii) the prime of the diagonally $G$-invariant subring of
$K\otimes_{\FF_q}K$ below $\Delta_K$ has the same residue field as $\Delta$, namely $k$. The upshot is
that the diagonal completion of
$k\otimes_{\FF_q}k$  coincides with the diagonal completion of the
diagonally
$G$-invariant subring of $K\otimes_{\FF_q}K$.
We now simply identify the latter with the former.
Under this identification the strange-looking formula
$$\varphi=\sum_{i=0}^\infty (T\otimes 1-1\otimes T)^{q^i}
-\prod_{i=0}^\infty \left(\frac{T\otimes 1}{1\otimes
T}\right)^{-q^i}.
$$
makes sense.

Now let's put global class field theory into the picture. 
Let
$\adeles$ be the adele ring of
$k$ and for each idele $a\in \adeles^\times$, let $\norm{a}\in q^\ZZ$ be the corresponding idele norm. Let 
$\rho:\adeles^\times\rightarrow G$
 be the
reciprocity law homomorphism of global class field theory, redefined the modern way, as in Tate's article \cite{TateBackground}, so that uniformizers map
to geometric Frobenius elements and hence, for all $c$ algebraic over $\FF_q$ and $a\in \adeles^\times$, we have
$c^{\rho(a)}=c^{\norm{a}}$. (We follow this rule throughout the paper.)
We make
$\adeles^\times$ act on the perfect ring
$\sqrt[q^\infty]{K}\otimes_{\FF_q}\sqrt[q^\infty]{K}$
by the rule
$$(x\otimes y)^{(a)}=x^{\rho(a)}\otimes y^{\norm{a}},$$
which makes sense since 
$\Gal(K/k)=\Gal(\sqrt[q^\infty]{K}/
\sqrt[q^\infty]{k})$.  Since the
$\adeles^\times$-action defined above commutes with the diagonal action
of
$G$, we may naturally view all the functions
$$\varphi^{(a)}\;\;\;(a\in
\adeles^\times)
$$
as elements of the perfection 
of the diagonal completion of
$k\otimes_{\FF_q}k$. 

We need a few more definitions before getting to the ``punchline''.
Let
$\ord_\Delta$ be the normalized additive valuation of
$k\otimes_{\FF_q}k$ giving the order of vanishing at the diagonal
ideal, and let this additive valuation be extended to the perfection
of the completion in evident fashion. Let
$\OO\subset\adeles$ be the maximal compact subring.  Let $\infty$ be
the unique place of
$k$ at which $T$ has a pole. Let $\tau\in \adeles^\times$ be the idele
whose component at $\infty$ is $T^{-1}$ and whose components elsewhere
are
$1$. Consider the
$\ZZ$-valued Schwartz function
$\Phi=\one_{\tau(T+1/T+\tau\OO)}-\one_{\tau(1+\tau\OO)}$
on $\adeles$,
where $\one_S$ is probabilist's notation for the function taking the
value $1$ on the set $S$ and $0$ elsewhere. 

We now
 get to the heart of the matter. 
The divisor 
formula (\ref{equation:ColemanDivisor}) can  be rewritten as
 the formula
\begin{equation}\label{equation:ApproximateConjecture1}
\ord_\Delta
\varphi^{(a)}=\sum_{x\in
k}\Phi(a^{-1}x)
\end{equation} holding for all $a\in \adeles^\times$,
plus a further assertion to the effect that no primes of the Dedekind
domain
$K\otimes_{\FF_q}K$ other than those of the form
$$\ker\left((x\otimes y\mapsto x^{\rho(a)}y^{\norm{a}}):
K\otimes_{\FF_q}K\rightarrow \sqrt[q^\infty]{K}\right)\;\;\;(a\in \adeles^\times)$$
enter into the prime factorization of
$\varphi$ in $K\otimes_{\FF_q}K$. The interpolation formula
(\ref{equation:ThakurInterpolation}) can  be rewritten as the
formula
\begin{equation}\label{equation:ApproximateConjecture2}
\ord_\Delta\left(\varphi^{(a)}-1\otimes\prod_{x\in
k^\times}x^{\Phi(a^{-1}x)}\right)>0
\end{equation} holding for all $a\in \adeles^\times$
such that $\norm{a}>1$. Formulas~(\ref{equation:ApproximateConjecture1})
and (\ref{equation:ApproximateConjecture2}) very strongly suggest that there's
more to the Coleman/Thakur phenomenon than just a few isolated examples.
 We now draw our
discussion of the example
$\varphi$ to a close.

We turn to the task of describing our conjecture, in order at least to convey some its flavor.   We continue with similar notation.
Let $k$ be any
global field with constant field
$\FF_q$, let $\adeles$ be the adele ring of $k$, and let $\rho:\adeles^\times\rightarrow \Gal(k^\ab/k)$ be the
reciprocity law homomorphism of global class field theory.
We define a {\em Coleman unit} to be an element of the
fraction field of the ``big ring''
 $\sqrt[q^\infty]{k^\ab}\otimes_{\FF_q^\ab}\sqrt[q^\infty]{k^\ab}$ which is (i)
invariant under the diagonal action $(x\otimes y)^\sigma=x^\sigma\otimes y^\sigma$ of $\Gal(k^\ab/k)$ and (ii) a unit at every maximal ideal of the big ring
not transformable under the action $(x\otimes y)^{(a)}=x^{\rho(a)}\otimes y^{\norm{a}}$
of $\adeles^\times$ to the diagonal ideal $\ker(x\otimes y\mapsto xy)$.
Our conjecture associates to every
$\ZZ[1/q]$-valued Schwartz-function $\Phi$ on $\adeles$
such that $\Phi(0)=0=\hat{\Phi}(0)$ a Coleman unit $\varphi$
for which  analogues of (\ref{equation:ApproximateConjecture1}) and
(\ref{equation:ApproximateConjecture2}) hold, the latter in a more
general formulation (not subject to any restriction on $\norm{a}$)
relating the ``leading Taylor coefficient''  to a gadget we call the
 {\em Catalan symbol}.  The values of the Catalan symbol are function field
elements roughly analogous in structure  to the Catalan numbers
$\frac{1}{n+1}\left(\begin{subarray}{c}
2n\\n
\end{subarray}\right)$. The Catalan symbol  is defined in terms of a further
gadget we call the {\em rational Fourier transform}. The
characteristic feature of the theory of rational Fourier transforms is
that the function
$(x\mapsto \delta_{x0}-\delta_{x1}):\FF_q\rightarrow\{-1,0,1\}$
is assigned the role usually played by a nonconstant complex-valued
character of $\FF_q$.

As mentioned above, we prove our
conjecture in genus zero. We believe that ideas from the author's
papers
\cite{AndersonAHarmonic} and \cite{AndersonFock} suitably combined
have a fighting chance to prove our conjecture in general. 
Perhaps the rather different set of ideas from
\cite{AndersonSpecialPoints} could also yield a proof. We hope for a
proof of the conjecture in the general case not so heavily computational
as the proof we give here in the genus zero case.

To make the point that our conjecture refines $\St(K/k,S)$, we provide a
deduction of the latter from the former in \S\ref{section:StarkRecovery}
below. In fact a large fraction of the paper is devoted to developing the
machinery needed to make the deduction of $\St(K/k,S)$ from our conjecture
run smoothly.
In particular, the technical tool we call the {\em adelic Stirling formula}
(Theorem~\ref{Theorem:AdelicStirling} below) is crucial for proving
our ``recipe'' (Theorem~\ref{Theorem:Recipe} below) for the Stark unit
conditional on our conjecture.

We conclude the introduction by remarking on the organization of the
paper. Since we have provided a table of contents, we needn't make
comments section-by-section. We just offer the following advice.
After glancing at
\S\ref{subsection:GeneralNotation} to be apprised of general notation
and terminology, the reader could very well start in
\S\ref{section:GlobalStirling}, because that's where the main story-line
of the paper begins; the preceding sections, which are more or less
independent of each other, could be
treated as references. But on the other hand, if the reader would take the
trouble to work patiently through the sections before
\S\ref{section:GlobalStirling}, he/she would
be prepared to hear the main story undistracted by technical
issues of an essentially nonarithmetical nature.

\section{General notation and
terminology}\label{subsection:GeneralNotation}
\subsection{}
We denote the cardinality of a set $S$ by $\card S$.
\subsection{}
Given sets $S\subset X$, 
let $\one_S$ denote the real-valued function taking the value $1$
on $S$ and $0$ on $X\setminus S$; we omit reference to $X$ in the
notation because it can always be inferred from context. We borrow this
notation from the probabilists.
\subsection{}
{\em Rings}
are always commutative with unit. Let $A^\times$ denote the group of
units of a ring $A$. Let $\FF_q$ denote a field of $q<\infty$ elements.

\subsection{}
A function $v$ on a ring $A$
taking values in $\RR\cup\{+\infty\}$ is called an {\em additive
valuation}
 if $v(a)=+\infty\Leftrightarrow a=0$,
$v(a+b)\geq\min(v(a),v(b))$ and
$v(ab)=v(a)+v(b)$ for all $a,b\in A$.  If
$A$ is a field and $v(A)=\ZZ\cup\{+\infty\}$, we say that $v$
is {\em normalized}. Additive valuations are said to be {\em
equivalent} if proportional. A {\em place} of a field $k$ is an
equivalence class of nontrivial additive
valuations of $k$.

\subsection{}
Given an integral domain $A$ of characteristic $p>0$, let
$A_\perf$ be the closure of $A$ under the extraction of $p^{th}$
roots, i.~e., the direct limit of the system 
$A\xrightarrow{x\mapsto x^p}A
\xrightarrow{x\mapsto x^p}\dots$ of rings and homomorphisms. 
If $A=A_\perf$, then we say that $A$ is {\em perfect}.

\subsection{}
Given a locally compact
totally disconnected topological space $X$, we denote by $\sch(X)$ the
{\em Schwartz space}
of locally constant compactly supported
complex-valued functions on $X$.
(In the literature this is sometimes called
instead the {\em Schwartz-Bruhat space}, cf.\ \cite{Ramakrishnan}.)

\subsection{Moore determinants}\label{subsection:MooreDeterminants}
Given an $\FF_q$-algebra $A$ and ring elements $x_1,\dots,x_n\in A$, 
we define the {\em Moore determinant}
$$\Moore(x_1,\dots,x_n)=
\left|\begin{array}{cccc}
x_1^{q^{n-1}}&x_1^{q^{n-2}}&\dots&x_1\\
\vdots&\vdots&&\vdots\\
x_n^{q^{n-1}}&x_n^{q^{n-2}}&\dots&x_n\\
\end{array}\right|=\det_{i,j=1}^n x_i^{q^{n-j}}\in A.$$
For this we have the well-known {\em Moore identity}
(see
\cite[Chap.\ 1]{Goss})
\begin{equation}\label{equation:MooreIdentity}
\Moore(x_1,\dots,x_n)=
\prod_{\ell=1}^n\sum_{c_{\ell+1}\in \FF_q}
\cdots \sum_{c_{n}\in \FF_q}\left(
x_\ell+\sum_{i=\ell+1}^{n}c_ix_i\right).
\end{equation}
For example, working over the field of rational functions in a variable
$T$ with coefficients in $\FF_q$,  we have
\begin{equation}\label{equation:DN}
\prod_{\begin{subarray}{c}
a\in \FF_q[T]\\
\deg a=N\\
a:\,\mbox{\tiny monic in}\,T
\end{subarray}}a=
\frac{\Moore(T^N,\dots,1)}{\Moore(T^{N-1},\dots,1)}=\prod_{i=0}^{N-1}(T^{q^N}-T^{q^i})
\end{equation}
for all nonnegative integers $N$, where $\deg a$ denotes the degree of $a$ as a polynomial in $T$. The last equality 
is justified by the Vandermonde identity.

\subsection{Ore determinants}
The following variant of the Moore determinant will also be needed.
Given an $\FF_q$-algebra $A$, ring elements $x_1,\dots,x_n\in A$, and an $\FF_q$-linear functional
$\xi$ defined on the $\FF_q$-span of $x_1,\dots,x_n$, put
$$\Ore(\xi, x_1,\dots,x_{n})=\left|\begin{array}{cccc}
\xi x_1&x_1^{q^{n-2}}&\dots&x_1\\
\vdots&\vdots&&\vdots\\
\xi x_n&x_n^{q^{n-2}}&\dots&x_n\\
\end{array}\right|\in A.$$
We call $\Ore(\xi, x_1,\dots,x_{n})$ the {\em Ore determinant} of $\xi,x_1,\dots,x_n$.
For this we have a variant of the Moore identity which,
for simplicity, we state under some special assumptions always
fulfilled in practice. We suppose now that $A$ is a field, and that
$x_1,\dots,x_n$ are $\FF_q$-linearly independent, in which case
$\Moore(x_1,\dots,x_n)\neq 0$.
Let $V$ be the $\FF_q$-linear span of $x_1,\dots,x_n$.
We suppose further that $\xi V\neq 0$. Then, so we claim,
\begin{equation}\label{equation:OreIdentity}
\frac{\Ore(\xi,x_1,\dots,x_n)}{\Moore(x_1,\dots,x_n)}
=\prod_{\begin{subarray}{c}
v\in V\\
\xi v=1
\end{subarray}}v^{-1}.
\end{equation}
We call this relation the {\em Ore identity}.
To prove the claim, first note that the ratio on the left side
does not change if we replace $x_1,\dots,x_n$ by any $\FF_q$-basis of $V$.
We may therefore assume without loss of generality that $\xi x_i=\delta_{i1}$
for $i=1,\dots,n$, in which case the claim follows from the Moore identity
(\ref{equation:MooreIdentity}). The claim is proved.
Relation (\ref{equation:OreIdentity}) is a convenient 
way to restate what for our purposes is the main point of
Ore-Elkies-Poonen duality (see \cite[\S4.14]{Goss}).

\section{The rational Fourier transform and the Catalan symbol (``toy'' versions)}
\label{section:Duality} 
We discuss in ``toy form'' some simple algebraic
notions later to be applied in the adelic context.

\subsection{Notation}
We fix vector spaces $V$ and $V^*$ over $\FF_q$ of the same finite 
dimension. We fix a perfect $\FF_q$-bilinear pairing
$$\langle\cdot,\cdot\rangle:V\times V^*\rightarrow \FF_q.$$
We fix a nonconstant homomorphism
$$\lambda:\FF_q\rightarrow U(1),$$
 where $U(1)$ is the group of complex
numbers of absolute value $1$.
We also fix a field $K$ containing a copy of $\FF_q$.

\subsection{Toy Fourier transforms}

Let $\sch(V)$ be the set of $\CC$-valued functions 
on $V$.  
 Given $\Phi\in \sch(V)$, put
$$\hat{\Phi}(v^*)=\FFF[\Phi](v^*)=\sum_{v\in
V}\Phi(v)\lambda(-\langle v,v^*\rangle)$$ for all $v^*\in V^*$, thus
defining the {\em toy Fourier
transform}
$$\hat{\Phi}=\FFF[\Phi]\in \sch(V^*).$$
Perhaps it is wasteful to define two notations for the Fourier
transform, but it is quite convenient, and we follow this practice
throughout the paper for all the versions of the Fourier transform that
we consider.

For example, given
$v_0\in V$ and an $\FF_q$-subspace $W\subset V$, we have
\begin{equation}\label{equation:ToyFourierExample}
\FFF[\one_{v_0+W}](v^*)=\card
W\cdot
\one_{W^\perp}(v^*)\cdot\lambda(-\langle v_0,v^*\rangle)
\end{equation}
for all $v^*\in V^*$,
where $W^\perp$ is the subspace of $V^*$ annihilated by $W$.
Further, we have
$$\Phi(v)=(\card
V)^{-1}\sum_{v^*\in V^*}\hat{\Phi}(v^*)\lambda(\langle v,v^*\rangle)
$$
for all $\Phi\in \sch(V)$ and $v\in V$; thus 
we can invert the toy Fourier transform explicitly.
So far all this is completely familiar.

\subsection{The function $\lambda_0$}\label{subsection:TrivialIdentities}
We begin to warp things a bit.
Put
$$\lambda_0=\left(x\mapsto\left\{\begin{array}{rl}
1&\mbox{if $x=0$}\\
-1&\mbox{if $x=1$}\\
0&\mbox{if $x\neq 0,1$}
\end{array}\right.\right):\FF_q\rightarrow\{-1,0,1\}\subset \CC.
$$
Then we have trivial identities
\begin{equation}\label{equation:lambdaIdentities}
\begin{array}{rcl}
\lambda_0(x)&=&q^{-1}\sum_{c\in
\FF_q^\times}(1-\lambda(c))\lambda(-cx),\\
\lambda(-x)&=&-\sum_{c\in\FF_q^\times}
\lambda(-c)\lambda_0(c^{-1}x),\\
\sum_{c\in\FF_q^\times}\lambda(cx)&=&
\sum_{c\in\FF_q^\times}\lambda_0(cx),\\
-\sum_{c\in
\FF_q^\times}\lambda(-c)\lambda(cx)
&=& q\lambda_0(x)-
\sum_{c\in \FF_q^\times}\lambda_0(cx)
\end{array}
\end{equation}
holding for all
$x\in \FF_q$. Roughly speaking, these identities say that although
$\lambda_0$ is not a character of $\FF_q$ (save in the case $q=2$), it
behaves sufficiently like a character that we can base an alternative
theory of Fourier transforms on it.

\subsection{Toy rational Fourier transforms}
\label{subsection:ToyRatDef}
 Given $\Phi\in \sch(V)$, put
$$\tilde{\Phi}(v^*)=\FFF_0[\Phi](v^*)=\sum_{v\in
V}\Phi(v)\lambda_0(\langle v,v^*\rangle)
$$ for all $v^*\in V^*$, thus defining the 
{\em toy rational
Fourier transform} 
$$\tilde{\Phi}=\FFF_0[\Phi]\in \sch(V^*).$$
Via (\ref{equation:lambdaIdentities}) we have
\begin{equation}\label{equation:ToyRelationOfFTs}
\begin{array}{rcl}
\hat{\Phi}(v^*)
&=&-\sum_{c\in \FF_q^\times}\lambda(-c)\tilde{\Phi}(c^{-1}v^*),\\
\tilde{\Phi}(v^*)&=&
q^{-1}\sum_{c\in \FF_q^\times}(1-\lambda(c))\hat{\Phi}(cv^*)
\end{array}
\end{equation}
for all $v^*\in V$; thus we can express toy Fourier transforms in terms
of their rational analogues and {\em vice versa}.
For example,
given
$v_0\in V$ and an $\FF_q$-subspace $W\subset V$, we have
\begin{equation}\label{equation:ToyRationalFourierExample}
\FFF_0[\one_{v_0+W}]=\card W\cdot 
\left(\one_{W^\perp\cap[\langle
v_0,\cdot\rangle=0]}-\one_{W^\perp\cap[\langle
v_0,\cdot\rangle=1]}\right)
\end{equation}
by combining
(\ref{equation:ToyFourierExample}),
(\ref{equation:lambdaIdentities}) and (\ref{equation:ToyRelationOfFTs}).
The notation
$[\langle v_0,\cdot\rangle=1]$ is probabilist-style 
shorthand for $\{v^*\in V^*\vert \langle v_0,v^*\rangle=1\}$.
We use similar notation below without further comment.
By combining the toy Fourier inversion formula
with (\ref{equation:lambdaIdentities})
and  (\ref{equation:ToyRelationOfFTs}) we have
$$
\begin{array}{rcl}
\card V\cdot \Phi(v)&=& \sum_{v^*\in V}\lambda(\langle
v,v^*\rangle)\left(-\sum_{c\in
\FF_q^\times}\lambda(-c)
\tilde{\Phi}(c^{-1}v^*)\right)\\
&=&\sum_{v^*\in V^*}\left(-\sum_{c\in
\FF_q^\times}\lambda(-c)\lambda(c\langle
v,v^*\rangle)\right)\tilde{\Phi}(v^*)\\
&=&\sum_{v^*\in V^*}\left(q\lambda_0(\langle v,v^*\rangle)-\sum_{c\in
\FF_q^\times}\lambda_0(c\langle v,v^*\rangle)
\right)\tilde{\Phi}(v^*)\end{array}
$$
for all $\Phi\in \sch(V)$ and $v\in V$; thus we can invert the toy
rational Fourier transform explicitly. 

Let $\ZZ[1/q]$ be the ring obtained from  $\ZZ$ by inverting $q$.
Let $\sch_0(V)$ be the group of $\ZZ[1/q]$-valued functions on $V$. Notice
that we have $\FFF_0[\sch_0(V)]=\sch_0(V^*)$.
It is because of the latter relation that $\FFF_0$ deserves to be called rational.

\subsection{Toy Catalan symbols} 
Recall that $K$ is a field containing $\FF_q$.
 Let
$\alpha:V\rightarrow K$ be an injective $\FF_q$-linear map. For all
$\Phi\in \sch_0(V)$, put
$$\left(\begin{array}{c}
\alpha\\
\Phi
\end{array}\right)=\prod_{0\neq v\in V}\alpha(v)^{\Phi(v)}\in
K_\perf^\times.
$$
The definition makes sense because $K_\perf^\times$ is a uniquely $q$-divisible group.
We call $\left(\begin{subarray}{c}
\cdot\\
\cdot
\end{subarray}\right)$ the  {\em toy Catalan symbol}.
It is clear that $\left(\begin{array}{c}
\alpha\\
\Phi
\end{array}\right)$ depends $\ZZ[1/q]$-linearly on $\Phi$.
We see the values of the toy Catalan symbol as being very roughly analogous in
structure to the {\em Catalan numbers}
$$\frac{1}{n+1}\left(\begin{array}{c}
2n\\
n
\end{array}\right)=\frac{(2n)!}{(n+1)!n!}=
\prod_{0\neq x\in \ZZ}x^{(\one_{(0,2n]}-\one_{(0,n+1]}-\one_{(0,n]})(x)},$$
whence the terminology.

\begin{Proposition}\label{Proposition:EasyMooreApp}
Let $\alpha:V\rightarrow K$ be an injective $\FF_q$-linear map.
Let
$W\subset V$ be an
$\FF_q$-subspace. Then the map
$$\left(v\mapsto 
\left\{
\begin{array}{cl}
\left(\begin{array}{c}
\alpha\\
\one_{v+W}
\end{array}\right)&\mbox{if $v\not\in W$}\\
0&\mbox{if $v\in W$}
\end{array}\right.\right):V\rightarrow K
$$
is $\FF_q$-linear.
\end{Proposition}
\proof Choose a basis $w_1,\dots,w_n$ for $W$ over $\FF_q$.
The $\FF_q$-linear map
$$\left(v\mapsto\frac{\Moore(\alpha(v),\alpha(w_1),\dots,
\alpha(w_n))}{\Moore(\alpha(w_1),\dots,\alpha(w_n))}\right):V\rightarrow
K$$ coincides with the map in question.
\qed
\begin{Lemma}
Let $\alpha:V\rightarrow K$ be an injective $\FF_q$-linear map.
Let
$W\subset V$ be an
$\FF_q$-subspace. Then the map
$$\left(v^*\mapsto\left\{\begin{array}{cl}
\left(\begin{array}{c}
\alpha\\
-\one_{W\cap [\langle \cdot,v^*\rangle=1]}
\end{array}\right)&\mbox{if $v^*\not\in W^\perp$}\\
0&\mbox{if $v^*\in W^\perp$}
\end{array}\right.\right):V^*\rightarrow K
$$
is $\FF_q$-linear.
\end{Lemma}
\proof We may assume without loss of generality that $W=V$. Choose an $\FF_q$-basis $v_1,\dots,v_n\in V$.
The $\FF_q$-linear map
$$\left(v^*\mapsto 
\frac{\Ore(v^*\circ \alpha^{-1},\alpha(v_1),\dots,\alpha(v_n))}
{\Moore(\alpha(v_1),\dots,\alpha(v_n))}\right):V^*\rightarrow K$$
coincides with the map in question.
\qed
\begin{Proposition}\label{Proposition:DualityClincher}
Let $\beta:V^*\rightarrow K$
be an injective
$\FF_q$-linear map. Let
$W\subset V$ be an
$\FF_q$-subspace. 
 Fix $v_0\in V\setminus W$. Then the map
$$\left(v\mapsto
\left\{\begin{array}{cl}
\left(\begin{array}{c}
\beta\\
\FFF_0[\one_{v+W}-\one_{v_0+W}]
\end{array}\right)&\mbox{if $v\not\in W$}\\
0&\mbox{if $v\in W$}
\end{array}\right.\right):V\rightarrow K$$
is $\FF_q$-linear.
\end{Proposition}
\noindent This delicate property of the toy rational Fourier transform
is the main formal justification for defining it as we have. 
In the adelic context there will be a similar phenomenon
(see Theorem~\ref{Theorem:GenericBehavior} below).
\proof For $v\in V\setminus W$ we have by
(\ref{equation:ToyRationalFourierExample}) and the definitions that
$$\begin{array}{cl}
&\left(\begin{array}{cl}
\beta\\
\FFF_0[\one_{v+W}-\one_{v_0+W}]
\end{array}\right)^{1/\card W}\\\\
=&\left(\begin{array}{c}
\beta\\
\one_{W^\perp\cap[\langle v,\cdot \rangle=0]}
-\one_{W^\perp\cap[\langle v,\cdot \rangle=1]}
-\one_{W^\perp\cap[\langle v_0,\cdot \rangle=0]}
+\one_{W^\perp\cap[\langle v_0,\cdot \rangle=1]}
\end{array}\right)\\\\
\stackrel{\star}{=}&\left(\begin{array}{c}
\beta\\
-\one_{W^\perp\cap[\langle v,\cdot \rangle\neq 0]}
-\one_{W^\perp\cap[\langle v,\cdot \rangle=1]}
+\one_{W^\perp\cap[\langle v_0,\cdot \rangle\neq 0]}
+\one_{W^\perp\cap[\langle v_0,\cdot \rangle=1]}
\end{array}\right)\\\\
=&\left(\begin{array}{c}
\beta\\
-\one_{W^\perp\cap[\langle v,\cdot \rangle=1]}
+\one_{W^\perp\cap[\langle v_0,\cdot \rangle=1]}
\end{array}\right)^q,
\end{array}$$
which by the preceding lemma proves the result. \qed
\subsection{Remark}
The star marks the main trick of the toy theory.
It is nothing more than the fact that for any three subsets $A$, $B$ and $C$
of a set $X$
we have 
$\one_{A\cap B}-\one_{A\cap C}=-\one_{A\setminus B}+\one_{A\setminus C}$.
We shall use the trick again to make explicit calculations 
in \S\ref{section:GenusZero}.

\section{A functor and related identities}
\label{section:FunctorsAndIdentities}
We work out identities in which on the left side appear functions on a
certain commutative affine group scheme over $\FF_q$ evaluated at
certain  special points, and in which on the right side appear
expressions interpretable as values of naturally occurring toy Catalan
symbols. These identities are used in the sequel for no purpose other than
as an inputs  to the proof in
\S\ref{section:GenusZero} of the genus zero case of the
conjecture to which we allude in the title of the paper. 
 Some similar identities were studied and applied in the papers
\cite{AndersonStickelberger}, \cite{AndersonAHarmonic} and  \cite{ABP}, but there are two
novel features here. Namely, (i) we work out connections with duality and (ii) we
allow for the possibility of infinite ramification at infinity.

\subsection{The functor $\HHH$ and associated structures}
\subsubsection{}\label{subsubsection:ResidueIntroduction}
Let $t$ be a variable. 
Given an $\FF_q$-algebra $R$, let $R((1/t))$ be the ring consisting of all power series
$\sum_i a_it^i$ with coefficients
$a_i\in R$ vanishing for $i\gg 0$, and put $\Res{t=\infty}\left(\sum_i a_it^idt\right)=-a_{-1}\in R$. 
Let $R[t]\subset R((1/t))$ (resp., $R[[1/t]]\subset R((1/t))$) be the subring consisting of 
power series $\sum_i a_it^i$
such that $a_i=0$ unless $i\geq 0$ (resp., $i\leq 0$).
When $R$ is a field, we denote by $R(t)$ the field of rational functions in $t$
with coefficients in $R$, and
we identify $R((1/t))$
with the completion of $R(t)$ at the infinite place.

\subsubsection{}
We define a representable group-valued functor $\HHH$ of
$\FF_q$-algebras
$R$ by the rule
$$\HHH(R)=\lim_{\leftarrow}(R[t]/m(t))^\times\times (1+(1/t)R[[1/t]]),$$
where the inverse limit is extended over monic $m(t)\in \FF_q[t]$ ordered
by the divisibility relation, and the group law is induced by
multiplication. The commutative affine group scheme
$\HHH$ is the natural one to consider in connection with the problem of
constructing the maximal abelian extension of a field of rational functions
in one variable with coefficients in $\FF_q$.
Later, in \S\ref{section:GenusZero}, we discuss and apply
the arithmetical properties of $\HHH$. But now we focus on
properties more in the line with classical symmetric function theory.

\subsubsection{}
By definition,
to give an $R$-valued point 
$P\in \HHH(R)$ is to give a power series
$$P_\infty(t)\in 1+(1/t)R[[1/t]]$$
and for each monic $m=m(t)\in \FF_q[t]$ a congruence class
$$P_m(t)\bmod{m(t)}\in (R[t]/m(t))^\times,$$
subject to the constraints that
$$m\vert n\Rightarrow P_m(t)\equiv P_n(t)\bmod{m(t)}$$
for all monic $m,n\in \FF_q[t]$. For definiteness, we always
choose the representative $P_m(t)\in R[t]$ to be of least
possible degree in its
congruence class modulo $m(t)$.

\subsubsection{}
We produce useful examples of $R$-valued points of $\HHH$ 
as follows. Let
$M(t)\in R[t]$ be a monic polynomial of degree $d$ such
that for all monic $m(t)\in \FF_q[t]$, the resultant of $M(t)$
and $m(t)$ is a unit of $R$.
We then define $[M(t)]\in \HHH(R)$
by setting $[M(t)]_\infty=M(t)/t^d$
and
$$[M(t)]_m=\mbox{remainder of $M(t)$ upon division by $m(t)$}
$$
for each monic $m=m(t)\in \FF_q[t]$. Given for $i=1,2$ a polynomial $M_i(t)\in R[t]$ 
for which $[M_i(t)]\in \HHH(R)$ is defined, note that
$[M_1(t)][M_2(t)]=[M_1(t)M_2(t)]$.

\subsubsection{}\label{subsubsection:Brackets}
 Given 
$x=x(t)\in \FF_q(t)$ and $P\in \HHH(R)$,
put
\begin{equation}\label{equation:HThetaDef}
\theta_x(P)=\Res_{t=\infty}x(t)(P_m(t)-P_\infty(t)/t)\,dt\,\in
R
\end{equation} where
$m=m(t)\in \FF_q[t]$ is any monic polynomial such that
$mx\in \FF_q[t]$. It is easy to verify that the right side of
(\ref{equation:HThetaDef}) is independent of $m$ and hence $\theta_x(P)$ is
well-defined. Clearly, $\theta_x(P)$ depends $\FF_q$-linearly on $x$. Now write
$$x=\langle x\rangle +\lfloor x\rfloor\;\;\;(\langle x\rangle
\in (1/t)\FF_q[[1/t]],\;\;\lfloor x\rfloor \in \FF_q[t]).$$
As the notation suggests, we think of $\langle x\rangle$ as the
``fractional part'' and $\lfloor x\rfloor$ as the ``integer part'' of
$x$. With $P$, $x$ and $m$ as above, we have
\begin{equation}\label{equation:HThetaDefBis}
\theta_x(P)=\Res_{t=\infty}(\langle x\rangle (t)P_m(t)-\lfloor x\rfloor(t)P_\infty(t)/t)\,dt
\end{equation}
after throwing away terms which do not contribute to the residue.
Note that the right side of (\ref{equation:HThetaDefBis}) remains unchanged
if we replaced $P_m(t)$ by any member of its congruence class in $R[t]$
modulo $m(t)$.
If moreover $P=[M(t)]$, where $M(t)\in R[t]$ is a monic polynomial of degree $d$ whose resultant
with each monic element of $\FF_q[t]$ is a unit of $R$, then
we have
\begin{equation}\label{equation:SimplestHHHExample}
\theta_x([M(t)])=
\Res_{t=\infty} (\langle x\rangle(t)-
\lfloor x\rfloor(t) t^{-d-1})M(t)dt,
\end{equation}
after replacing $P_m(t)=[M(t)]_m$ in (\ref{equation:HThetaDefBis}) by $M(t)$.

In the next three propositions we study the values of the natural transformations $\theta_x$ at certain special points of $\HHH$.

\begin{Proposition}\label{Proposition:JacobiTrudiRedux}
Let $T$ be a variable independent of $t$.
Let $N$ be a nonnegative integer. Fix nonzero $x\in \FF_q(t)$
and put 
$$X=X(t)=\langle x\rangle +t^N\lfloor x\rfloor\in \FF_q(t)^\times.$$
Consider the point
$$P=\left[\prod_{i=0}^N(t-T^{q^i})\right]^{-1}\in
\HHH(\FF_q(T)).$$ Then we have
\begin{equation}\label{equation:CompleteSymmetric}
\theta_x(P)= 
\prod_{\begin{subarray}{c}
a\in \FF_q[T]\\
\deg a<N
\end{subarray}}
\frac{X(T)+a}{T^N+a}\in
\FF_q(T)^\times.
\end{equation}
\end{Proposition}
\proof 
To simplify writing, put $K=\FF_q(T)$ and $T_i=T^{q^i}$. Let $X_0,\dots,X_N$ be independent variables, independent also of $t$ and $T$, and put
$$F(X_0,\dots,X_N)=\frac{\left|\begin{array}{ccc}
X_N&\dots&X_0\\
T_N^{N-1}&\dots&T_0^{N-1}\\
\vdots&&\vdots\\
1&\dots&1
\end{array}\right|}{\left|\begin{array}{ccc}
T_N^{N}&\dots&T_0^{N}\\
\vdots&&\vdots\\
1&\dots&1
\end{array}\right|}\in K[X_0,\dots,X_N].$$
Choose monic $m=m(t)\in \FF_q[t]$ such that $mx\in \FF_q[t]$ and write
$$\prod_{i=0}^N(t-T_i)\cdot
P_m(t)=1-W(t)m(t)\;\;\;(W(t)\in K[t]).$$ 
Then we have, so we claim,
$$\begin{array}{cl}
&\theta_x(P)\\
=&\displaystyle
\Res_{t=\infty}
\left(\langle x\rangle(t)\left(\frac{1-W(t)m(t)}{\prod_{i=0}^N(t-T_i)}
\right)-\lfloor
x\rfloor(t)\, t^{-1}\prod_{i=0}^N(1-T_i/t)^{-1}\right)\,dt\\\\
=&\displaystyle-\Res_{t=\infty}
\left(\frac{\langle x\rangle(t)W(t)m(t)+t^N\lfloor
x\rfloor(t)}{\prod_{i=0}^N(t-T_i)}dt\right)\\\\
=&F(X(T_0),\dots,X(T_N))\\\\
=&\displaystyle\frac{\Moore(X(T),T^{N-1},\dots,1)}
{\Moore(T^N,\dots,1)}.
\end{array}
$$
This chain of equalities is justified as follows.
We get the first equality
by  plugging into  version (\ref{equation:HThetaDefBis}) of the
definition of $\theta$. We get the
second equality after throwing away the term
$\frac{\langle x\rangle(t)}{\prod_{i=0}^N(t-T_i)}dt$, which does not
contribute to the residue,  and then rearranging in evident fashion.
Since we have
$$\frac{1}{ \prod_{i=0}^N(t-T_i)}=F\left(\frac{1}{t-T_0},\dots,\frac{1}{t-T_N}\right),$$
we get the penultimate equality by ``sum of residues
equals zero'' for meromorphic differentials on the $t$-line over $K$. The last equality follows directly from the definitions. The claim is proved. Via the Moore determinant identity (\ref{equation:MooreIdentity}),
 the result follows. \qed

\begin{Proposition}\label{Proposition:DualJacobiTrudiRedux}
Let $T$ and $x$  be as in the preceding proposition.
 Let $\nu$ be a nonnegative integer.
Put $$X=X(t)=\langle x\rangle+\lfloor t^{-\nu-1}x\rfloor\in \FF_q(t).$$
Consider the point
$$P=\left[\prod_{i=0}^{\nu-1}(t-T^{q^i})\right]\in
\HHH(\FF_q(T)).$$
Consider the $\FF_q$-linear functionals
$$\left.\begin{array}{rcr}
\xi&=&(a(T)\mapsto \Res_{t=\infty}(\langle
x\rangle(t)-t^{-\nu-1}\lfloor x\rfloor(t))\,a(t)\,dt)\\
\xi_1&=&(a(T)\mapsto -\Res_{t=\infty}t^{-\nu-1}\,a(t)\,dt)
\end{array}\right\}:V\rightarrow \FF_q,$$
where
$$V=(\textup{\mbox{$\FF_q$-span of $1,T,\dots,T^\nu$}})\subset \FF_q(T).$$
Then we have
\begin{equation}\label{equation:ElementarySymmetric}
\theta_x(P)=\left\{\begin{array}{cl}
\displaystyle\prod_{\begin{subarray}{c}
a\in V\\
\xi_1a=1
\end{subarray}}a\bigg /
\prod_{\begin{subarray}{c}
a\in V\\
\xi a=1
\end{subarray}}a\in \FF_q(T)^\times&\mbox{if $\xi\neq 0$,}\\\\
0&\mbox{if $\xi=0$.}
\end{array}\right.
\end{equation}
Moreover, we have
\begin{equation}\label{equation:ThetaVanishing1}
\theta_x(P)=0\Leftrightarrow  \langle x\rangle-t^{-\nu-1}\lfloor
x\rfloor\in t^{-\nu-2}\FF_q[[1/t]]+\FF_q[t]
\end{equation}
 and 
 \begin{equation}\label{equation:ThetaVanishing2}
 \theta_x(P)=0\Rightarrow X\neq 0.
 \end{equation}
\end{Proposition}
\proof Write
$$\prod_{i=0}^{\nu-1}(t-T^{q^i})=\sum_{i=0}^\nu (-1)^i e_{i}t^{\nu-i}.$$
Then $e_i$ for $i=1,\dots,\nu$ is the $i^{th}$ elementary symmetric function of 
 $T,T^q,\dots,T^{q^{\nu-1}}$,
and $e_0=1$.
By specializing the well-known representation of the $i^{th}$ elementary symmetric
function as a quotient of determinants (see \cite[I,3]{Macdonald}) to the present case
we have
$$e_i=\frac{\Moore(T^{\nu},\dots,\widehat{T^{\nu-i}},\dots,1)}{\Moore(T^{\nu-1},\dots,1)}$$
for $i=0,\dots,\nu$; here the term bearing the ``hat'' is to be omitted,
and if $\nu=0$, the denominator (an empty determinant) is put to $1$.
By means of the Ore identity (\ref{equation:OreIdentity}) we can express the right side of the claimed identity as the ratio of
Ore determinants 
$$\frac{\Ore(\xi,T^\nu,\dots,1)}{\Ore(\xi_1,T^\nu,\dots,1)},$$
and the latter  (note that $\xi_1(T^i)=\delta_{i\nu}$ 
for $i=0,\dots,\nu$) can be brought
to the form $$\sum_{i=0}^\nu (-1)^{i}e_{i}\xi(T^{\nu-i})$$
by straightforward manipulation of determinants.
But this last equals $\theta_x(P)$ by (\ref{equation:SimplestHHHExample}).
Thus (\ref{equation:ElementarySymmetric}) is proved.
It follows that $\xi=0\Leftrightarrow \theta_x(P)=0$. Moreover, $\xi=0$
if and only if the condition on the right side of (\ref{equation:ThetaVanishing1}) is fulfilled.
Therefore (\ref{equation:ThetaVanishing1}) holds.

We turn, finally,
to the proof of (\ref{equation:ThetaVanishing2}). Consider
the Laurent expansion
$\sum_i b_it^i$ of $x$ at $t=\infty$, where $b_i\in \FF_q$ and $b_i=0$ for $i\gg 0$.
Since $x\neq 0$, not all the coefficients $b_i$ vanish. Moreover,
by (\ref{equation:ThetaVanishing1}), we have $b_{i-\nu-1}=b_i$ for $i=0,\dots,\nu$.  Finally, by definition
of $X$, 
the Laurent expansion of $X$ at $t=\infty$ is 
$\sum_{i<0}b_it^i+\sum_{i\geq 0}b_{i+\nu+1}t^{i}$. The latter Laurent expansion does not vanish identically since to form it we have merely suppressed some repetitions of digits
in the Laurent expansion $\sum_i b_it^i$. Thus relation (\ref{equation:ThetaVanishing2}) holds.
 \qed

\begin{Proposition}\label{Proposition:PeerIntoTheHole}
Let $T$, $x$, $\nu$, $X$ and $P$ be as
in the preceding proposition.
Assume now that $\theta_x(P)=0$ (and hence $X\neq 0$).
Put
$$\varpi=T\otimes 1-1\otimes T\in \FF_q(T)\otimes_{\FF_q}\FF_q(T).$$
(Note that the ring on the right is a principal ideal domain of which
$\varpi$ is a prime element.)    
Consider the point
$$\tilde{P}=\left[t-T^{q^\nu}\otimes 1\right]^{-1}
\left[\prod_{i=0}^\nu(t-1\otimes T^{q^i})\right]\in
\HHH(\FF_q(T)\otimes_{\FF_q}
\FF_q(T))$$
(which reduces modulo $\varpi$ to $P$).
We have
\begin{equation}\label{equation:ThetaLead2}
\theta_x(\tilde{P})\equiv
(1\otimes C)
\cdot \varpi^{q^\nu}
\bmod \varpi^{q^\nu+1},
\end{equation}
where
$$C=X(T^{q^\nu}) \prod_{i=0}^{\nu-1}(T^{q^\nu}-T^{q^i})\in \FF_q(T)^\times.$$
\end{Proposition}
\proof Put $K=\FF_q(T)$ to simplify writing.
Fix monic $m=m(t)\in \FF_q[t]$ such that $mx \in \FF_q[t]$.
Put $$Q=[t-T^{q^\nu}]^{-1}\in \HHH(K),$$
noting that
$$Q_\infty(t)=\left(1-T^{q^\nu}/t\right)^{-1}\in 1+(1/t)K[[1/t]],$$
and write
$$(t-T^{q^\nu}) Q_m(t)=1-W(t)m(t)\;\;\;(W(t)\in K[t]).$$  
Let 
$$1\otimes Q_m(t)\in (K\otimes_{\FF_q}K)[t],\;\;\;1\otimes Q_\infty(t)\in 1+(1/t)(K\otimes_{\FF_q}K)[[1/t]]$$ be the results of applying
the homomorphism 
$$(x\mapsto 1\otimes x):K\rightarrow K\otimes_{\FF_q}K$$ coefficient by coefficient to $Q_m(t)$ and $Q_\infty(t)$, respectively. Then, so we claim, we have congruences
$$\begin{array}{rcl}
\tilde{P}_m(t)&\equiv&\displaystyle (1+\varpi^{q^\nu}(1\otimes
Q_m(t)))\prod_{i=0}^{\nu-1}(t-1\otimes T^{q^i})
\bmod{(m(t),\varpi^{q^\nu+1})},\\\\
\tilde{P}_\infty(t)&\equiv&\displaystyle
\left(1+\varpi^{q^\nu}(1\otimes Q_\infty(t))/t\right)
\prod_{i=0}^{\nu-1}(1-(1\otimes T^{q^i})/t)\bmod{\varpi^{q^\nu+1}}.
\end{array}
$$
To verify the first congruence, multiply both sides by
$$t-T^{q^\nu}\otimes 1=(t-1\otimes T^{q^\nu})-\varpi^{q^\nu}$$ and
see that  both sides reduce to
$\prod_{i=0}^{\nu}(t-1\otimes T^{q^i})$ modulo
$(m(t),\varpi^{q^\nu+1})$. To verify the second congruence multiply both sides
by
$$1-\frac{T^{q^\nu}\otimes 1}{t} =1-\frac{1\otimes T^{q^\nu}}{t}-\frac{\varpi^{q^\nu}}{t}$$
and see that both sides reduce to $\prod_{i=0}^\nu\left(1-\frac{1\otimes T^{q^i}}{t}\right)$
modulo $\varpi^{q^\nu+1}$. Thus the claim is proved. 

Continuing our calculation, we now plug into version
(\ref{equation:HThetaDefBis}) of the definition of
$\theta$, taking into account our hypothesis that $\theta_x(P)=0$.
We
find that
$$\begin{array}{rcl}
C&=&\displaystyle \Res_{t=\infty}
\left(\langle x\rangle (t) Q_m(t)\prod_{i=0}^{\nu-1}(t-T^{q^i})\right.\\\\
&&\;\;\;\;\;\;\;\;\;\;\displaystyle \left.-
\lfloor x\rfloor(t)Q_\infty(t)t^{-2}\prod_{i=0}^{\nu-1}(1-T^{q^i}/t)\right)dt\\\\
&=&\displaystyle\Res_{t=\infty}
\frac{-\langle x\rangle(t)W(t)m(t)+\langle x\rangle(t)
-t^{-\nu-1}\lfloor
x\rfloor(t)}{(t-T^{q^\nu})}\prod_{i=0}^{\nu-1}(t-T^{q^i})dt,
\end{array}
$$
where $C\in K$ is the unique coefficient for which
(\ref{equation:ThetaLead2}) holds.
By (\ref{equation:ThetaVanishing1})
and another application of our hypothesis that $\theta_x(P)=0$, 
we have
$$C=-\Res_{t=\infty}
\frac{\langle x\rangle(t)W(t)m(t)+\lfloor t^{-\nu-1}
x\rfloor(t)
}{(t-T^{q^\nu})}\prod_{i=0}^{\nu-1}(t-T^{q^i})dt.
$$
Finally, we get the claimed value for $C$ by applying ``sum of residues equals zero''
for meromorphic differentials on the $t$-line over $K$.
\qed

\section{Two group-theoretical lemmas}
We prove a couple of technical results which in the sequel
are used for no purpose other than as inputs to the
proof of Theorem~\ref{Theorem:Recipe}. We
``quarantine'' the results here since they are general
group-theoretical facts whose proofs
would be a distraction from the main story.

\begin{Lemma}\label{Lemma:IdealTheoryInGroupRings}
Let $q>1$ be an integer.
Let $\Gamma$ be an abelian group equipped with a nonconstant homomorphism
$\norm{\cdot}:\Gamma\rightarrow q^\ZZ\subset \RR^\times$. Let $\ZZ[\Gamma]$ be the integral group ring of $\Gamma$ and let
$\JJ\subset\ZZ[\Gamma]$ be the kernel of the ring
homomorphism $\ZZ[\Gamma]\rightarrow \RR$ induced by $\norm{\cdot}$.
For every subgroup
$\Gamma'\subset\Gamma$, let
$\II(\Gamma')\subset \ZZ[\Gamma]$ be the ideal 
generated by differences of
elements of
$\Gamma'$.  
Let $\Pi\subset \Gamma$ be a subgroup of finite index and 
put $\Pi_1=\Pi\cap\ker\norm{\cdot}$. 
Then we have
$\II(\Pi)\cap \JJ=\II(\Pi_1)+\II(\Pi)\cdot \JJ$.
\end{Lemma}
\proof Since $\II(\Pi_1)\subset \II(\Pi)\cap\JJ$,  we may
pass to the quotient $\ZZ[\Gamma/\Pi_1]=\ZZ[\Gamma]/\II(\Pi_1)$ in
order to carry out our analysis of ideals. After replacing $\Gamma$ by $\Gamma/\Pi_1$,
we may simply assume that $\Pi_1=\{1\}$, in which case the function $\norm{\cdot}$
maps $\Pi$ isomorphically to a subgroup
of $q^\ZZ$. We cannot have $\Pi=\{1\}$ lest $[\Gamma:\Pi]=\infty$,
and hence $\Pi$ is a free abelian group of rank $1$. Let $\pi_0\in \Pi$ be a generator.
Then $\II(\Pi)$ is the principal ideal of $\ZZ[\Gamma]$ generated by $1-\pi_0$.
Suppose now that we are given $f\in \II(\Pi)\cap \JJ$.
Then we can write $f=(1-\pi_0)g$ for some $g\in \ZZ[\Gamma]$ and
since $1-\norm{\pi_0}\neq 0$, 
we must have $g\in \JJ$.  Therefore we have
$f\in \II(\Pi)\cdot \JJ$. \qed

\begin{Lemma}\label{Lemma:AbelianRootRedux}
Let $k$ be a  field. Fix an algebraic closure $\bar{k}/k$.
Let $k^\ab$ be the abelian closure of $k$ in $\bar{k}$.
Let $K/k$ be a finite subextension of $k^\ab/k$. Let $\mu(K)$
be the group of roots of unity in $K$ and assume that
$\card\mu(K)<\infty$. Put
$G=\Gal(K/k)$ and let $J\subset \ZZ[G]$ be the ideal of the integral
group ring annihilating $\mu(K)$. 
Let $\epsilon:\ZZ[G]\rightarrow \bar{k}^\times$
be a homomorphism of abelian groups such that
$\epsilon\vert_J\in \Hom_G(J,K^\times)$. 
Then $\epsilon$ takes all its values in $(k^\ab)^\times$.
\end{Lemma}
\noindent The lemma is a variant of
\cite[Lemma 6]{Stark} and \cite[p.\ 83, Prop.\ 1.2]{Tate}. 

\proof 
Let $k^\sep$ be the separable algebraic closure of $k$ in $\bar{k}$. 
Put $e=[\ZZ[G]:J]=\card\mu(K)$. Since $e\ZZ[G]\subset J$
and $e$ is prime to the characteristic of $\bar{k}$, the homomorphism
$\epsilon$ takes all its values in $(k^\sep)^\times$,
and moreover
the restriction map
$$(\psi\mapsto \psi\vert_J):\Hom(\ZZ[G],(k^\sep)^\times)\rightarrow
\Hom(J,(k^\sep)^\times)$$
is surjective. Put
$\gggg=\Gal(k^\sep/k)$. Regard
$\ZZ[G]$ and $J$ as (left) $\gggg$-modules by inflation.
Given two (left)
$\gggg$-modules $A$ and $B$, let the group  $\Hom(A,B)$
of homomorphisms of abelian groups be equipped with (left)
$\gggg$-module structure by the rule $(\sigma h)(a)=\sigma
(h(\sigma^{-1}a))$. Fix a generator $\zeta\in \mu(K)$.
From the exact sequence of $G$-modules
$$0\rightarrow J\rightarrow
\ZZ[G]\xrightarrow{\zeta\mapsto \zeta^{\abold}}\mu(K)\rightarrow 0$$ we
deduce an exact sequence 
$$0\rightarrow \ZZ/e\ZZ
\rightarrow  \Hom(\ZZ[G],(k^\sep)^\times)
\rightarrow\Hom(J,(k^\sep)^\times)\rightarrow 0$$
of $\gggg$-modules, where we make the evident identification
$$\ZZ/e\ZZ=\Hom_G(\mu(K),(k^\sep)^\times).$$
Since $\gggg$ acts trivially on $\ZZ/e\ZZ$, we can extract an exact
sequence
$$
\Hom_{\gggg}(\ZZ[G],(k^\sep)^\times)
\rightarrow
\Hom_{\gggg}(J,(k^\sep)^\times)\xrightarrow{\delta}
\Hom_{\mbox{\tiny loc.\ const.}}(\gggg,\ZZ/e\ZZ)
$$
of abelian groups from the long exact sequence in Galois cohomology. 
The group on the right is the group of locally constant homomorphisms
from $\gggg$ to $\ZZ/e\ZZ$.
When we make the boundary map $\delta$ explicit by a diagram-chase, we
find that
$$\epsilon(\sigma^{-1}\abold)^\sigma=
(\zeta^{\abold})^{\delta[\epsilon\vert_J](\sigma)}\epsilon(\abold)$$
for all $\abold\in \ZZ[G]$ and $\sigma\in \gggg$.
(In this nonabelian setting, even though
we employ exponential notation, we remain steadfastly leftist---we
follow the rule
$(x^\tau)^\sigma=x^{\sigma\tau}$.)
It follows that $\epsilon$ commutes with every
 $\sigma\in
\Gal(k^\sep/k^\ab)$, and since every such $\sigma$ acts trivially on
$\ZZ[G]$, it follows finally that $\epsilon$ takes all its values in
$(k^\ab)^\times$. 
\qed

\section{The local Stirling formula}\label{section:LocalStirling}
In this section we work in the setting of harmonic analysis on a
nonarchimedean local field. The main result here
is the {\em local Stirling formula} (Theorem~\ref{Theorem:LocalStirling}), which is a technical result
used in the sequel for no purpose other than as  an input to the proof of
the adelic Stirling formula (Theorem~\ref{Theorem:AdelicStirling}).
To motivate the local Stirling formula, we prove a corollary
(Corollary~\ref{Corollary:LocalStirling}, which is not needed in the sequel)
 which begins to explain the relationship with the classical Stirling formula.   See \cite{Ramakrishnan} for background on local harmonic analysis.

\subsection{Data}
\begin{itemize}
\item
Let $k$ be a nonarchimedean local field.
(It is not necessary to assume that $k$ is of positive characteristic.)
\item Let $\ee:k\rightarrow U(1)$ be a nonconstant
continuous homomorphism from the additive group of $k$ 
to the group of complex numbers of absolute value $1$. 
\end{itemize}
All the constructions in \S\ref{section:LocalStirling}
proceed naturally from these choices. 

\subsection{Notation}
\begin{itemize}
\item Let $\OO$ be the maximal compact subring of $k$.
\item Fix $\kappa\in k^\times$ such that $\kappa^{-1}\OO=\{\xi\in k\mid
\xi\OO\subseteq
\ker\ee\}$.
\item Let $\mu$ be Haar measure on $k$, normalized by $\mu\OO\cdot
\mu(\kappa^{-1}\OO)=1$.
\item Let $\mu^\times$ be Haar
measure on $k^\times$, normalized by $\mu^\times
\OO^\times=1$.
\item Let $q$ be the cardinality of the residue field of $\OO$.
\item For each $a\in k$,
put $\norm{a}=\frac{\mu(a\OO)}{\mu\OO}$ and $\ord
a=-\frac{\log\norm{a}}{\log q}$.
\item Let $\sch(k)$ be the Schwartz space of functions on $k$.
\end{itemize}
\subsection{Basic rules of calculation} 
To help reconcile the present system of notation to whatever
system the reader might be familiar with, we recall a few facts routinely
used below. 
\begin{itemize}
\item 
$\int f(x)d\mu(x)=\norm{a}\int f(ax)d\mu(x)$ for
$\mu$-integrable $f$
and $a\in k^\times$.
\item 
$\mu\OO=\norm{\kappa}^{1/2}$ and $\mu
\OO^\times=\frac{q-1}{q}\norm{\kappa}^{1/2}$.
\item $\norm{\cdot}$ is an absolute value of $k$ with respect
to which $k$ is complete.
\item $\ord$ is a normalized additive valuation of $k$.
\item 
$\int f(x)d\mu(x)=\frac{q-1}{q}\norm{\kappa}^{1/2}\int
f(t)\norm{t}d\mu^\times(t)$ for $\mu$-integrable $f$.
\end{itemize}

\subsection{Fourier transforms}\label{subsection:LocalFourierTransforms}
\subsubsection{}
Given a
complex-valued
$\mu$-integrable function
$f$ on $k$, put 
$$\FFF[f](\xi)=\hat{f}(\xi)=\int f(x)\ee(-x\xi)d\mu(x)$$
for all $\xi\in k$, thus defining
the {\em Fourier transform} $\hat{f}$, denoted also by
$\FFF[f]$, which is again a complex-valued function on $k$.  The Fourier
transform
$\hat{f}$ is continuous and tends to $0$ at infinity.
If $f$ is compactly supported,
then
$\hat{f}$ is locally constant. If $f$ is locally constant, then $\hat{f}$ is
compactly supported.  
For example, we
have 
\begin{equation}\label{equation:LocalFourierExample}
\FFF[\one_{x+b\OO}](\xi)=
\norm{\kappa}^{1/2}\norm{b}\one_{b^{-1}\kappa^{-1}\OO}(\xi)\ee(-x\xi)
\end{equation}
for all $b\in k^\times$ and $x,\xi\in k$.

\subsubsection{}
The Schwartz class $\sch(k)$ is stable
under Fourier transform and the {\em
Fourier inversion formula} states that
$$\Phi(x)=\int \ee(x\xi)\hat{\Phi}(\xi)d\mu(\xi)$$
for all $\Phi\in \sch(k)$. 
Our normalization of $\mu$ is chosen to
make the inversion formula hold in the latter particularly simple form;
in general there would be a positive constant depending on $\mu$ but
independent of
$\Phi$ multiplying the right side. 
The Fourier inversion formula implies the {\em squaring rule}
\begin{equation}\label{equation:LocalSquaring}
\FFF^2[\Phi](x)=\Phi(-x)
\end{equation}
for all $\Phi\in \sch(k)$.
We mention also the very frequently
used {\em scaling rule}
\begin{equation}\label{equation:LocalScaling}
\FFF[f^{(a)}]=\norm{a}\hat{f}^{(a^{-1})}\;\;\;(f^{(a)}(x)=
f(a^{-1}x))
\end{equation}
which holds for all $a\in k^\times$ and $\mu$-integrable $f$.
\subsubsection{}
We define a {\em lattice} $L\subset k$ to be a cocompact discrete
subgroup. We remark that lattices exist in $k$ only if $k$ is of positive characteristic.
We remark also that the notion of lattice figures 
only in Corollary~\ref{Corollary:LocalStirling} below, and otherwise
is unused in the sequel.
 For each lattice $L$ put $L^\perp=\{\xi\in k\vert \ee(x\xi)=1\;\mbox{for all}\;x\in L\}$. For all lattices $L$ again $L^\perp$ is a lattice and $(L^\perp)^\perp=L$.
The {\em Poisson summation formula} states that
\begin{equation}\label{equation:PoissonSummation}
\sum_{x\in L}\Phi(x)=
\mu(k/L)^{-1}\sum_{\xi\in L^\perp}\hat{\Phi}(\xi)\end{equation}
for all lattices $L\subset k$ and  Schwartz functions $\Phi\in \sch(k)$,
where $\mu(k/L)$ denotes the {\em covolume} of $L$ in $k$ with respect
to $\mu$.

\subsection{The linear functional $\MMM^{(a)}$ and linear operators $\LLL^{\pm}$}

\subsubsection{}\label{subsubsection:MDef}
For each $\Phi\in \sch(k)$ and $a\in k^\times$ put
\begin{equation}\label{equation:MDef}
\MMM^{(a)}[\Phi]=\int
(\Phi(t)-\one_{a\OO}(t)\Phi(0))d\mu^\times(t).
\end{equation} To see that $\MMM^{(a)}[\Phi]$ is a well-defined complex number, let
$G(t)$ temporarily denote the integrand on the right side and note the
following:
\begin{itemize}
\item $G(t)$ is defined and locally constant on
$k^\times$. 
\item $G(t)$ vanishes for $\max(\norm{t},\norm{t}^{-1})$ sufficiently
large. 
\end{itemize}
Therefore the integral on the right side of (\ref{equation:MDef})
converges.
We have a scaling rule
\begin{equation}\label{equation:MScale1}
\MMM^{(a)}[\Phi^{(b^{-1})}]=\MMM^{(ab)}[\Phi]=\MMM^{(a)}[\Phi]+\Phi(0)\ord
b
\end{equation}
for all $b\in k^\times$, cf.\ scaling rule
(\ref{equation:LocalScaling}) for the Fourier transform. 

To motivate the definition of $\MMM^{(a)}$ we remark that Lemma~\ref{Lemma:LFEConsequence} below can be reinterpreted as the assertion that
the linear functional $$\Phi\mapsto \Phi(0)\left(-\ord\kappa+
\frac{1}{q-1}\right)+\MMM^{(1)}[\Phi]$$
on $\sch(k)$ is the Fourier transform
of $\ord x$
in the sense of the theory of distributions.

\subsubsection{}\label{subsubsection:LDef}
For each $\Phi\in \sch(k)$ and $x\in k^\times$, put
\begin{equation}\label{equation:LDef}
\LLL^{\pm}[\Phi](x)=\int \frac{H(t^{\mp 1})
\Phi(xt)-\frac{1}{2}
\one_{\OO^\times}(t)\Phi(x)}{\norm{t-1}}\norm{t}d\mu^\times(t),
\end{equation}
where
$$H=\one_\OO-\frac{1}{2}\one_{\OO^\times}.$$
 To see that $\LLL^{\pm}[\Phi](x)$  is a
well-defined complex number, let $F^{\pm}(t)$ temporarily denote the integrand on the
right side of (\ref{equation:LDef}), and note the following:

\begin{itemize}
\item $F^{\pm}(t)$ is defined and locally constant on
$k^\times\setminus\{1\}$.
\item $F^{\pm}(t)=0$ for
$\norm{t-1}$ sufficiently small.
\item $F^{\pm}(t)=0$ for $\norm{t}$ sufficiently large.
\item $F^+(t)=0$ for $\norm{t}$ sufficiently small.
\item $F^-(t)=\norm{t}\Phi(0)$ for $\norm{t}$
sufficiently small.
\end{itemize}
Therefore the integral on the right side of
(\ref{equation:LDef}) converges. 
It is not difficult to verify that
\begin{equation}\label{equation:EulerMascheroni}
\lim_{\norm{x}\rightarrow 0} \LLL^-[\Phi](x)=
\Phi(0)\int_{\norm{t}<1}\norm{t}d\mu^\times(t)=\frac{\Phi(0)}{q-1}.
\end{equation}
We obtain  another 
integral representation 
\begin{equation}\label{equation:LDefBis}
\LLL^+[\Phi](x)=\int \frac{H(t)
\Phi(xt^{-1})-\frac{1}{2}
\one_{\OO^\times}(t)\Phi(x)}{\norm{1-t}}d\mu^\times(t)
\end{equation}
for $\LLL^+$
by substituting $t^{-1}$ for $t$ in the integral on the right
side of (\ref{equation:LDef}).  

In contrast to the case of the linear functional $\MMM^{(a)}$ we cannot easily give the motivation for the definition of the operators $\LLL^{\pm}$. 
The best we can say at present is that $\LLL^+$ makes the following theorem hold, and that $\LLL^-$ is indispensable to the proof of the theorem.
A more conceptual characterization of the operators $\LLL^{\pm}$ would be nice to have.

\begin{Theorem}[The local Stirling formula]\label{Theorem:LocalStirling}
There exists a unique linear operator
$$\KKK:\sch(k)\rightarrow\sch(k)$$
such that 
\begin{equation}\label{equation:KDef}
\KKK[\Phi](x)= 
\left\{\begin{array}{cl}
{
-\Phi(x)\left(\frac{1}{2}\ord\kappa+\ord x\right)+\LLL^+[\Phi](x)}
&\mbox{if $x\neq 0$,}\\\\
-\frac{1}{2}\Phi(0)\ord\kappa+\MMM^{(1)}[\Phi]&\mbox{if
$x=0$,}\end{array}\right.
\end{equation}
for all $\Phi\in \sch(k)$ and $x\in k$. Moreover, we have
\begin{equation}\label{equation:AntiCommutation}
\KKK[\FFF[\Phi]]=-\FFF[\KKK[\Phi]]
\end{equation}
for all $\Phi\in \sch(k)$.
\end{Theorem}
\noindent For the proof we need three lemmas.

\begin{Lemma}\label{Lemma:LFEConsequence}
We have
\begin{equation}\label{equation:LFEConsequence}
\begin{array}{cl}
&{
\int\Phi(x)
\left(\frac{1}{2}\ord \kappa+\ord x\right)d\mu(x)}\\= &{
-\hat{\Phi}(0)\left(\frac{1}{2}\ord\kappa+\ord
a\right) +\hat{\Phi}(0)/(q-1)+\MMM^{(a)}[\hat{\Phi}]}
\end{array}
\end{equation}
for all $\Phi\in\sch(k)$ and $a\in k^\times$.
\end{Lemma}
\proof By scaling rule (\ref{equation:MScale1}) the  right side of
(\ref{equation:LFEConsequence}) is independent of $a$.  For every $b\in k^\times$ we have
$$\hat{\Phi}(0)\ord
b=\int(\norm{b}^{-1}\Phi^{(b)}(x)-\Phi(x))
\left(\frac{1}{2}\ord \kappa+\ord x\right)d\mu(x),$$
and hence by the scaling rules
(\ref{equation:LocalScaling}) and 
(\ref{equation:MScale1}), the difference of
left and right sides of (\ref{equation:LFEConsequence}) remains unchanged
if we replace
$\Phi$ by
$\norm{b}^{-1}\Phi^{(b)}$ (and hence $\hat{\Phi}$ by
$\hat{\Phi}^{(b^{-1})}$).  
Clearly, the left side of (\ref{equation:LFEConsequence}) remains
unchanged if we replace $\Phi$ by $\Phi^{(u)}$ for any $u\in \OO^\times$.
The same holds for the right side by the scaling rules (\ref{equation:LocalScaling})
and (\ref{equation:MScale1}).
Consequently both sides of (\ref{equation:LFEConsequence}) remain
unchanged if we replace
$\Phi$ by  the averaged function
$\int_{\OO^\times}\Phi^{(u)}d\mu^\times(u)$, which is a finite linear
combination of functions of the form $\one_{b\OO}$ with $b\in k^\times$.
Taking into account all the preceding
reductions, we may now assume without loss of generality that
$a=1$ and
$\Phi=\one_{\OO}$. Then equation (\ref{equation:LFEConsequence})
is easy to check by direct calculation. We omit further details.
\qed

\begin{Lemma}\label{Lemma:LFEConsequence2}
We have
\begin{equation}\label{equation:LFEConsequence2}
\begin{array}{cl}
&{ 
\int\ee(-x\xi) \Phi(x)\left(\frac{1}{2}\ord\kappa+\ord x\right)
d\mu(x)}\\ =&{-\hat{\Phi}(\xi)\left(\frac{1}{2}\ord\kappa+\ord
\xi\right)+(\LLL^++\LLL^-)[\hat{\Phi}](\xi)}
\end{array}
\end{equation}
for all $\Phi\in\sch(k)$ and $\xi\in k^\times$.
\end{Lemma}
\proof We calculate as follows:\\
$$\begin{array}{cl}
&{ \hat{\Phi}(\xi)\left(\frac{1}{2}\ord\kappa+\ord \xi\right)+
\int\ee(-x\xi) \Phi(x)\left(\frac{1}{2}\ord\kappa+\ord
x\right) d\mu(x)}\\
=&{\hat{\Phi}(\xi)/(q-1)+\int
(\hat{\Phi}(\xi+t)-\hat{\Phi}(\xi)
\one_{\OO}(\xi^{-1}t))d\mu^\times(t)}\\=&{
\hat{\Phi}(\xi)/(q-1)+\int
(\hat{\Phi}(\xi(t+1))-\hat{\Phi}(\xi)\one_\OO(t))d\mu^\times(t)}\\
=&{\hat{\Phi}(\xi)/(q-1)+\int
\frac{\hat{\Phi}(\xi(t+1)))-\hat{\Phi}(\xi)
\one_{\OO}(t)}{\norm{t}}\norm{t}d\mu^\times(t)}\\=&{
\hat{\Phi}(\xi)/(q-1)+\int
\frac{\hat{\Phi}(\xi
t)-\hat{\Phi}(\xi)\one_\OO(t)}{\norm{t-1}}\norm{t}d\mu^\times(t)}\\
=&{
\int
\frac{\hat{\Phi}(\xi t)-\hat{\Phi}(\xi)\one_{\OO^{\times}}
(t)}{\norm{t-1}}
\norm{t}d\mu^\times(t)=(\LLL^++\LLL^-)[\hat{\Phi}](\xi)}\end{array}
$$\\
We get the first equality by applying (\ref{equation:LFEConsequence})
with $a$ replaced by $\xi$ and $\Phi(x)$ replaced by
$\ee(-x\xi)\Phi(x)$. The
rest of the calculation is routine.
\qed

\begin{Lemma}\label{Lemma:MollificationEstimate}
Fix $\Phi\in \sch(k)$. Fix $a\in k^\times$
such that $\Phi$ is supported in $a\OO$.
Fix $0\neq b\in a\OO$ such that $\Phi$ is constant on cosets of $b\OO$.
Fix $C>0$ such that
$|\Phi(x)|\leq C$ for all $x\in k$. Then we have
\begin{equation}\label{equation:MollificationEstimate}
\left|\frac{H(t)\Phi(xt^{-1})-\frac{1}{2}\one_{\OO^\times}(t)\Phi(x)}{\norm{1-t}}\right|\leq
C\norm{a/b}\one_{ta\OO}(x)\one_\OO(t)
\end{equation}
 for all $x\in k^\times$ and $t\in k^\times\setminus\{1\}$. 
  \end{Lemma}
\noindent For convenience in applying the lemma, note that
\begin{equation}\label{equation:IntegrationOfEstimate}
\int \one_{ta\OO}(x)\one_\OO(t)d\mu^\times(t)=\ord(x/a)\one_\OO(x/a)
\end{equation}
for all $a,x\in k^\times$.
\proof By hypothesis we have
$$\norm{1-t}<\norm{b/a}\Rightarrow H(t)\Phi(xt^{-1})-(1/2)\one_{\OO^\times}(t)\Phi(x)=0$$
 and hence 
$$\begin{array}{cl}
&\displaystyle\left|\frac{H(t)\Phi(xt^{-1})-(1/2)\one_{\OO^\times}(t)\Phi(x)}{\norm{1-t}}\right|
\\\\
\leq&C \norm{a/b}(H(t)\one_{ta\OO}(x)+(1/2)\one_{\OO^\times}(t)\one_{a\OO}(x)),
\end{array}
$$
whence the result.
\qed

\subsection{Proof of the theorem}
The right side of (\ref{equation:KDef}) 
defines on $k$ a complex-valued function $\KKK[\Phi]$ 
depending linearly  on $\Phi$.  By
Lemma~\ref{Lemma:MollificationEstimate} the function $\LLL^+[\Phi]$ is $\mu$-integrable
 and compactly supported,
and hence so is the function $\KKK[\Phi]$. 
Similarly, $\KKK[\FFF[\Phi]]$ is
$\mu$-integrable and compactly supported.
Since $\KKK[\Phi]$ is $\mu$-integrable and compactly supported,
it follows that the Fourier transform $\FFF[\KKK[\Phi]]$ is well-defined and locally constant. 
We have
\begin{equation}\label{equation:MollificationEstimateBis}
\FFF[\LLL^+[\Phi]](\xi)=\LLL^-[\FFF[\Phi]](\xi)\;\;\;(\xi\in k^\times),
\end{equation}
since by another application of
Lemma~\ref{Lemma:MollificationEstimate} it is justified to
reverse the order of integration, and by limit formula (\ref{equation:EulerMascheroni})
it follows that
\begin{equation}\label{equation:MollificationEstimateTer}
\FFF[\LLL^+[\Phi]](0)=\frac{\hat{\Phi}(0)}{q-1}.
\end{equation}
 By (\ref{equation:MollificationEstimateTer}) above,  Lemma~\ref{Lemma:LFEConsequence}, and the definitions we have
$$\FFF[\KKK[\Phi]](0)=-\KKK[\FFF[\Phi]](0).$$
By (\ref{equation:MollificationEstimateBis}) above,
Lemma~\ref{Lemma:LFEConsequence2}, and the definitions we have
$$\FFF[\KKK[\Phi]](\xi)=-\KKK[\FFF[\Phi]](\xi)\;\;(\xi\in k^\times).$$
Therefore we have $\FFF[\KKK[\Phi]]=-\KKK[\FFF[\Phi]]$,
hence $\KKK[\FFF[\Phi]]$ is both compactly supported and locally constant,
and hence $\KKK[\FFF[\Phi]]\in \sch(k)$. 
Since $\Phi\in \sch(k)$ was arbitrarily specified and the operator $\FFF$ is invertible,
it follows that the operator $\KKK$ stabilizes 
the  Schwartz space $\sch(k)$ and anticommutes with the Fourier transform $\FFF$, as claimed.
\qed

\bigskip
The following corollary will not be needed in the sequel---we include it just to motivate the theorem.
\begin{Corollary}
\label{Corollary:LocalStirling}
Fix a  Schwartz function $\Phi\in\sch(k)$ and a lattice $L\subset k$.
Put
$$\vartheta(t):=\sum_{x\in L}\Phi(t^{-1}x),\;\;\;
\vartheta^*(t):=\sum_{\xi\in L^\perp}\hat{\Phi}(t^{-1}\xi)\;\;\;(t\in k^\times).$$ 
Then the following hold: 
\begin{equation}\label{equation:LocalStirling1}
\vartheta(t)=\mu(k/L)^{-1}\norm{t}\vartheta^*(t^{-1})
\end{equation}
\begin{equation}\label{equation:LocalStirling2}
\vartheta(t)=\Phi(0)\;\mbox{for $\norm{t}$ sufficiently small.}
\end{equation}
\begin{equation}\label{equation:LocalStirling3}
\vartheta(t)=\mu(k/L)^{-1}\hat{\Phi}(0)\norm{t}\;\mbox{for $\norm{t}$ sufficiently large.}
\end{equation}
\begin{equation}\label{equation:LocalStirling4}
\mbox{$\vartheta(t)$ is a locally constant function of $t$.}
\end{equation}
We have
\begin{equation}\label{equation:LocalStirling}
\begin{array}{cl}
&\;\;\;\;\vartheta(a)\ord \kappa\\\\
&\displaystyle+\sum_{0\neq x\in L}\Phi(a^{-1}x)\ord x+
\mu(k/L)^{-1}\sum_{0\neq \xi\in L^\perp}\norm{a}\hat{\Phi}(a\xi)\ord \xi\\\\
&\displaystyle-\left(-\Phi(0)\ord a+\int(\Phi(t)-\Phi(0)\one_\OO(t))d\mu^\times(t))\right)\\\\
&\displaystyle -\mu(k/L)^{-1}\norm{a}
\left(\hat{\Phi}(0)\ord a+\int
\left(\hat{\Phi}(t)-\hat{\Phi}(0)\one_{\OO}(t)\right)d\mu^\times(t)\right)\\\\
=&\displaystyle\int\frac{\vartheta(at)-\left\{\begin{array}{cl}
\Phi(0)&\mbox{if $\norm{t}<1$}\\
\vartheta(a)&\mbox{if $\norm{t}=1$}\\
\mu(k/L)^{-1}\hat{\Phi}(0)\norm{at}&\mbox{if $\norm{t}>1$}
\end{array}\right.}{\norm{1-t}}d\mu^\times(t)
\end{array}
\end{equation}
for all $a\in k^\times$.
\end{Corollary}
\proof[Proof of the corollary]
Functional equation (\ref{equation:LocalStirling1}) 
follows from  scaling formula (\ref{equation:LocalScaling}) 
and the Poisson summation formula (\ref{equation:PoissonSummation}). 
Statements (\ref{equation:LocalStirling2},\ref{equation:LocalStirling3},\ref{equation:LocalStirling4})
are easy to verify. We omit the details. Statements (\ref{equation:LocalStirling2},\ref{equation:LocalStirling3},\ref{equation:LocalStirling4}) granted, it is clear that the integral
on the right side of (\ref{equation:LocalStirling}) converges.
 We turn to the proof of equation (\ref{equation:LocalStirling}).
First, we reduce to the case $a=1$ by observing that  under replacement
of the pair $(\Phi,a)$ by the pair $(\Phi^{(a)},1)$, neither the left side nor the right side
of (\ref{equation:LocalStirling}) change. 
Next, we note that
$$\sum_{x\in L}\KKK[\Phi](x)+\mu(k/L)^{-1}\sum_{\xi\in L^\perp}\KKK[\hat{\Phi}](\xi)=0$$
by the Poisson summation formula (\ref{equation:PoissonSummation})
and the anticommutation relation (\ref{equation:AntiCommutation}). 
Clearly, it is possible to rearrange the terms above
to put the left side in coincidence with the left side of (\ref{equation:LocalStirling})
in the case $a=1$, leaving a sum $A+\mu(k/L)^{-1}A^*$ on the right side, where
$$A=\sum_{0\neq x\in L}\LLL^+[\Phi](x),\;\;\;
A^*=\sum_{0\neq \xi\in L^\perp}
\LLL^+[\hat{\Phi}](\xi).$$
Using the presentation (\ref{equation:LDefBis}) of the operator $\LLL^+$
and carrying the sum under the integral, we have
$$A=\int \frac{H(t)
(\vartheta(t)-\Phi(0))-\frac{1}{2}
\one_{\OO^\times}(t)(\vartheta(1)-\Phi(0))}{\norm{1-t}}d\mu^\times(t).
$$
Exchange of sum and integral is justified by
Lemma~\ref{Lemma:MollificationEstimate}
and the remark (\ref{equation:IntegrationOfEstimate})
immediately following.
We have a similar representation for $A^*$, which after
the substitution of $t^{-1}$ for $t$
takes the form  $$A^*=\int \frac{H(t^{-1})
(\vartheta^*(t^{-1})-\hat{\Phi}(0))-\frac{1}{2}
\one_{\OO^\times}(t)(\vartheta^*(1)-\hat{\Phi}(0))}{\norm{t-1}}\norm{t}d\mu^\times(t).
$$
Finally, by exploiting functional equation (\ref{equation:LocalStirling1})
we can bring  $A+\mu(k/L)^{-1}A^*$ to the form of the right
side of (\ref{equation:LocalStirling}) in the case $a=1$. We omit the remaining details of bookkeeping.  \qed

\subsection{Remark}
Because the sum
$$
\sum_{0\neq \xi\in L^\perp}\norm{a}\hat{\Phi}(a\xi)\log\norm{\xi}
$$
vanishes for $\norm{a}\gg 1$, equation (\ref{equation:LocalStirling}) provides a precise asymptotic description of the sum
\begin{equation}\label{equation:StirlingDiscussionTwo}
\sum_{0\neq x\in L}\Phi(a^{-1}x)\log\norm{x}
\end{equation}
as $\norm{a}\rightarrow\infty$.
In the simple special case $\Phi=\one_\OO$, the latter sum is analogous to the sum
\begin{equation}\label{equation:StirlingDiscussionThree}
\begin{array}{cl}
&\displaystyle\frac{1}{2}\sum_{0\neq x\in \ZZ}\one_{[-1,1]}(n^{-1}x)\log |x|\\\\=&\log n!=
\displaystyle
n\log n -n+\frac{1}{2}\log n+\frac{1}{2}\log(2\pi)+o_{n\rightarrow\infty}(1)
\end{array}
\end{equation}
on the real line.
So it is reasonable to view equation (\ref{equation:LocalStirling}) as an analogue and generalization of Stirling's
formula in the nonarchimedean setting. 
Given the key role played by Theorem~\ref{Theorem:LocalStirling} 
in deriving (\ref{equation:LocalStirling}),
we choose to regard Theorem~\ref{Theorem:LocalStirling} itself as an
analogue and generalization of the classical Stirling formula. This explains our terminology.

But the analogy of Theorem~\ref{Theorem:LocalStirling}
with Stirling's formula is rather imperfect. A closer look reveals some complications.
As it happens, for Theorem~\ref{Theorem:LocalStirling}
and Corollary~\ref{Corollary:LocalStirling} we can derive  direct analogues on the real line, using standard methods of the theory of tempered distributions, and these results in turn lead to a precise description of the asymptotic
behavior 
as $n\rightarrow\infty$ of sums of the form
$$
\sum_{0\neq x\in \ZZ}\varphi(n^{-1}x)\log|x|,
$$
where $\varphi(x)$ is a Schwartz function on the real line. For example, we can derive in this way
the asymptotic formula
\begin{equation}\label{equation:StirlingDiscussionFive}
\begin{array}{cl}
&\displaystyle\sum_{0\neq x\in \ZZ}\exp(-\pi(x/n)^2)\log |x|\\\\
=&n\log n-(\log 2+\frac{1}{2}\log \pi+\gamma/2)n+\log 2\pi+o_{n\rightarrow\infty}(1),
\end{array}
\end{equation}
where $\gamma$ is the Euler-Mascheroni constant.
From the point of view of harmonic analysis, the sum (\ref{equation:StirlingDiscussionTwo}) in the simple special case $\Phi=\one_\OO$ is more closely analogous to the ``soft-edged'' sum (\ref{equation:StirlingDiscussionFive}) than to the ``hard-edged'' sum (\ref{equation:StirlingDiscussionThree}). These complications taken into account, we may only say that a ``soft'' analogy  exists between Stirling's formula
and Theorem~\ref{Theorem:LocalStirling}. We shall discuss analogues of Theorem~\ref{Theorem:LocalStirling}  for $\RR$ and $\CC$
and their applications in detail on another occasion.

\section{The adelic Stirling formula}\label{section:GlobalStirling}
From now on in this paper we work in the setting of
harmonic analysis on the adele ring of a global field of positive
characteristic, a setting for which the (unfortunately) rather elaborate
notation is set out in detail in tables immediately below.  The
main result of this section is the {\em adelic Stirling formula}
(Theorem~\ref{Theorem:AdelicStirling}),
which in form resembles Corollary~\ref{Corollary:LocalStirling}.
See \cite{Ramakrishnan} for background on global harmonic analysis.

\subsection{Data}
\begin{itemize}
\item Let $k$ be a global field of positive characteristic,
of genus $g$, and with constant field $\FF_q$.
\item Fix a nonzero K\"{a}hler differential
$\omega\in\Omega=\Omega_{k/\FF_q}$.
\item Fix a nonconstant character $\lambda:\FF_q\rightarrow U(1)$,
as in \S\ref{section:Duality}.
\end{itemize}
All constructions below proceed naturally from these
choices. 

\subsection{Notation (local)}
Let $v$ be any place of $k$. 
\begin{itemize}
\item Let $k_v$ be the completion of $k$ at $v$.
\item Let $\OO_v$ be the maximal compact subring of $k_v$.
\item Let $\FF_v$ be the residue field of $\OO_v$.
\item Let $\Res_v:\Omega\otimes_k k_v\rightarrow\FF_v$ be the residue
map at $v$.
\item Put $\ee_v=(x\mapsto\lambda(\tr_{\FF_v/\FF_q}\Res_v \omega
x)):k_v\rightarrow U(1)$. 
\item Let $\mu_v$ be Haar measure on $k_v$.
\item Normalize $\mu_v$ by the rule $\mu_v\OO_v\cdot
\mu_v(\kappa_v^{-1}\OO_v)=1$. 
\item Let $q_v$ be the cardinality of $\FF_v$.
\item For each $a\in k_v$, put
$\norm{a}_v=\frac{\mu_v(a\OO_v)}{\mu_v\OO_v}$ and
$\ord_v a=-\frac{\log \norm{a}_v}{\log q_v}$.
\item Let $\mu_v^\times$ be Haar measure on $k_v^\times$.
\item Normalize $\mu_v^\times$
by the rule $\mu_v^\times \OO_v^\times=1$.
\item Let $\ord_v\omega$ be the order of vanishing of
$\omega$ at $v$.
\item Choose $\kappa_v\in k_v^\times$ such that $\ord_v\kappa_v=\ord_v\omega$,
and hence
$\kappa_v^{-1}\OO_v=\{\xi\in k_v\vert
\xi\OO_v\subset \ker\ee_v\}$.
\item Let $\sch(k_v)$ be the Schwartz class of functions on $k_v$.
\end{itemize}

\subsection{Notation (global)}\label{subsection:GlobalNotation}
\begin{itemize}
\item Let $\adeles$ be the adele ring of
$k$. 
\item Let $\OO\subset \adeles$ be the maximal compact subring.
\item Recall that each adele $x\in\adeles$ is a family $x=[x_v]$
indexed by places $v$ such that $x_v\in k_v$
for all
$v$ and $x_v\in \OO_v$ for all but finitely many $v$.
Moreover, we have $\OO=\prod_v \OO_v$.
\item As usual, we regard $k$ as diagonally embedded in $\adeles$;
thus $k$ becomes a cocompact discrete subgroup of $\adeles$.
\item Let $\mu$ be a Haar measure on $\adeles$.
\item Normalize $\mu$ by the rule
$\mu\OO=q^{1-g}$. Equivalently: $\mu(\adeles/k)=1$.
\item Recall that each idele $a\in\adeles^\times$ is
a family $a=[a_v]$ indexed by places $v$ such that
$a_v\in k_v^\times$ for all $v$ and $a_v\in \OO_v^\times$
for all but finitely many $v$. Moreover, we have
$\OO^\times=\prod_v
\OO_v^\times$.
\item For each $a\in\adeles^\times$, put $\norm{a}=\frac{\mu
(a\OO)}{\mu\OO}$. 
\item Let $\mu^\times$ be a Haar measure on $\adeles^\times$.
\item Normalize $\mu^\times$ by the rule $\mu^\times \OO^\times=1$.
\item Put $R_\omega=\left([x_v]\mapsto \sum_v 
\tr_{\FF_v/\FF_q}\Res_v x_v\omega\right):\adeles\rightarrow \FF_q$.
\item For all $x,y\in \adeles$, put $\langle
x,y\rangle=R_\omega(xy)\in \FF_q$.
\item Put $\ee=\lambda\circ R_\omega=([x_v]\mapsto \prod_v \ee_v(x_v)):
\adeles\rightarrow U(1)$.
\item Put $\kappa=[\kappa_v]\in \adeles^\times$.
\item We write  $\ord_v\omega=\ord_v\kappa=\ord_v\kappa_v$ for all places $v$.
\item For each place $v$ let 
$$i_v:k_v\rightarrow \adeles,\;\;\;i_v^\times:k_v^\times
\rightarrow \adeles^\times$$
be the unique maps such that 
$$(i_v(x))_w=\left\{\begin{array}{rl}
x&\mbox{if $w=v$,}\\
0&\mbox{otherwise,}
\end{array}\right.\;\;(i_v^\times(t))_w=\left\{\begin{array}{rl}
t&\mbox{if $w=v$,}\\
1&\mbox{otherwise,}
\end{array}\right.$$ for all $x\in k_v$,
$t\in k_v^\times$ and places $w$.
\item Let $\sch(\adeles)$ be the Schwartz class of functions on
$\adeles$.
\end{itemize}
\subsection{Basic rules of calculation}\label{subsection:BasicAdelicRules}
To help reconcile the present system of notation to whatever
system the reader might be familiar with, we recall a few facts routinely
used below. 
\begin{itemize}
\item $\int f(x)d\mu(x)=\norm{a}\int f(ax)d\mu(x)$ 
all $\mu$-integrable $f$
and $a\in \adeles^\times$.
\item $\mu=\bigotimes \mu_v$ and $\mu^\times=\bigotimes \mu_v^\times$.
\item $\norm{a}=\prod_v \norm{a_v}_v$ for all $a=[a_v]\in
\adeles^\times$.
\item $\norm{x}=1$ for all $x\in k^\times$. (Artin's product
formula.)
\item $R_\omega x=0$ for all $x\in k$. (Sum of residues equals
zero.)
\item $\norm{\kappa}^{1/2}=q^{1-g}=\mu\OO$.
\item $\kappa^{-1}\OO=\{\xi\in \adeles\vert \xi\OO\subset \ker
R_\omega\}=\{\xi\in \adeles\vert \xi\OO\subset \ker\ee\}$.
\item For all $a,b\in \adeles^\times$ such that $a\OO\subset b\OO$,
the spaces $b\OO/a\OO$ and $a^{-1}\kappa^{-1}\OO/b^{-1}\kappa^{-1}\OO$
are paired perfectly by $\langle \cdot,\cdot\rangle$. Moreover we
have
$k^\perp=k$ with respect to $\langle\cdot,\cdot \rangle$. (Serre
duality.)
\end{itemize}

\subsection{The adelic Fourier transform and the theta symbol}
\subsubsection{}
Given a $\mu$-integrable complex-valued
function
$f$ on
$\adeles$, put
$$\FFF[f](\xi)=\hat{f}(\xi)=
\int f(x)\lambda(-\langle x,\xi\rangle)d\mu(x)=\int
f(x)\ee(-x\xi)d\mu(x)$$ for all $\xi\in \adeles$, thereby defining 
the {\em Fourier transform} $\hat{f}$, also denoted by $\FFF[f]$, which
is a complex-valued continuous function on $\adeles$ tending to $0$ at
infinity. If $f$ is compactly supported, then $\hat{f}$ is locally
constant. If $f$ is locally constant, then $\hat{f}$ is compactly
supported.  For example, we have
\begin{equation}\label{equation:GlobalFourierExample}
\FFF[\one_{x+b\OO}](\xi)=
q^{1-g}\norm{b}\one_{b^{-1}\kappa^{-1}\OO}(\xi)
\lambda(-\langle x,\xi\rangle)
\end{equation}
for all $b\in \adeles^\times$ and $x,\xi\in \adeles$.
\subsubsection{}
The Schwartz class $\sch(\adeles)$ is stable
under the Fourier transform
and the {\em Fourier inversion formula} states that
\begin{equation}\label{equation:AdelicFourierInversion}
\Phi(x)=\int \hat{\Phi}(\xi)\ee(x\xi)d\mu(\xi)
\end{equation}
for all $\Phi\in \sch(\adeles)$. Our normalization of $\mu$
is chosen to make the Fourier inversion formula hold in the
simple form stated above; in general there would be a positive constant
depending on $\mu$ but independent of $\Phi$ multiplying the right side.
The Fourier inversion formula implies the {\em squaring rule}
\begin{equation}\label{equation:AdelicFourierSquared}
\FFF^2[\Phi](x)=\Phi(-x)
\end{equation}
for all $\Phi\in \sch(\adeles)$. We mention also
the {\em scaling rule}
\begin{equation}\label{equation:AdelicScaling}
\FFF[f^{(a)}]=\norm{a}\hat{f}^{(a^{-1})}\;\;\;
(f^{(a)}(x)=f(a^{-1}x))
\end{equation}
holding for all
$a\in \adeles^\times$ and $\mu$-integrable $f$.

\subsubsection{}
The {\em Poisson summation formula} 
states that
$$\sum_{x\in
k}\Phi(a^{-1}x)=\norm{a}\sum_{\xi\in
k}\hat{\Phi}(a\xi)$$ for all
$a\in \adeles^\times$ and $\Phi\in \sch(\adeles)$. 
Note that since both $\Phi$ and $\hat{\Phi}$ are compactly supported,
only
finitely many nonzero terms occur in the sums on either side of the
formula.  The normalization $\mu\OO=q^{1-g}=\mu(\adeles/k)$ of Haar
measure
$\mu$ on $\adeles$ was chosen to make the Poisson summation formula
hold in the particularly simple form above; in general the right side
would be multiplied by a positive constant depending on $\mu$
but independent of $\Phi$.

\subsubsection{}\label{subsubsection:ThetaSymbolDef}
For all $a\in \adeles^\times$ and $\Phi\in \sch(\adeles)$ 
put
$$\Theta(a,\Phi)=\sum_{x\in k}\Phi(a^{-1}x),$$
thereby defining the {\em theta symbol} $\Theta(\cdot,\cdot)$.
We have
\begin{equation}\label{equation:ThetaIdeleClassDependence}
\Theta(ax,\Phi)=\Theta(a,\Phi)
\end{equation}
for all $x\in k^\times$. In other words, $\Theta(a,\Phi)$ depends
only on the image of $a$ in the {\em idele class group}
$\adeles^\times/k^\times$. Clearly, we have a {\em scaling rule}
\begin{equation}\label{equation:ThetaScaling}
\Theta(a,\Phi^{(b)})=\Theta(ab,\Phi)=\Theta(1,\Phi^{(ab)})
\end{equation}
for all $b\in \adeles^\times$.  We have a {\em functional equation}
\begin{equation}\label{equation:ThetaPoisson}
\Theta(a,\Phi)=\norm{a}\Theta(a^{-1},\hat{\Phi}),
\end{equation}
which is just a rewrite of the Poisson summation formula.

\begin{Proposition}\label{Proposition:ThetaAsymptotics}
Fix $\Phi\in \sch(\adeles)$. Fix $b\in \adeles^\times$
and $f\in \adeles^\times \cap \OO$ such that
$\Phi$ is constant on cosets of $b\OO$ and supported in $f^{-1}b\OO$.
Then: (i) We have
$$\Theta(a,\Phi)=
\left\{\begin{array}{rl}
\norm{a}\hat{\Phi}(0)&\mbox{if $\norm{ab}>q^{2g-2}$,}\\\\
\Phi(0)&\mbox{if $\norm{ab}<\norm{f}$,}
\end{array}\right.
$$
for all $a\in \adeles^\times$. (ii) Moreover, $\Theta(a,\Phi)$ depends in locally constant fashion on $a$.
\end{Proposition}
\proof From example (\ref{equation:GlobalFourierExample}) and Serre duality it follows that
the Fourier transform $\hat{\Phi}$ is constant on cosets of
$fb^{-1}\kappa^{-1}\OO$ and supported in  $b^{-1}\kappa^{-1}\OO$.
We therefore have
\begin{equation}\label{equation:ThetaAsymptotics}
\sum_{x\in k\cap af^{-1}b\OO}\Phi(a^{-1}x)=
\Theta(a,\Phi)=
\norm{a}\sum_{\xi\in k\cap
a^{-1}b^{-1}\kappa^{-1}\OO}\hat{\Phi}(a\xi).
\end{equation}
via (\ref{equation:ThetaPoisson}). Statement (i) follows via the Artin product formula.
Statement (ii) follows from the observation that
the map $a\mapsto k\cap af^{-1}b\OO$
from $\adeles^\times$ to finite subsets of $k$ is locally constant. \qed
\begin{Theorem}[The adelic Stirling formula]
\label{Theorem:AdelicStirling}
For all ideles $a\in \adeles^\times$, Schwartz functions $\Phi\in
\sch(\adeles)$ and places
$v$ of $k$ we have
\begin{equation}\label{equation:AdelicStirling}
\begin{array}{cl}
&\displaystyle \Theta(a,\Phi)\ord_v \omega+
\sum_{x\in k^\times} \Phi(a^{-1}x)\ord_v x+
\sum_{\xi\in k^\times}\norm{a}\hat{\Phi}(a\xi)\ord_v \xi\\&\\
=&{\displaystyle
\int\frac{\Theta(i_v^\times(t)a,\Phi)-\left\{\begin{array}{cl}
\Phi(0)&\mbox{if $\norm{t}_v<1$}\\
\Theta(a,\Phi)&\mbox{if $\norm{t}_v=1$}\\
\norm{a}\norm{t}_v\hat{\Phi}(0)&\mbox{if $\norm{t}_v>1$}
\end{array}\right.}{\norm{1-t}_v}d\mu_v^\times(t)}\\&\\
&{\displaystyle -\Phi(0)\ord_v
a_v+\int\left(\Phi(i_v(t))-\Phi(0)\one_{\OO_v}(t)\right)d\mu_v^\times(t)}\\&\\
&{\displaystyle
+\norm{a}\left(\hat{\Phi}(0)\ord_v a_v+\int(
\hat{\Phi}(i_v(t))-\hat{\Phi}(0)
\one_{\OO_v}(t))d\mu_v^\times(t)\right).}
\end{array}
\end{equation}
\end{Theorem}
\noindent 
Note that Proposition~\ref{Proposition:ThetaAsymptotics} guarantees convergence of the integral  involving the theta symbol on the right side of (\ref{equation:AdelicStirling}).
The pattern set by the
proof of  Corollary~\ref{Corollary:LocalStirling} 
will guide us in proving (\ref{equation:AdelicStirling}). The proof will be completed in
\S\ref{equation:EndOfAdelicStirlingProof} after some preparation.  
 The rationale for the Stirling formula terminology here is essentially
the same as that offered in connection
with Theorem~\ref{Theorem:LocalStirling} and Corollary~\ref{Corollary:LocalStirling}.
Namely, to the extent that sums of the form $\sum_{x\in k^\times} \Phi(a^{-1}x)\ord_v x$
may be regarded as analogues of $\log n!$, formula (\ref{equation:AdelicStirling})
may be regarded as an analogue of the classical Stirling formula. 
It is an interesting problem
to devise and interpret a version of the adelic Stirling formula for number fields.
We hope to have progress to report on this problem in the near future.

\subsection{Recollection of the tensor
decomposition of
$\sch(\adeles)$}
\subsubsection{} Given any family of vector spaces $\{V_i\}_{i\in I}$,
where each $V_i$ is equipped with a {\em neutral element}
$0\neq e_i\in V_i$, there is a natural way to form the (possibly infinite)
tensor product
$\bigotimes_{i\in I}V_i$. Namely, the latter is by definition spanned by
symbols of the form 
$$\bigotimes_{i\in I} v_i\;\;\;
(
\mbox{$v_i\in V_i$ for all $i$, and $v_i=e_i$ for all but finitely
many
$i$}),$$
 subject to the obvious
relations. Moreover, given for each index $i$ a linear
endomorphism $L_i$ of $V_i$ such that $L_ie_i=e_i$ for all but finitely
many
$i$, there is a natural tensor product $\bigotimes_{i\in I}L_i$ of
operators, namely, that sending each symbol $\bigotimes_{i\in
I} v_i$ to $\bigotimes_{i\in I} (L_iv_i)$. 

\subsubsection{}
In the general sense above, as is well known,
$\sch(\adeles)$ can naturally be identified with $\bigotimes_v
\sch(k_v)$, where for each place $v$ the space $\sch(k_v)$ is
equipped with the neutral element $\one_{\OO_v}$. 
(See, for example, \cite[p.\ 260]{Ramakrishnan}. But note that 
in the definition presented there the condition ``$f_v\vert_{\OO_v}=1$''
should be strengthened to ``$f_v=\one_{\OO_v}$''. 
Otherwise the function $f$ thus defined might fail to be compactly supported.)

\subsubsection{}
We have at
least one natural example of a linear operator on $\sch(\adeles)$
that factors as an infinite tensor product, namely the 
Fourier transform $\FFF$. More precisely, for each place
$v$, let
$$\FFF_v:\sch(k_v)\rightarrow\sch(k_v)$$
be the local Fourier transform defined by the rule
$$\FFF_v[\Phi](\xi)=\int \Phi(x)\ee_v(-x\xi)d\mu_v(x).$$
Since $\FFF_v[\one_{\OO_v}]=\one_{\OO_v}$ for all but finitely many $v$,
it makes sense to form the tensor product of the operators $\FFF_v$ over
all places $v$. One can easily verify that this tensor product
does indeed coincide with $\FFF$.
\subsubsection{}
Let a place $v$ and a linear endomorphism
of $L$ of $\sch(k_v)$ be given. 
Let $\adeles^v$ be
the ``coordinate hyperplane'' of $\adeles$ consisting of families
$x=[x_w]$ such that
$x_v=0$. Given $y\in \adeles^v$
and $\Phi\in \sch(\adeles)$, let $\Phi_{y,v}\in
\sch(k_v)$ be defined by the rule 
$$\Phi_{y,v}(x)=\Phi(i_v(x)+y)$$
for all $x\in k_v$. Then, so we claim, there exists a unique linear
endomorphism
$\tilde{L}$ of $\sch(\adeles)$ such that
$$\tilde{L}[\Phi](i_v(x)+y)=L[\Phi_{y,v}](x)$$
for all $\Phi\in \sch(\adeles)$, $x\in k_v$ and $y\in \adeles^v$.
To prove the claim one just has to check that $\tilde{L}$ is the
tensor product of the family 
$$L_w=\left\{\begin{array}{rl}
L&\mbox{if $v=w$,}\\
1&\mbox{if $v\neq w$,}
\end{array}\right.
$$
extended over all places $w$. We omit these details. We call
$\tilde{L}$ the {\em canonical prolongation} of $L$.

\subsection{The operators
$\MMM_v$,  $\LLL_v$ and $\KKK_v$}
Fix a place $v$ of $k$. We keep the notation introduced in the preceding
discussion of the tensor decomposition of $\sch(\adeles)$.

\subsubsection{}
Given $\Phi\in \sch(\adeles)$, we define a function $\MMM_v[\Phi]$ on
$\adeles^v$ by the rule
$$\MMM_v[\Phi](y)=
\int(\Phi_{y,v}(t)-\one_{\OO_v}(t)\Phi_{y,v}(0))d\mu_v^\times(t)$$
for all $y\in \adeles^v$. The function
$\MMM_v[\Phi]$ is well-defined by a repetition of the argument
given in
\S\ref{subsubsection:MDef}. Note in particular that
\begin{equation}\label{equation:MSubVAlternativeDef}
\MMM_v[\Phi](0)=\int
(\Phi(i_v(t))-\one_{\OO_v}(t)\Phi(0))d\mu_v^\times(t)
\end{equation}
for all $\Phi\in \sch(\adeles)$.

\subsubsection{}
 Given $\Phi\in \sch(\adeles)$, we define a function $\LLL_v[\Phi]$
on $\adeles\setminus \adeles^v$ by the rule 
$$\LLL_v[\Phi](i_v(x)+y)=\displaystyle\int\frac{H_v(t)\Phi_{y,v}(t^{-1}x)-\frac{1}{2}
\one_{\OO_v^\times}(t)\Phi_{y,v}(x)} {\norm{1-t}_v}d\mu_v^\times(t)
$$
for all $x\in k_v^\times$ and
$y\in\adeles^v$, where
$$H_v=\one_{\OO_v}-\frac{1}{2}\one_{\OO_v^\times}.$$
 The function $\LLL_v[\Phi]$
is well-defined by a repetition of the argument given in
\S\ref{subsubsection:LDef}. Equivalently, we have
\begin{equation}\label{equation:LSubVAlternativeDef}
\LLL_v[\Phi](z)=\int\frac{H_v(t)\Phi(i_v^\times(t^{-1})z)-\frac{1}{2}\one_{\OO_v^\times}(t)\Phi(z)}
{\norm{1-t}_v}d\mu_v^\times(t)
\end{equation}
for all $z\in \adeles\setminus \adeles^v$. Note that the operator $\LLL_v$
preserves supports in the sense that for all
$a\in
\adeles^\times$ and $\Phi\in \sch(\adeles)$,
if $\Phi$ vanishes outside $a\OO$, 
then $\LLL_v[\Phi]$  vanishes outside  $a\OO\setminus \adeles^v$.

\subsubsection{}
Theorem~\ref{Theorem:LocalStirling} and the canonical
prolongation process now yield a unique linear operator
$$\KKK_v:\sch(\adeles)\rightarrow\sch(\adeles)$$
 such that
$$\begin{array}{cl}
&\displaystyle\KKK_v[\Phi](i_v(x)+y)\\\\
=&
\displaystyle\left\{\begin{array}{ll}
-\Phi(i_v(x)+y)(\frac{1}{2}\ord_v\omega+\ord_vx)+\LLL_v[\Phi](i_v(x)+y)
&\mbox{if $x\neq 0$,}\\ -\frac{1}{2}\Phi(y)\ord_v\omega+
\MMM_v[\Phi](y)&\mbox{if $x=0$,}
\end{array}\right.
\end{array}
$$
for all $\Phi\in \sch(\adeles)$, $x\in k_v$ and $y\in \adeles^v$.
Moreover, since $\FFF$ factors as the tensor product of local
Fourier transforms, and the operator on $\sch(k_v)$
defined by Theorem~\ref{Theorem:LocalStirling} anticommutes with
$\FFF_v$, we necessarily have
$$\KKK_v[\FFF[\Phi]]=-\FFF[\KKK_v[\Phi]]$$
for all $\Phi\in \sch(\adeles)$, i.~e., $\KKK_v$ anticommutes
with the adelic Fourier transform $\FFF$.

\begin{Lemma}\label{Lemma:IntegralAndSums}
Fix $\Phi\in \sch(\adeles)$. Fix a place $v$ of $k$. Put
$$
\vartheta(t)=\Theta(i_v^\times(t),\Phi)
$$
for all $t\in k_v^\times$. Then we have
\begin{equation}\label{equation:KeyIdentity}
\begin{array}{cl}
&\displaystyle\sum_{x\in
k^\times}\LLL_v[\Phi](x)+
\sum_{\xi\in k^\times}\LLL_v[\hat{\Phi}](\xi)\\\\
=&\displaystyle
\int\frac{\vartheta(t)-\left\{\begin{array}{cl}
\Phi(0)&\mbox{if $\norm{t}_v<1$}\\
\vartheta(1)&\mbox{if $\norm{t}_v=1$}\\
\norm{t}_v\hat{\Phi}(0)&\mbox{if $\norm{t}_v>1$}
\end{array}\right.}{\norm{1-t}_v}d\mu_v^\times(t).
\end{array}
\end{equation}
\end{Lemma}
\noindent Convergence of the sums on the left and the integral on the right
will be established in the course of the proof. 
\proof The integrand on the right side of
(\ref{equation:KeyIdentity}) is 
locally constant  and is supported on a compact subset of $k^\times_v
\setminus\{1\}$ by
Proposition~\ref{Proposition:ThetaAsymptotics}. So the
integral on the right side of (\ref{equation:KeyIdentity})
converges. Now fix
$a\in
\adeles^\times$ such that both
$\Phi$ and
$\hat{\Phi}$ are supported in $a\OO$,  hence
both $\LLL_v[\Phi]$ and $\LLL_v[\hat{\Phi}]$ are supported in
$a\OO\setminus \adeles^v$,  and hence the quantities
$$\Phi(i_v(t)^{-1}x),\;\;\;\hat{\Phi}(i_v(t)^{-1}\xi),\;\;\;
\LLL_v[\Phi](x),\;\;\;\LLL_v[\hat{\Phi}](\xi)$$
vanish for all $x,\xi\in k\setminus
a\OO$ and $t\in k_v^\times\cap\OO_v$.
Convergence of the sums on the left side of (\ref{equation:KeyIdentity})
follows, and moreover it is justified to carry the summations
under the integral (\ref{equation:LSubVAlternativeDef}) representing
the operator $\LLL_v$. We
find that the left side of (\ref{equation:KeyIdentity}) equals
$$\int\frac{H_v(t)(\vartheta(t)-\Phi(0))
-\frac{1}{2}\one_{\OO_v^\times}(t)(\vartheta(1)-\Phi(0))}
{\norm{1-t}_v}d\mu_v^\times(t)
$$
$$
+
\int\frac{H_v(t)(\norm{t}_v\vartheta(t^{-1})-\hat{\Phi}(0))
-\frac{1}{2}\one_{\OO_v^\times}(t)(\vartheta(1)-\hat{\Phi}(0))}
{\norm{1-t}_v}d\mu_v^\times(t),
$$
 whence the result after making a substitution of $t^{-1}$
for $t$ in the second integral and then collecting like
terms.
 \qed

\subsection{Proof of the
theorem}\label{equation:EndOfAdelicStirlingProof} Neither side of
(\ref{equation:AdelicStirling}) changes
under the replacement of $(a,\Phi)$ by $(1,\Phi^{(a)})$. We may therefore
assume without loss of generality that $a=1$.
Since
$\KKK_v$ anti-commutes with $\FFF$, the Poisson summation formula
gives the identity
$$\sum_{x\in k}\KKK_v[\Phi](x)+\sum_{\xi\in
k}\KKK_v[\hat{\Phi}](\xi)=0.$$ After rearranging the
terms in this last identity and applying 
(\ref{equation:ThetaPoisson}), (\ref{equation:MSubVAlternativeDef}) and
(\ref{equation:KeyIdentity}), we obtain 
(\ref{equation:AdelicStirling}) in the case $a=1$.  The proof of the
theorem is complete.
 \qed
\section{The rational Fourier transform and the Catalan symbol (adelic versions)}
\label{section:RFT}
We pick up again the ideas introduced in \S\ref{section:Duality}.

\subsection{The rational Fourier transform}

\subsubsection{}
With $\lambda_0:\FF_q\rightarrow\{-1,0,1\}\subset\CC$ as in
\S\ref{subsection:TrivialIdentities}, put
$$\ee_0=\lambda_0\circ R_\omega:\adeles\rightarrow
\{-1,0,1\}\subset \CC$$ and then given any $\Phi\in \sch(\adeles)$,
put
$$\FFF_0[\Phi](\xi)=\tilde{\Phi}(\xi)=\int
\Phi(x)\ee_0(x\xi)d\mu(x)=
\int \Phi(x)\lambda_0(\langle x,\xi\rangle)d\mu(x)
$$
for all $\xi\in \adeles$, thereby defining  
the {\em rational Fourier transform}
$\tilde{\Phi}$ of
$\Phi$, also denoted by
$\FFF_0[\Phi]$. Via the trivial identities (\ref{equation:lambdaIdentities}) we have
\begin{equation}\label{equation:RelationOfFTs}
\begin{array}{rcl}
\hat{\Phi}(\xi)
&=&-\sum_{c\in \FF_q^\times}\lambda(-c)\tilde{\Phi}(c^{-1}\xi),
\\
\tilde{\Phi}(\xi)&=&
q^{-1}\sum_{c\in \FF_q^\times}(1-\lambda(c))\hat{\Phi}(c\xi)
\end{array}
\end{equation}
for all $\xi\in \adeles$.
Thus  $\hat{\Phi}$ and  $\tilde{\Phi}$ are expressible each
in terms of the other. It follows that $\FFF_0$ preserves
the Schwartz space $\sch(\adeles)$, just as does the standard
Fourier transform $\FFF$. 
It follows from the trivial identities 
(\ref{equation:lambdaIdentities}) and the example
(\ref{equation:GlobalFourierExample}) that
\begin{equation}\label{equation:RationalFourierExample}
\FFF_{0}[\one_{x+b\OO}]=
q^{1-g}\norm{b}\left(\one_{b^{-1}\kappa^{-1}\OO\cap[\langle
x,\cdot\rangle=0]}-\one_{b^{-1}\kappa^{-1}\OO\cap[\langle
x,\cdot\rangle=1]}\right)
\end{equation}
for all $b\in \adeles^\times$ and $x,\xi\in \adeles$.
It follows in turn via Serre duality that for all $b\in
\adeles^\times$,
$f\in
\adeles^\times
\cap\OO$ and
$\Phi\in
\sch(\adeles)$ such that $\Phi$ is constant on cosets of $b\OO$ and
supported in
$f^{-1}b\OO$, the rational Fourier transform
$\tilde{\Phi}$ is constant on cosets of
$fb^{-1}\kappa^{-1}\OO$ and supported in $b^{-1}\kappa^{-1}\OO$.

\subsubsection{}
Fix $\Phi\in \sch(\adeles)$. 
By equations (\ref{equation:RelationOfFTs}) relating standard and rational
Fourier transforms, the inversion
formula (\ref{equation:AdelicFourierInversion}) for the standard
Fourier transform, and the trivial identities (\ref{equation:lambdaIdentities}),
 we have
$$\begin{array}{rcl}
\Phi(x)&=&{ \int \ee(x\xi)\left(-\sum_{c\in
\FF_q^\times}\lambda(-c)
\tilde{\Phi}(c^{-1}\xi)\right)d\mu(\xi)}\\
&=&{
\int\left(-\sum_{c\in
\FF_q^\times}\lambda(-c)\ee(cx\xi)\right)\tilde{\Phi}(\xi)d\mu(\xi)}\\
&=&{\int\left(q\ee_0(x\xi)-\sum_{c\in
\FF_q^\times}\ee_0(cx\xi)
\right)\tilde{\Phi}(\xi)d\mu(\xi)}
\end{array}
$$
for all $x\in \adeles$; thus we can invert the rational Fourier transform.

\subsubsection{}
Given $\Phi\in \sch(\adeles)$ and $x\in \adeles$ we
put
$$\NNN[\Phi](x)=\sum_{c\in \FF_q^\times}\Phi(cx),$$
thereby defining
$$\NNN[\Phi]\in \sch(\adeles).$$
Then the inversion formula above for the rational Fourier transform can be rewritten more
compactly  in the form 
$$(q-\NNN)\FFF_0^2=1.$$
The related operator identities
$$\NNN\FFF_0=\FFF_0\NNN=\FFF\NNN=\NNN\FFF,\;\;\;
q=(q-\NNN)(1+\NNN)$$
and the {\em squaring rule}
\begin{equation}\label{equation:RelatedOpIdentitiesBis}
\FFF_{0}^2=q^{-1}(1+\NNN)
\end{equation}
are easy to verify.  We omit the
proofs. The last formula should be compared with 
the squaring rule (\ref{equation:AdelicFourierSquared}) for the adelic
Fourier transform. (If $q=2$, then the two formulas are the same.)

\subsubsection{}
We have for the rational Fourier transform a {\em scaling rule} of exactly
the same form as that obeyed by the usual Fourier transform, namely
\begin{equation}\label{equation:AdelicRationalScaling}
\FFF_0[\Phi^{(a)}]=\norm{a}\tilde{\Phi}^{(a^{-1})}
\end{equation}
for all $\Phi\in \sch(\adeles)$ and $a\in \adeles^\times$.
Moreover, since $\NNN\FFF=\NNN\FFF_0$
we have
$$
\sum_{x\in k}\Phi(a^{-1}x)=
\sum_{x\in k}\Phi^{(a)}(x)=
\sum_{\xi\in k}\FFF_0[\Phi^{(a)}](\xi)
=\norm{a}\sum_{\xi\in k}\tilde{\Phi}(a\xi)
$$
for all $a\in \adeles^\times$ and $\Phi\in \sch(\adeles)$,
i.~e., the Poisson summation formula continues to hold
with the rational Fourier transform replacing the standard one.
Consequently we have a {\em functional equation}
\begin{equation}\label{equation:ThetaRationalFE}
\Theta(a,\Phi)=\norm{a}\Theta(a^{-1},\tilde{\Phi})
\end{equation}
for all $a\in \adeles^\times$ and $\Phi\in \sch(\adeles)$.

\subsection{The Catalan symbol}
\subsubsection{}\label{subsubsection:CatalanDef}
Recall that $\ZZ[1/q]$ (as in the discussion of
toy rational Fourier transforms and toy Catalan symbols)
is the ring consisting of those rational numbers
$x$ such that $q^nx\in \ZZ$ for $n\gg 0$.
Let $\sch_0(\adeles)$ denote the space of
$\ZZ[1/q]$-valued Schwartz functions on $\adeles$,
and let $\sch_{00}(\adeles)$ denote the subgroup
consisting of functions $\Phi$ such that $\Phi(0)=0=\tilde{\Phi}(0)$.
Clearly, the groups
$\sch_0(\adeles)\supset \sch_{00}(\adeles)$ are stable under the
action of the operators
$\FFF_0$ and
$\NNN$. Note that $\Theta(a,\Phi)\in \ZZ[1/q]$ 
for all $a\in \adeles^\times$ and $\Phi\in \sch_0(\adeles)$.
Note that $x^n$ makes
sense for all
$x\in
k_\perf^\times$
and $n\in \ZZ[1/q]$. Note that
$\norm{a}\in q^\ZZ$ for all
$a\in
\adeles^\times$  and hence 
$x^{\norm{a}}$ makes sense for all $x\in k_\perf$.
For all places $v$ of $k$ and $x\in k_\perf^\times$ we define
$\ord_v x=q^{-n}\ord_v x^{q^n}\in \ZZ[1/q]$ for all sufficiently large
$n$. 
\subsubsection{}\label{subsubsection:CatalanFormalProperties} For all
$\Phi\in
\sch_{0}(\adeles)$ and $a\in \adeles^\times$, put
$$\begin{array}{lcl}
\left(\begin{array}{c}a\\\Phi\end{array}\right)_+&=&\displaystyle
\prod_{x\in k^\times}x^{\Phi(a^{-1}x)}\in
k^\times_\perf,\\\\
\left(\begin{array}{c}a\\\Phi\end{array}\right)&=&\displaystyle
\left(\begin{array}{c}a\\\Phi\end{array}\right)
_+\left(\begin{array}{c}a^{-1}\\\tilde{\Phi}\end{array}\right)_+^{\norm{a}}\in
k^\times_\perf.
\end{array}$$ 
In the  infinite products only finitely many terms differing from $1$
occur, so these objects are well-defined.
We call $\left(\begin{subarray}{c}
\cdot\\
\cdot
\end{subarray}\right)$ the {\em Catalan symbol},
and $\left(\begin{subarray}{c}
\cdot\\
\cdot
\end{subarray}\right)_+$ the {\em partial Catalan symbol}.
The rationale for the terminology is the same 
as the one offered  in the case of the toy Catalan symbol, namely the values of the Catalan
and partial Catalan symbol bear a certain structural resemblance to the Catalan numbers
$\frac{1}{n+1}\left(\begin{subarray}{c}
2n\\n
\end{subarray}\right)$.

\begin{Proposition}\label{Proposition:GenericCatalanBehavior}
Fix $\Phi\in \sch_{0}(\adeles)$.
Fix $b\in \adeles^\times$ and $f\in \adeles^\times\cap\OO$ such that
$\Phi$ is constant on cosets of $b\OO$ and supported in $f^{-1}b\OO$. 
Then
\begin{equation}
\left(\begin{array}{c}
a\\
\Phi
\end{array}\right)=
\left\{\begin{array}{ll}
\left(\begin{array}{c}
a^{-1}\\
\tilde{\Phi}
\end{array}\right)_+^{\norm{a}}&\mbox{if $\norm{ab}<\norm{f}$}\\
\left(\begin{array}{c}
a\\
\Phi
\end{array}\right)_+&\mbox{if $\norm{ab}>q^{2g-2}$}
\end{array}\right.\end{equation}
for all $a\in \adeles^\times$.
\end{Proposition}
\proof In parallel to formula (\ref{equation:ThetaAsymptotics}),
we have
$$
\left(\begin{array}{c}
a\\
\Phi
\end{array}\right)_+=
\prod_{x\in k^\times\cap af^{-1}b\OO}x^{\Phi(a^{-1}x)},\;\;
\left(\begin{array}{c}
a^{-1}\\
\tilde{\Phi}
\end{array}\right)_+=
\prod_{\xi\in k^\times\cap
a^{-1}\kappa^{-1}b^{-1}\OO}\xi^{\tilde{\Phi}(a\xi)},
$$
whence the result by Artin's product formula and the definitions.
\qed

\subsection{Formal properties of the Catalan symbol}
\subsubsection{}
Given $a,b\in \adeles^\times$, $x\in k^\times$, and $\Phi\in \sch_0(\adeles)$, we claim that
\begin{equation}\label{equation:IdeleClassInvarianceBis}
\left(\begin{array}{c}
ax\\
\Phi
\end{array}\right)=
\left(\begin{array}{c}
a\\
\Phi
\end{array}\right)x^{-\Phi(0)+\norm{a}\tilde{\Phi}(0)},
\end{equation}
\begin{equation}\label{equation:CatalanScaling}
\left(\begin{array}{c}
ab\\
\Phi
\end{array}\right)=\left(\begin{array}{c}
a\\
\Phi^{(b)}
\end{array}\right)=\left(\begin{array}{c}
1\\
\Phi^{(ab)}
\end{array}\right),
\end{equation}
and
\begin{equation}\label{equation:GFE2}
\left(\begin{array}{c}
a^{-1}\\
\tilde{\Phi}\end{array}\right)^{\norm{a}}=
(-1)^{\Theta(a,\Phi)-\Phi(0)}\left(\begin{array}{c}
a\\
\Phi
\end{array}\right).
\end{equation}
Note in particular that if $\Phi\in \sch_{00}(\adeles)$, then 
$\left(\begin{array}{c}
a\\
\Phi
\end{array}\right)$ depends only on the image of $a$ in the idele class
group $\adeles^\times/k^\times$. The last two relations we call the {\em scaling rule}
and {\em functional equation}, respectively, obeyed by the Catalan symbol.
\subsubsection{}\label{subsubsection:RoutineCatalanCalculations}
Here are the proofs of the claims. The calculation
$$
\begin{array}{rcl}
\left(\begin{array}{c}
ax\\
\Phi
\end{array}\right)&=&
\prod_{y\in k^\times}y^{\Phi(a^{-1}x^{-1}y)}\cdot \prod_{\eta\in
k^\times} \eta^{\tilde{\Phi}(ax\eta)\norm{ax}}\\\\
&=&
\prod_{y\in k^\times}(xy)^{\Phi(a^{-1}y)}\cdot \prod_{\eta\in
k^\times} (x^{-1}\eta)^{\tilde{\Phi}(a\eta)\norm{ax}}\\\\
&=&\left(\begin{array}{c}
a\\
\Phi
\end{array}\right)\cdot x^{\Theta(a,\Phi)-\Phi(0)-
\norm{ax}\Theta(a^{-1},\tilde{\Phi})+\norm{ax}\tilde{\Phi}(0)}\\\\
&=&\left(\begin{array}{c} a\\
\Phi
\end{array}\right)x^{-\Phi(0)+\norm{a}\tilde{\Phi}(0)}
\end{array}
$$
proves (\ref{equation:IdeleClassInvarianceBis}).
At the last equality we applied
Artin's product formula 
and functional equation (\ref{equation:ThetaRationalFE}).  The
calculation
$$
\begin{array}{rcl}
\left(\begin{array}{c}
ab\\
\Phi
\end{array}\right)&=&
\displaystyle
\prod_{x\in k^\times}x^{\Phi(a^{-1}b^{-1}x)}
\cdot \prod_{\xi\in k^\times} \xi^{\tilde{\Phi}(ab\xi)\norm{ab}}\\\\
&=&\displaystyle\prod_{x\in k^\times} x^{\Phi^{(b)}(a^{-1}x)}
\cdot \prod_{\xi\in k^\times}
\xi^{\FFF_0[\Phi^{(b)}](a\xi)\norm{a}}=\left(\begin{array}{c} a\\
\Phi^{(b)}
\end{array}\right)
\end{array}
$$
proves
(\ref{equation:CatalanScaling}). At the second equality 
we applied  scaling
rule (\ref{equation:AdelicRationalScaling}). 
Finally, the calculation
$$\begin{array}{rcl}
\left(\begin{array}{c}
a^{-1}\\
\tilde{\Phi}
\end{array}\right)^{q\norm{a}}
&=&
\displaystyle
\prod_{x\in k^\times}x^{\tilde{\Phi}(ax)\norm{a}q}
\cdot \prod_{\xi\in k^\times}
\xi^{q\FFF_0^2[\Phi](a^{-1}\xi)}\\\\
&=&
\displaystyle
\prod_{x\in k^\times}x^{\tilde{\Phi}(ax)\norm{a}q}
\cdot \prod_{\xi\in k^\times}
\xi^{(1+\NNN)[\Phi](a^{-1}\xi)}\\\\
&=&\displaystyle\left(\begin{array}{c}
a\\
\Phi
\end{array}\right)^q \prod_{c\in \FF_q^\times}\prod_{\xi\in
k^\times}
c^{\Phi(a^{-1}\xi)}\\\\&=&
\left(\begin{array}{c}
a\\
\Phi
\end{array}\right)^q(-1)^{\Theta(a,\Phi)-\Phi(0)}
\end{array}
$$
proves (\ref{equation:GFE2}). At the second equality we applied
formula (\ref{equation:RelatedOpIdentitiesBis}) for $\FFF_0^2$.
\subsubsection{}
For all $a\in \adeles^\times$, $x\in k^\times$ and $\Phi\in \sch_{00}(\adeles)$
we have
\begin{equation}\label{equation:PartialCatalanScaling}
\left(\begin{array}{c}
ax\\
\Phi
\end{array}\right)_+=\left(\begin{array}{c}
a\\
\Phi^{(x)}
\end{array}\right)_+=x^{\Theta(a,\Phi)}\left(\begin{array}{c}
a\\
\Phi
\end{array}\right)_+.
\end{equation}
We note this for convenient reference. The proof is more or less the same
as for (\ref{equation:IdeleClassInvarianceBis}).

\subsubsection{} The rational Fourier transform and the Catalan symbol
do in fact depend on the choice of differential $\omega$, even though we
have suppressed reference to $\omega$ in the notation. To clarify this
dependence let us temporarily write
$\FFF_{0,\omega}$ and $\left(\begin{subarray}{c}
\cdot\\
\cdot
\end{subarray}\right)_\omega$. Then for all $t\in k^\times$, $\Phi\in \sch(\adeles)$
and $\xi\in \adeles$ we have
\begin{equation}
\FFF_{0,t\omega}[\Phi](\xi)=\FFF_{0,\omega}[\Phi](t\xi)=\FFF_{0,\omega}[\Phi]^{(t^{-1})}(\xi).
\end{equation}
In turn, for all $a\in \adeles^\times$, $\Phi\in \sch_{00}(\adeles)$
and $t\in k^\times$,  we have
\begin{equation}\label{equation:CatalanChangeOfDifferential}
\left(\begin{array}{c}
a\\
\Phi
\end{array}\right)_{t\omega}=t^{-\Theta(a,\Phi)}\left(\begin{array}{c}
a\\
\Phi
\end{array}\right)_\omega
\end{equation}
by (\ref{equation:PartialCatalanScaling}) 
above, functional equation (\ref{equation:ThetaRationalFE})
obeyed by the theta symbol, and the definitions.

\subsubsection{}
Rewritten in terms of the Catalan symbol and the rational Fourier
transform, with attention restricted to Schwartz functions taken
from the group
$\sch_{00}(\adeles)$, the adelic Stirling formula
(\ref{equation:AdelicStirling}) takes a greatly simplified form,
namely
\begin{equation}\label{equation:AdelicStirlingBis}
\begin{array}{cl}
&{\displaystyle 
\Theta(a,\Phi)\ord_v \omega+\ord_v\left(\begin{array}{c}
a\\
\Phi
\end{array}\right)}\\&\\ =&{\displaystyle
\int\frac{\Theta(i_v^\times(t)a,\Phi)-\one_{\OO_v^\times}(t)\Theta(a)}{\norm{1-t}_v}d\mu_v^\times(t)}\\&\\
&{\displaystyle +\int\Phi(i_v(t))d\mu_v^\times(t)+\norm{a}\int
\tilde{\Phi}(i_v(t))d\mu_v^\times(t)}
\end{array}
\end{equation}
for all $a\in \adeles^\times$, $\Phi\in \sch_{00}(\adeles)$
and places $v$ of
$k$.  Clearly, the
right side of (\ref{equation:AdelicStirlingBis}) is 
independent of the differential $\omega$. As a consistency check, note that 
(\ref{equation:CatalanChangeOfDifferential})  forces the
left side of (\ref{equation:AdelicStirlingBis}) to be independent of
$\omega$.

\begin{Theorem}\label{Theorem:GenericBehavior}
Fix
$$x_0\in \adeles,\;\;\;a,b\in
\adeles^\times,\;\;f\in \adeles^\times\cap\OO$$ such that
$$x_0\not\in b\OO,\;\;
fx_0\in b\OO,\;\;
\left\{\begin{array}{rcl}
\norm{ab}&>&q^{2g-2}\\
&\mbox{or}&\\
\norm{ab}&<&\norm{f}.
\end{array}\right.$$ Then the map
$$
\left(x+b\OO\mapsto 
\left\{
\begin{array}{cl}
\left(\begin{array}{c}
a\\
\one_{x+b\OO}-\one_{x_0+b\OO}\end{array}\right)&\mbox{if $x\not\in
b\OO$}\\ 0&\mbox{if $x\in b\OO$}
\end{array}\right.
\right):f^{-1}b\OO/b\OO\rightarrow k_\perf
$$
is $\FF_q$-linear.
\end{Theorem}
\noindent Thus, besides obvious $\ZZ[1/q]$-linearity
 ($\sch_0(\adeles)$ to $k_\perf^\times$), the Catalan symbol has a
``hidden'' $\FF_q$-linearity. 
As the proof shows, the theorem is essentially just a rehash of 
Propositions~\ref{Proposition:EasyMooreApp} and 
\ref{Proposition:DualityClincher}.
\proof By scaling rule (\ref{equation:CatalanScaling}), we may assume
without loss of generality that
$a=1$. Let
$\gamma_{b,f,x_0}$ denote the map in question (with $a=1$). Our task is to
prove that $\gamma_{b,f,x_0}$ is $\FF_q$-linear. To do so, we distinguish
two cases, namely: (i) $\norm{b}>q^{2g-2}$ and (ii) $\norm{b}<\norm{f}$.

We turn to case (i). Put
$$V=f^{-1}b\OO\cap k,\;\;\;W=b\OO\cap k.$$
For all
$x\in f^{-1}b\OO$, the quantity
$\card((x+b\OO)\cap k)=\Theta(1,\one_{x+b\OO})$ is positive and independent of $x$
by Proposition~\ref{Proposition:ThetaAsymptotics}.
It follows that the inclusion-induced natural map
$V/W\rightarrow f^{-1}b\OO/b\OO$
is bijective. It suffices to prove the $\FF_q$-linearity of the map
$\tilde{\gamma}_{b,f,x_0}$ obtained by following the isomorphism
$V/W\rightarrow f^{-1}b\OO/b\OO$ by $\gamma_{b,f,x_0}$. 
Choose $x_1\in V$ such that
$x_0+b\OO=x_1+b\OO$. By
Proposition~\ref{Proposition:GenericCatalanBehavior} and the
definition of the toy Catalan symbol we have
$$\tilde{\gamma}_{b,f,x_0}(x+W)=\left(\begin{array}{c}
1\\
\one_{x+b\OO}-\one_{x_0+b\OO}
\end{array}\right)_+=\left(\begin{array}{c}
\alpha\\
\one_{x+W}
\end{array}\right)\bigg/\left(\begin{array}{c}
\alpha\\
\one_{x_1+W}
\end{array}\right)$$
for all $x\in V\setminus W$,
where $\alpha$ is the inclusion $V\rightarrow k_\perf$. 
Therefore the map
$\tilde{\gamma}_{b,f,x_0}$ is $\FF_q$-linear
by Proposition~\ref{Proposition:EasyMooreApp}.
Thus case (i) is proved.

We turn to case (ii).
In the obvious way let us now identify the
space $\sch(f^{-1}b\OO/b\OO)$ of complex-valued functions on the finite set
$f^{-1}b\OO/b\OO$ with the
subspace of
$\sch(\adeles)$ consisting of functions constant on cosets of $b\OO$
and supported in $f^{-1}b\OO$. By example
(\ref{equation:RationalFourierExample}) and the remark following we
have
$$\FFF_0[\sch(f^{-1}b\OO/b\OO)]
\subset\sch(b^{-1}\kappa^{-1}\OO/
fb^{-1}\kappa^{-1}\OO).$$
But the spaces $f^{-1}b\OO/b\OO$ and
$b^{-1}\kappa^{-1}\OO/ fb^{-1}\kappa^{-1}\OO$ are Serre dual with respect
to the pairing
$\langle\cdot,\cdot\rangle$ and so we also have at our disposal
a toy rational Fourier transform
$$\sch(f^{-1}b\OO/b\OO)\xrightarrow{\FFF_0^\toy}\sch(b^{-1}\kappa^{-1}\OO/
fb^{-1}\kappa^{-1}\OO).$$
By comparing the examples (\ref{equation:ToyRationalFourierExample}) and
(\ref{equation:RationalFourierExample}) it can be seen that there exists
an integer
$\ell$ such that
$$q^\ell \FFF^\toy_0[\Phi]=\FFF_0[\Phi]$$
for all $\Phi\in \sch(f^{-1}b\OO/b\OO)\subset \sch(\adeles)$.
Fix $\xi_0\in b^{-1}\kappa^{-1}\OO\setminus
fb^{-1}\kappa^{-1}\OO$ arbitrarily.
By case (i) already proved,
the map
$$\gamma_{b^{-1}\kappa^{-1}f\OO,f,\xi_0}:b^{-1}\kappa^{-1}\OO/
fb^{-1}\kappa^{-1}\OO\rightarrow k_\perf$$
 is $\FF_q$-linear. Let $\beta$ be the $\FF_q$-linear map obtained by multiplying 
 $\gamma_{b^{-1}\kappa^{-1}f\OO,f,\xi_0}$ by the factor
 $\left(\begin{array}{c}
 1\\
 \one_{\xi_0+fb^{-1}\kappa^{-1}\OO}\end{array}\right)_+$.
 By Proposition~\ref{Proposition:GenericCatalanBehavior},
 the map $\beta$ takes the form
 $$\xi+fb^{-1}\kappa^{-1}\OO
\mapsto 
\left\{\begin{array}{cl}
\left(\begin{array}{c}
1\\
\one_{\xi+fb^{-1}\kappa^{-1}\OO}
\end{array}\right)_+&\mbox{if $\xi\not\in fb^{-1}\kappa^{-1}\OO$,}\\
0&\mbox{otherwise.}
\end{array}\right.$$
Note that $\beta$ is injective. By Proposition~\ref{Proposition:GenericCatalanBehavior}
and the definition of the toy Catalan symbol we have
$$\begin{array}{rcl}
\gamma_{b,f,x_0}(x+b\OO)&=&
\left(\begin{array}{c}
1\\
\FFF_0[\one_{x+b\OO}-\one_{x_0+b\OO}]
\end{array}\right)_+\\\\
&=&
\left(\begin{array}{c}
\beta\\
\FFF_0^\toy[\one_{x+b\OO}-\one_{x_0+b\OO}]
\end{array}\right)^{q^\ell}
\end{array}$$
for all $x\in f^{-1}b\OO\setminus b\OO$.
Therefore the map $\gamma_{b,f,x_0}$ is $\FF_q$-linear by
Proposition~\ref{Proposition:DualityClincher}. Thus case (ii) is proved.
\qed

\section{Formulation of a conjecture}
In this section we state the conjecture (Conjecture~\ref{TheConjecture})
to which we allude in the title of the paper.
The conjecture links theta and Catalan symbols to two-variable
algebraic functions of a certain special type. Following the conjecture
we make some amplifying remarks.

\subsection{The ring $\kk$}
\subsubsection{}
Put
$$\kk_0=k\otimes_{\FF_q}k,\;\;\;\Delta_0=
\ker\left((x\otimes y\mapsto xy):k\otimes_{\FF_q}k\rightarrow k\right).$$
It is not difficult to verify that the ring $\kk_0$ is a noetherian integrally closed domain of dimension one, and hence a Dedekind domain. Clearly, the ideal
$\Delta_0$ is maximal in $\kk_0$, and by construction has residue field
canonically identified with $k$. Let
$\hat{\kk}_0$ be the completion of
$\kk_0$ with respect to $\Delta_0$, and let $\hat{\Delta}_0$ be the
closure of
$\Delta_0$ in $\hat{\kk}_0$.  Then
$\hat{\kk}_0$ is a discrete valuation ring with residue field canonically
identified with
$k$ and with maximal ideal $\hat{\Delta}_0$.
Put $$\kk=(\hat{\kk}_0)_\perf,\;\;\;
\Delta=\bigcup \sqrt[q^n]{\hat{\Delta}_0}.$$ Then
$\kk$ is a (nondiscrete) valuation ring with residue field
canonically identified with
$k_\perf$ and with maximal ideal $\Delta$.

\subsubsection{}
We identify the universal derivation
$d:k\rightarrow\Omega$ with the map $$
(x\mapsto x\otimes 1-1\otimes
x\bmod{\Delta_0^2}):k\rightarrow\Delta_0/\Delta_0^2=\hat{\Delta}_0/\hat{\Delta}_0^2
,$$
thus fixing an identification
$$\Omega=\hat{\Delta}_0/\hat{\Delta}_0^2.$$
We fix a generator $\varpi$ of the principal ideal
$\hat{\Delta}_0$ such that 
$$\omega\equiv\varpi\bmod{\hat{\Delta}_0^2}.$$
We call $\varpi$ a {\em lifting} of $\omega$.
\subsubsection{}
Every nonzero element $\varphi$ of the fraction field of $\hat{\kk}_0$
has a unique $\hat{\Delta}_0$-adic expansion of the form
$$\varphi=\sum_{i=i_0}^\infty (1\otimes a_i)\varpi^i\;\;\;\;(
i_0\in \ZZ, \;\;a_i\in k,\;\;a_{i_0}\neq 0)$$
and in terms of this expansion we define
$$\ord_\Delta \varphi=i_0,\;\;\;\lead_\Delta=a_{i_0}.$$ 
Also put $\ord_\Delta 0=+\infty$. 
More generally, for all nonzero $\varphi$ in the fraction field of $\kk$,
we define 
$$\ord_\Delta \varphi=q^{-n}\ord_\Delta \varphi^{q^n}\in \ZZ[1/q],\;\;\;
\lead_\Delta \varphi=(\lead_{\Delta}\varphi^{q^n})^{q^{-n}}\in
k_\perf^\times$$ for all sufficiently large $n$.
The function
 $\ord_\Delta$ is an additive valuation of the fraction field of $\kk$
independent of the choice of differential
$\omega$ and lifting
$\varpi$.  The function $\lead_\Delta$
is a homomorphism from the multiplicative group of the fraction field of
$\kk$ to $k_\perf^\times$ depending only on the differential
$\omega$, not on the lifting $\varpi$. 

\subsubsection{} We clarify the dependence of $\lead_\Delta$ on the
choice of differential $\omega$ as follows.
Let us temporarily write
$\lead_{\Delta,\omega}$ to stress the $\omega$-dependence. Then we
have
\begin{equation}\label{equation:LeadChangeOfDifferential}
\lead_{\Delta, t\omega}\varphi=t^{-\ord_\Delta
\varphi}\lead_{\Delta,\omega} \varphi
\end{equation}
for all $t\in k^\times$ and nonzero elements $\varphi$ of the fraction
field of $\kk$. This ought to be compared with equation
(\ref{equation:CatalanChangeOfDifferential}) above.

\subsubsection{}\label{subsubsection:Dagger1}
The exchange-of-factors automorphism 
$$(x\otimes y\mapsto y\otimes x):\kk_0\rightarrow\kk_0
$$ 
preserves
 $\Delta_0$,
hence has a unique
$\hat{\Delta}_0$-adically continuous extension to an automorphism of
$\hat{\kk}_0$, and hence has a unique extension to an automorphism of the
fraction field of
$\kk$. The result of applying the latter
automorphism  to an element $\varphi$ of the fraction field of $\kk$
we denote by $\varphi^\dagger$.
We have
\begin{equation}\label{equation:DaggerSymmetryBis}
(\varphi^\dagger)^\dagger=\varphi,\;\;\;\ord_\Delta
\varphi^\dagger=\ord_\Delta\varphi,\;\;\;
\lead_\Delta\varphi^\dagger=(-1)^{\ord_\Delta\varphi}\lead_\Delta
\varphi
\end{equation}
for all nonzero $\varphi$ in the fraction field of $\kk$.

\subsection{The ring $\KK$}

\subsubsection{}\label{subsubsection:Renormalization} For any abelian
extension
$K/k$, let 
$$\rho=\rho_{K/k}:\adeles^\times \rightarrow \Gal(K/k)$$ be the
reciprocity law homomorphism of global class field theory defined according to the nowadays
commonly followed convention of Tate's paper \cite{TateBackground}. It is important
to bear in mind two facts concerning this convention,
which differs from the older convention of, say,
Artin-Tate \cite{ArtinTate}. Firstly, if
$v$ is a place of
$k$ unramified in $K$, and $\tau\in k_v^\times$ is a uniformizer at $v$,
then
$\rho_{K/k}(i_v^\times(\tau))\in \Gal(K/k)$ is a geometric
Frobenius element (inverse of the usual Artin symbol) at $v$.
Secondly,  we have
$$C^{\rho_{K/k}(a)}=C^{\norm{a}}$$ for all $a\in \adeles^\times$ and
constants
$C\in K$.

\subsubsection{}
Let $k^\ab$ be the abelian closure of $k$ in the algebraic closure $\bar{k}$. 
Let $\FF_q^\ab$ be the abelian (also the algebraic) closure of $\FF_q$ in $\bar{k}$.
Note that
$\Gal(k^\ab/k)=\Gal(k^\ab_\perf/k_\perf)$ and that the natural action
of $\Gal(k^\ab/k)$ on $k^\ab_\perf$ commutes with the $q^{th}$ power
automorphism of $k^\ab_\perf$.  The group
$$\left\{[\sigma,\tau]\in\Gal(k^\ab/k)\times\Gal(k^\ab/k)\left|
\sigma\vert_{\FF_q^\ab}=
\tau\vert_{\FF_q^\ab}\right.\right\}$$
acts naturally 
on
the integral domains
$$k^\ab\otimes_{\FF^\ab_q}k^\ab,\;\;\;
k^\ab_\perf\otimes_{\FF^\ab_q}k^\ab_\perf$$
by the rule 
$$(x\otimes y)^{[\sigma,\tau]}=x^\sigma\otimes y^\tau,$$
with fixed rings $\kk_0$ and $\kk$, respectively.
Put
$$\delta \Gal(k^\ab/k)=\{[\sigma,\sigma]\vert \sigma\in
\Gal(k^\ab/k)\},$$
$$\KK_0=(k^\ab\otimes_{\FF^\ab_q}k^\ab)^{\delta\Gal(k^\ab/k)},\;\;\;
\KK=(k^\ab_\perf\otimes_{\FF^\ab_q}k^\ab_\perf)^{\delta\Gal(k^\ab/k)}.$$
Note that $\KK=(\KK_0)_\perf$. 
Note that $\Gal(k^\ab/\FF_q^\ab)$ may be identified with the Galois groups
of the etale ring extensions $\KK_0/\kk_0$ and 
$\KK/\kk$ under the map $\sigma\mapsto [\sigma,1]$.

\subsubsection{}\label{subsubsection:Dagger2}
Let $a\in \adeles^\times$ be given. The automorphism
$$\left(x\otimes y\mapsto
x^{\rho(a)}\otimes y^{\norm{a}}\right):
k^\ab_\perf\otimes_{\FF^\ab_q}k^\ab_\perf
\rightarrow k^\ab_\perf\otimes_{\FF^\ab_q}k^\ab_\perf
$$
commutes with the action of the group
$\delta\Gal(k^\ab/k)$, hence descends to an automorphism
of $\KK$ and hence has a unique extension to an automorphism of the
fraction field of
$\KK$. The result of applying the latter automorphism to an element
$\varphi$ of the fraction field of $\KK$ we denote by $\varphi^{(a)}$.
Thus we equip the fraction field of
$\KK$ with an action of $\adeles^\times$ factoring through the idele
class group $\adeles^\times/k^\times$.
\subsubsection{} The exchange-of-factors automorphism
$$(x\otimes y\mapsto
y\otimes x):k^\ab_\perf\otimes_{\overline{\FF}_q}k^\ab_\perf
\rightarrow k^\ab_\perf\otimes_{\overline{\FF}_q}k^\ab_\perf
$$
commutes with the action of the group $\delta\Gal(k^\ab/k)$,
hence descends to an automorphism of $\KK$,
and hence has a unique extension to an automorphism of the fraction field
of $\KK$. The result of applying the latter automorphism to an element
$\varphi$ of the fraction field of $\KK$ we denote by $\varphi^\dagger$. 
We have
\begin{equation}\label{equation:DaggerSymmetry}
(\varphi^\dagger)^{(a)}=((\varphi^{(a^{-1})})^\dagger)^{\norm{a}}
\end{equation}
for all $\varphi$ in the fraction field of $\KK$
and $a\in \adeles^\times$.
\begin{Proposition}
(i) There exists a unique $\kk_0$-algebra embedding 
$$\iota:\KK\rightarrow \kk$$ such that
$$\iota^{-1}(\Delta)=
\KK\cap\ker\left((x\otimes y\mapsto xy):k^\ab_\perf
\otimes_{\FF_q^\ab}k^\ab_\perf
\rightarrow k^\ab_\perf\right).$$
(ii) Moreover, $\iota$ commutes with $\dagger$.
\end{Proposition}
\noindent After the proof of this proposition,
in order to keep notation simple, we just identify
$\KK$ with its image in $\kk$ under $\iota$, and dispense
with the notation $\iota$ altogether. 
In this way
$\ord_\Delta
\varphi$ and $\lead_\Delta \varphi$ are defined for all nonzero $\varphi$
in the fraction field of $\KK$, and moreover
both (\ref{equation:DaggerSymmetryBis}) and
(\ref{equation:DaggerSymmetry}) hold.
\proof Part (i) granted, the $\kk_0$-algebra embedding
$\varphi\mapsto (\iota(\varphi^\dagger))^\dagger$ has the property
uniquely characterizing $\iota$ and hence coincides with $\iota$.
It is enough to prove part (i). 
In turn, it is enough to prove that there exists a unique
$\kk_0$-algebra embedding
$\iota_0:\KK_0\rightarrow \hat{\kk}_0$
such that
$$\iota_0^{-1}(\hat{\Delta}_0)=
\KK_0\cap\ker\left((x\otimes y\mapsto xy):k^\ab\otimes_{\FF_q^\ab}
k^\ab\rightarrow k^\ab\right),$$
because once this is proved, the induced homomorphism
$$(\iota_0)_\perf:(\KK_0)_\perf\rightarrow (\hat{\kk}_0)_\perf$$
is the only possibility for $\iota$.

Now let $K/k$ be any finite abelian extension
and consider the following naturally associated objects:
\begin{itemize}
\item Let $G=\Gal(K/k)$.
\item Let $\FF$ be the
constant field of
$K$.  
 \item Let 
$G'=\{[\sigma,\tau]\in G\times G\vert\sigma\vert_{\FF}=
\tau\vert_\FF\}$.
\item Let $G'$ act on $K\otimes_\FF K$ by the rule
$(x\otimes y)^{[\sigma,\tau]}=x^\sigma\otimes y^\tau$.
\item Let $\delta G=\{[\sigma,\sigma]\vert \sigma \in G\}\subset G'$.
\item Put $R=(K\otimes_\FF K)^{\delta G}$.
\item Put $I=R\cap \ker\left((x\otimes y\mapsto
xy):K\otimes_\FF K\rightarrow K\right)$.
\end{itemize}
It is enough to prove that there exists a unique
$\kk_0$-algebra homomorphism
$\iota_{K/k}:R\rightarrow \hat{\kk}_0$
such that
$\iota_{K/k}^{-1}(\hat{\Delta}_0)=I$,
because once this is proved the limit
$\lim_{\rightarrow}
\iota_{K/k}:\KK_0\rightarrow \hat{\kk}_0$
extended over all finite subextensions $K/k$
of $k^\ab/k$ is the only possibility
for $\iota_0$. Now the ring $K\otimes_\FF K$ is a Dedekind
domain finite etale and abelian over $\kk_0$.
It follows that the $\delta G$-fixed subring $R$ has these same
properties by descent. 
The prime $I\subset R$ lies above the prime $\Delta_0\subset\kk_0$
and we have $R/I=\kk_0/\Delta_0=k$. Since $R$ is etale over $\kk_0$,
existence and uniqueness of
$\iota_{K/k}$ follow now by the infinitesimal lifting criterion. 
\qed

\subsection{Coleman units}
We declare a nonzero element $\varphi$ of the fraction field
of $\KK$ to be a {\em Coleman unit} if for every maximal
ideal $M\subset \KK$ not of the form $\{\psi\in \KK\vert \ord_\Delta
\psi^{(a)}>0\}$ for some
$a\in \adeles^\times$, we have
$\varphi\in \KK_M^\times$, where $\KK_M$ is the local ring of $\KK$ at $M$.
The
Coleman units form a group under multiplication. For every
Coleman unit
$\varphi$ and $a\in \adeles^\times$, again $\varphi^{(a)}$ is a Coleman
unit, and so also is $\varphi^\dagger$. The group of Coleman units is
closed under the extraction of
$q^{th}$ roots.  The inspiration for this definition
comes from Coleman's paper \cite{Coleman}.

\begin{Conjecture}\label{TheConjecture}
For every $\Phi\in \sch_{00}(\adeles)$,
there exists a unique Coleman unit $\varphi$ such that 
\begin{equation}\label{equation:ConjecturalAttachment}
\ord_\Delta
\varphi^{(a)}=\Theta(a,\Phi),\;\;\;
\lead_\Delta\varphi^{(a)}=\left(\begin{array}{c}
a\\
\Phi
\end{array}\right)
\end{equation}
for all $a\in \adeles^\times$.
\end{Conjecture}
\subsection{Amplification} \label{subsection:CCalculus}
\subsubsection{} We have, so we claim, the following {\em uniqueness
principle}: for every $\psi$ in the fraction field of $\KK$,
\begin{equation}\label{equation:UniquenessPrinciple}
\card\{\norm{a}\;\vert\; 
a\in
\adeles^\times,\;
\ord_\Delta\psi^{(a)}\neq 0\}=\infty\Rightarrow \psi=0.
\end{equation}
In any case, we can find a Dedekind domain $R$ between $\KK$ and $\kk_0$ which is finite over $\kk_0$ and to which $\psi$ belongs. Since the maximal ideals of
$\kk_0$ of the form 
$$\ker\left((x\otimes y\mapsto xy^{q^n}):\kk_0\rightarrow
k_\perf\right)\;\;\;(n\in \ZZ)$$
are distinct, the hypothesis implies that $\psi$
has nonzero valuation at infinitely many maximal ideals of $R$ and hence
vanishes identically. The claim is proved. By the uniqueness principle, for each
$\Phi\in\sch_{00}(\adeles)$, there can be at most one element $\varphi$
of the fraction field of $\KK$ satisfying
(\ref{equation:ConjecturalAttachment}) for all $a\in
\adeles^\times$. So the
uniqueness asserted in the conjecture is clear. Only existence is at
issue.

\subsubsection{}\label{subsubsection:SimNotation}
Let nonzero $\varphi$ in the fraction field of $\KK$,
$\Phi\in
\sch_{00}(\adeles)$, and $a\in \adeles^\times$ be given.
We write
$\varphi\sim_a\CCC[\Phi]$ if (\ref{equation:ConjecturalAttachment})
holds for the data $\varphi$, $a$ and $\Phi$.
We write $\varphi=\CCC[\Phi]$
if $\varphi\sim_a\CCC[\Phi]$ for all $a\in \adeles^\times$ and $\varphi$ is a Coleman unit.
In the situation $\varphi=\CCC[\Phi]$, as the notation is meant to suggest, we think of $\varphi$
as the image of $\Phi$ under a transformation $\CCC$. 
From the latter point of view,
the conjecture asserts the well-definedness of this transformation
$\CCC$ on the whole of the space $\sch_{00}(\adeles)$.

\subsubsection{}
Let $\varphi$ in the fraction field of $\KK$ and
$\Phi\in \sch_{00}(\adeles)$ be given.
Comparison with (\ref{equation:ConjecturalAttachment}) of
the formula (\ref{equation:CatalanChangeOfDifferential})
explaining the dependence of $\left(\begin{subarray}{c}\cdot\\\cdot\end{subarray}\right)$
on $\omega$ and the formula (\ref{equation:LeadChangeOfDifferential}) explaining
the dependence of $\lead_\Delta$  on $\omega$ shows that 
if $\varphi=\CCC[\Phi]$, then no matter what nonzero differential $\omega\in \Omega_{k/\FF_q}$ we
choose to define rational Fourier transforms, Catalan symbols,
and ``leading Taylor coefficients at $\Delta$'', it remains the case that
$\varphi=\CCC[\Phi]$.  In short, $\CCC$ is independent of $\omega$.
So if the conjecture is true for one
choice of  $\omega$, it is true for all. 

\subsubsection{}
Grant the conjecture for this paragraph, so that the transformation
$\CCC$ is defined. Clearly, we have
\begin{equation}
\CCC[n_1\Phi_1+n_2\Phi_2]=\CCC[\Phi_1]^{n_1}\CCC[\Phi_2]^{n_2}
\end{equation}
for all $\Phi_1,\Phi_2\in \sch_{00}(\adeles)$
and $n_1,n_2\in \ZZ[1/q]$.
Comparison with 
(\ref{equation:ConjecturalAttachment}) of the scaling rule
(\ref{equation:ThetaScaling}) for the theta symbol and
the scaling rule (\ref{equation:CatalanScaling})
for the Catalan symbol shows that 
\begin{equation}\label{equation:CCCEquivariance}
\CCC[\Phi^{(a)}]=\CCC[\Phi]^{(a)}
\end{equation} 
for all $a\in \adeles^\times$ and $\Phi\in \sch_{00}(\adeles)$.
Comparison with (\ref{equation:ConjecturalAttachment})
of the functional equation (\ref{equation:ThetaRationalFE}) for the theta symbol, the functional equation
(\ref{equation:GFE2}) for the Catalan symbol, 
and the properties (\ref{equation:DaggerSymmetryBis}) and
(\ref{equation:DaggerSymmetry}) of the dagger operation on the fraction field of $\KK$ shows that
\begin{equation}\label{equation:ConjecturalAttachmentBis}
\CCC[\tilde{\Phi}]=\CCC[\Phi]^\dagger
\end{equation}
for all $\Phi\in \sch_{00}(\adeles)$. 
For any fixed $b\in \adeles^\times$ and $x_0\in \adeles\setminus b\OO$ the
map
\begin{equation}\label{equation:SecretLinearity}
\left(x+b\OO\mapsto
\left\{\begin{array}{rl}
\CCC\left[\one_{x+b\OO}-\one_{x_0+b\OO}\right]&\mbox{if $x\not\in
b\OO$}\\
0&\mbox{if $x\in b\OO$}
\end{array}\right.\right):\adeles/b\OO\rightarrow \KK
\end{equation}
is $\FF_q$-linear by the ``hidden'' $\FF_q$-linearity of the Catalan symbol (Theorem~\ref{Theorem:GenericBehavior}) and the  uniqueness principle
(\ref{equation:UniquenessPrinciple}).

\subsubsection{} We return now to the general discussion of the conjecture, no longer assuming that it holds. It is easy to verify that the set of $\Phi\in \sch_{00}(\adeles)$ for which $\CCC[\Phi]$ is defined is a $\ZZ[1/q][\adeles^\times]$-submodule of $\sch_{00}(\adeles)$. 
It follows that to prove the conjecture it is enough to fix a set of
generators for $\sch_{00}(\adeles)$ as a $\ZZ[1/q][\adeles^\times]$-module
 and for each generator $\Phi$ to prove the existence of $\CCC[\Phi]$. For
example, for any fixed
$b\in
\adeles^\times$ and $x_0\in \adeles\setminus b\OO$, the family
$\left\{\left.\one_{x+b\OO}-\one_{x_0+b\OO}\right| x\in
\adeles\setminus b\OO\right\}$ is large enough to generate
$\sch_{00}(\adeles)$ as a
$\ZZ[1/q][\adeles^\times]$-module.

\subsubsection{} Here are some trivial but handy cases in
which we can easily prove that the transformation $\CCC$ is defined.
Fix
$\Phi\in\sch_0(\adeles)$. Fix $x\in k^\times$.
Note that $\Phi^{(x)}-\Phi\in\sch_{00}(\adeles)$. We
claim that
\begin{equation}\label{equation:TrivialSolitonCalculation}
\CCC[\Phi^{(x)}-\Phi]=x^{-\Phi(0)}\otimes x^{\tilde{\Phi}(0)}.
\end{equation}
In any case, the right side is a unit of $\KK$
and {\em a fortiori} a Coleman unit; further, we have
$$\Theta(a,\Phi^{(x)}-\Phi)=0$$ for all $a\in \adeles^\times$
by (\ref{equation:ThetaIdeleClassDependence})
and (\ref{equation:ThetaScaling}); and finally we have
$$\left(\begin{array}{c}
a\\
\Phi^{(x)}-\Phi
\end{array}\right)=x^{-\Phi(0)+\norm{a}\tilde{\Phi}(0)}
$$
for all $a\in \adeles^\times$ by
(\ref{equation:IdeleClassInvarianceBis})
and (\ref{equation:CatalanScaling}).
This is enough to prove the claim.

\section{Proof of the conjecture in genus zero}
\label{section:GenusZero}
We assume in \S\ref{section:GenusZero} that $g=0$ and under this
additional assumption we are going to prove
Conjecture~\ref{TheConjecture}. 

\subsection{Notation and reductions} 
 Every function field  of genus zero
with field of constants $\FF_q$
 has a place with residue field
$\FF_q$, and as we noted in \S\ref{subsection:CCalculus},
the conjecture is invariant under change of differential.
We therefore may assume without loss of generality that
$$k=\FF_q(T),\;\;\;\omega=dT,$$
where $T$ is transcendental over $\FF_q$.
 We take 
$$\varpi=T\otimes 1-1\otimes T\in \kk_0$$
as a lifting of $\omega$. Note that $\varpi$ is already a prime element
of the principal ideal domain $\kk_0$. Note also that $\varpi$ is
a Coleman unit.
 Let $\infty$ be the unique place of $k$ at which $T$ has a pole and put $$\tau=i_{\infty}^\times(T^{-1})\in
\adeles^\times,$$
 noting that $$\norm{\tau}=q^{-1}.$$
 Put 
$$\UUU=\left\{a=[a_v]\in \OO^\times \left|\;
\norm{a_\infty-1}_\infty<1\right.\right\}.$$
Then every $a\in \adeles^\times$ factors uniquely in the form
$$a=uz\tau^{-N}\;\;\;(u\in \UUU,\;z\in k^\times,\;N\in \ZZ).$$
We use this factorization repeatedly below. 
 For each $x\in k$, write 
$$x=\lfloor x\rfloor+\langle x\rangle\;\;\;(\lfloor x\rfloor\in
\FF_q[T],\;\;\norm{\langle x\rangle}_\infty<1)$$
in the unique possible way, in parallel with
the definitions made in
\S\ref{subsubsection:Brackets}. 
 Put
$$\deg a=-\ord_\infty a=(\mbox{degree of $a$ as a polynomial in $T$})$$ for $a\in \FF_q[T]$.
 For each $x\in k^\times$, put
$$\Psi_x=\one_{x+\tau\OO}-\one_{1+\tau\OO}
\in \sch_{00}(\adeles).$$ 
The set $\{\Psi_x\}$ generates $\sch_{00}(\adeles)$
as a $\ZZ[1/q][\adeles^\times]$-module, and so to prove
the conjecture in the case at hand it suffices to prove that
$\CCC[\Psi_x]$ exists for every $x\in k^\times$.

\subsection{Calculation of theta and Catalan symbols}
Let $N\in \ZZ$ and $x\in k^\times$ be given. We 
calculate 
$\Theta(\tau^{-N},\Psi_x)$ and $\left(\begin{array}{c}
\tau^{-N}\\
\Psi_x
\end{array}\right)$.

\subsubsection{}
Assume at first that $N\geq 0$, which is the relatively easy case.
 We shall handle the case $N<0$ presently.
We have
$$\begin{array}{rcl}
k\cap \tau^{-N}(x+\tau\OO)&
=&\{z\in k\vert \;
\langle
x  \rangle=\langle z\rangle,\;\;
\lfloor x\rfloor=\lfloor T^{-N}z\rfloor
\}\\
&=&\langle x\rangle+T^N\lfloor x\rfloor
+\{a\in \FF_q[T]\vert\deg a<N\}.
\end{array}$$
Therefore we have
\begin{equation}\label{equation:PositiveThetaCalc}
\Theta(\tau^{-N},\Psi_x)=0,
\end{equation}
\begin{equation}\label{equation:PositiveCatalanCalc}
\left(\begin{array}{c}
\tau^{-N}\\\Psi_x\end{array}\right)=\left(\begin{array}{c}
\tau^{-N}\\\Psi_x\end{array}\right)_+
=\prod_{\begin{subarray}{c}
a\in \FF_q[T]\\
\deg a<N
\end{subarray}}
\frac{a+\langle x\rangle+T^N\lfloor x\rfloor}{a+T^N},
\end{equation}
by Proposition~\ref{Proposition:GenericCatalanBehavior}
at the first equality of (\ref{equation:PositiveCatalanCalc}),
and at the other two equalities by direct appeal to the definitions.

\subsubsection{}
We now take a close look at the rational Fourier transform $\tilde{\Psi}_x$
in order to prepare for handling the case $N<0$.
By specializing the example (\ref{equation:RationalFourierExample})
of an adelic rational Fourier transform to the present case
(we may take $\kappa=\tau^{-2}$ because $\omega=dT$ has a double pole at $\infty$ and no other poles or zeroes), we have
$$\FFF_0[\one_{x+\tau\OO}]=
\one_{\tau\OO\cap[\langle x,\cdot\rangle=0]}
-\one_{\tau\OO\cap[\langle x,\cdot\rangle=1]}
$$
and hence
\begin{equation}\label{equation:TildePsi}
\tilde{\Psi}_x=(1+\NNN)
\left[\Psi^\star_x\right]\;\;\;(\Psi^\star_x=-\one_{\tau\OO\cap[\langle
x,\cdot\rangle=1]}+\one_{\tau\OO\cap[\langle
1,\cdot\rangle=1]}),
\end{equation}
where to get (\ref{equation:TildePsi})
from the preceding formula we reused the ``starred''
trick from the proof of Proposition~\ref{Proposition:DualityClincher}.
It follows (without any restriction on $N$) that
\begin{equation}\label{equation:RationalPsiStar}
\Theta(\tau^{-N},\Psi_x)=q^N\Theta(\tau^{N},\tilde{\Psi}_x)
=q^{N+1}\Theta(\tau^N,\Psi^\star_x),
\end{equation}
by functional equation (\ref{equation:ThetaRationalFE}) at the first equality, and  scaling rule (\ref{equation:ThetaScaling}) at the second. It follows in turn (again without restriction on $N$) that 
\begin{equation}\label{equation:CarlitzRationalFourierTer}
\left(\begin{array}{c}
\tau^{N}\\
\tilde{\Psi}_x
\end{array}\right)_+=(-1)^{\Theta(\tau^{-N},\Psi_x)}\left(\begin{array}{c}
\tau^{N}\\
\Psi^\star_x
\end{array}\right)_+^q,
\end{equation}
by scaling rule (\ref{equation:PartialCatalanScaling}) followed by an application of (\ref{equation:RationalPsiStar}).
\subsubsection{}
We now handle the remaining case $N<0$ of our calculation. Put
$$\nu=-N-1,$$
noting that $\nu\geq 0$.
Put
$$V=\{a\in \FF_q[T]\vert \deg a\leq \nu\}=
\tau^{N+1}\OO\cap k,$$
$$\left.
\begin{array}{rcl}
\xi&=&(a\mapsto -\langle x,\tau^{-N}a\rangle)\\
\xi_1&=&(a\mapsto -\langle 1,\tau^{-N}a\rangle)
\end{array}\right\}:V\rightarrow\FF_q,$$
noting that 
$$\xi(a)=\Res_\infty(a
(\langle x\rangle-T^{-\nu-1}\lfloor x\rfloor)\,dT)),\;\;
\xi_1(a)=-\Res_\infty(a T^{-\nu-1}\,dT).$$
The latter relations are verified by applying ``sum of residues equals zero'' and then discarding terms not contributing to the residue at $T=\infty$. 
Since
\begin{equation}\label{equation:Intermediate1}
k\cap\tau^{N}(\tau\OO\cap [\langle x,\cdot\rangle=1])
=\{a\in V\vert \xi(a)=-1\},
\end{equation}
we have via (\ref{equation:RationalPsiStar}) that
\begin{equation}\label{equation:NegativeThetaCalc}
\begin{array}{cl}
&\Theta(\tau^{\nu+1},\Psi_x)\\\\
=&\left\{\begin{array}{cl}
0&\mbox{if $\xi\neq 0$,}\\
1&\mbox{if $\xi=0$,}
\end{array}\right.\\\\
=&\left\{\begin{array}{cl}
0&\mbox{if $\langle x\rangle-T^{-\nu-1}\lfloor x\rfloor\not\in T^{-\nu-2}\FF_q[[1/T]]+\FF_q[T]$,}\\
1&\mbox{if $\langle x\rangle-T^{-\nu-1}\lfloor x\rfloor\in T^{-\nu-2}\FF_q[[1/T]]+\FF_q[T]$.}
\end{array}\right.
\end{array}
\end{equation}
It follows that
\begin{equation}\label{equation:Intermediate2}
\begin{array}{cl}&k\cap \tau^{-N}(x+\tau\OO)\\
=&\{z\in k\vert \;
\langle
x  \rangle=\langle z\rangle,\;\;
\lfloor x\rfloor=\lfloor T^{\nu+1}z\rfloor
\}
\\
=&
\left\{\begin{array}{cl}
\emptyset&\mbox{if 
$\Theta(\tau^{-N},\Psi_x)=0$,}\\
\{\langle x\rangle+\lfloor T^{-\nu-1} x \rfloor
\}&\mbox{if 
$\Theta(\tau^{-N},\Psi_x)=1$.}
\end{array}\right.
\end{array}
\end{equation}
Finally, it follows via (\ref{equation:CarlitzRationalFourierTer},
\ref{equation:Intermediate1}, \ref{equation:NegativeThetaCalc},
\ref{equation:Intermediate2})  and the definitions that
\begin{equation}\label{equation:NegativeCatalanCalc}
\left(\begin{array}{c}
\tau^{\nu+1}\\\Psi_x\end{array}\right)^{q^\nu}=
\left\{\begin{array}{cl}
\displaystyle \prod_{\begin{subarray}{c} a\in V\\
\xi_1(a)=1
\end{subarray}}a\bigg/
\prod_{\begin{subarray}{c} a\in V\\
\xi(a)=1
\end{subarray}}a
&\mbox{if $\Theta(\tau^{\nu+1},\Psi_x)=0$,}\\\\
(\langle x\rangle+\lfloor T^{-\nu-1}x\rfloor)^{q^\nu}D_{\nu}&\mbox{if
$\Theta(\tau^{\nu+1},\Psi_x)=1$,}
\end{array}\right.
\end{equation}
 where 
$$D_\nu=-\prod_{\begin{subarray}{c}
a\in V\\
\xi_1(a)=-1
\end{subarray}}a=\prod_{\begin{subarray}{c}
a\in V\\
\xi_1(a)=1
\end{subarray}}a
=\prod_{i=0}^{\nu-1}(T^{q^\nu}-T^{q^i}),$$
cf.\ equation (\ref{equation:DN}).
The calculation is complete.

\subsection{The special point $\Pbold$}
We study a certain $\hat{k}_0$-valued point of the
group-scheme $\HHH$ introduced in \S\ref{section:FunctorsAndIdentities}.

\subsubsection{}
Given $u\in \UUU$, 
the set $k\cap (u+f\tau\OO)$ consists of a single element $a\in \FF_q[T]$ such that $\deg a<\deg f$. That noted,
it is clear that
there is a unique group isomorphism $$(u\mapsto [u]):\UUU\iso
\HHH(\FF_q)$$
such that
$$[u]_\infty(t)=u_\infty(t)\;\;\;(u_\infty=u_\infty(T)\in
1+(1/T)\FF_q[[1/T]]\subset\OO_\infty^\times)$$ and 
$$[u]_f(t)=a(t)\;\;\;(a=a(T)\in k\cap (u+f\tau\OO))$$
for every monic $f=f(T)\in \FF_q[T]$. 
 Let $k^\sep/k$ be a separable algebraic closure of $k$ 
of which $k^\ab/k$ is a subextension.  Let $\Frob:\HHH\rightarrow\HHH$ be the $q^{th}$ power Frobenius
endomorphism. 
According to geometric class field theory (see \cite{Serre}) the set of solutions $X\in \HHH(k^\sep)$ of the equation
\begin{equation}\label{equation:LangTorsor}
\Frob X=[t-T]X
\end{equation} 
forms an $\HHH(\FF_q)$-torsor, and for any solution $X$
we have an explicit reciprocity law
\begin{equation}\label{equation:ExplicitReciprocity}
X^{\rho(uz\tau^{-N})}=[u]X\;\;\;(u\in
\UUU,\;z\in k^\times, \;N\in \ZZ).
\end{equation}
In particular,
we always have $X\in
\HHH(k^\ab)$.

\subsubsection{}
Put
$$\Pbold=\prod_{i=0}^\infty
\left( [t-T^{q^i}\otimes
1]^{-1} [t-1\otimes T^{q^i}]\right)\in \HHH(\hat{\kk}_0).$$
The product is convergent because 
$$[t-T^{q^i}\otimes 1]^{-1}
[t-1\otimes T^{q^i}]\equiv 1\bmod{\varpi^{q^i}}$$
for all $i$. 
Since the Lang isogeny $(\Frob-1):\HHH\rightarrow\HHH$
is etale, the functional equation and congruence
\begin{equation}\label{equation:BigPCharacterization}
\Frob
\Pbold=[t-T\otimes 1][t-1\otimes T]^{-1}\Pbold,\;\;\;\Pbold\equiv
1\bmod{\varpi}
\end{equation} characterize $\Pbold$ uniquely
in $\HHH(\hat{\kk}_0)$.

\subsubsection{}
Let $X\in \HHH(k^\ab)$ be any solution of (\ref{equation:LangTorsor}),
let 
$$X\otimes 1\in \HHH(k^\ab\otimes_{\FF_q^\ab}k^\ab)$$ be the
image of
$X$ under the map induced by the ring homomorphism
$$(x\mapsto x\otimes 1):k^\ab\rightarrow
k^\ab\otimes_{\FF_q^\ab}k^\ab,$$ let 
$$1\otimes X\in \HHH(k^\ab\otimes_{\FF_q^\ab}k^\ab)$$ be
defined analogously, and put
$$X\otimes X^{-1}=(X\otimes 1)(1\otimes X)^{-1}.$$
Then, so we claim, we have
\begin{equation}
\label{equation:BoldPFactorization}
\Pbold=X\otimes X^{-1}
\in \HHH(\KK_0)\subset\HHH(\hat{\kk}_0).
\end{equation}
The equality holds because (i) $X\otimes X^{-1}$
satisfies the conditions (\ref{equation:BigPCharacterization})
uniquely characterizing $\Pbold$, and  (ii) $X\otimes X^{-1}$ is 
$\delta\Gal(k^\ab/k)$-invariant by the explicit reciprocity law
(\ref{equation:ExplicitReciprocity}).
 The
claim is proved.
\subsubsection{}
By combining the observations of the preceding three paragraphs, for any $a\in \adeles^\times$,
we can 
calculate the image
$\Pbold^{(a)}$ of $\Pbold$ under the map
$\HHH(\KK_0)\rightarrow\HHH(\KK)$ induced by the ring homomorphism
$(\varphi\mapsto
\varphi^{(a)}):\KK_0\rightarrow \KK$. We have
\begin{equation}\label{equation:BoldPTwist}
\Pbold^{(uz\tau^{-N-1})}
=[u]\Pbold
\left\{\begin{array}{ll}
\displaystyle\prod_{i=0}^{N}[t-1\otimes T^{q^i}]^{-1}&\mbox{if $N\geq
0$,}\\\\
\displaystyle\prod_{i=1}^{\nu}[t-1\otimes T^{q^{-i}}]&\mbox{if
$N<0$,}
\end{array}\right.
\end{equation}
for all $u\in \UUU$, $z\in
k^\times$ and
$N\in
\ZZ$, where $N=-\nu-1$. 

\begin{Lemma}\label{Lemma:StraightTheta}
For all $x=x(T)\in k^\times$ and $a\in \adeles^\times$ we have
\begin{equation}\label{equation:StraightTheta}
(\theta_{x(t)}(\Pbold))^{(\tau^{-1})}\sim_a\CCC[\Psi_x].
\end{equation}
\end{Lemma}
\noindent Recall that $\theta$ was defined in \S\ref {subsubsection:Brackets}
and $\sim$ in \S\ref{subsubsection:SimNotation}.
\proof Write
$$a=uz\tau^{-N}\in \adeles^\times\;\;\;(u\in \UUU,\;z\in
k^\times,\;N\in
\ZZ)$$
in the unique possible way. Without loss of generality we may assume that $z=1$. Fix monic $f=f(T)\in \FF_q[T]$ such that $fx\in \FF_q[T]$
and let $\tilde{x}=\tilde{x}(T)\in k$ be uniquely characterized by the relations
$$[u]_f(t)\langle x(t)\rangle \equiv \langle\tilde{x}(t)\rangle\bmod{\FF_q[t]},$$
$$
[u]_\infty(t)\lfloor x(t)\rfloor\equiv\lfloor \tilde{x}(t)\rfloor\bmod{(1/t)\FF_q[[1/t]]},$$
in which case (equivalently)
$$u(x+\tau\OO)=\tilde{x}+\tau\OO.$$ 
Let $\psi_x\in \KK_0$ be defined by the expression on the left side of 
 (\ref{equation:StraightTheta}).
Then we have
$$\psi_x^{(a)}=\theta_{x(t)}(\Pbold^{(a\tau^{-1})})=\theta_{x(t)}([u]\Pbold^{(\tau^{-N-1})})=
\theta_{\tilde{x}(t)}(\Pbold^{(\tau^{-N-1})})=
\psi_{\tilde{x}}^{(\tau^{-N})}$$
and $$\Psi_x^{(u\tau^{-N})}=\Psi_{\tilde{x}}^{(\tau^{-N})}.$$
After replacing $x$ by $\tilde{x}$,  we may simply assume
that $a=\tau^{-N}$. All the preceding reductions taken into account, it will be enough to prove that
\begin{equation}\label{equation:SimClincher}
\ord_\Delta\psi_x^{(\tau^{-N})}=\Theta(\tau^{-N},\Psi_x),\;\;\;
\lead_\Delta\psi_x^{(\tau^{-N})}=\left(\begin{array}{c}
\tau^{-N}\\
\Psi_x
\end{array}\right).
\end{equation}
If $N\geq 0$, then we have via
(\ref{equation:BigPCharacterization}) and (\ref{equation:BoldPTwist})
that
$$\psi_x^{(\tau^{-N})}=
\theta_{x(t)}\left(\Pbold^{(\tau^{-N-1})}\right)\equiv
\theta_{x(t)}\left(\prod_{i=0}^N[t-1\otimes T^{q^i}]^{-1}\right)
\bmod{\hat{\Delta}_0},$$
and this noted,  (\ref{equation:SimClincher}) follows
by comparing Proposition~\ref{Proposition:JacobiTrudiRedux} to formulas
(\ref{equation:PositiveThetaCalc}) and
(\ref{equation:PositiveCatalanCalc}). If $N=-\nu-1<0$, then we 
have via
(\ref{equation:BigPCharacterization}) and (\ref{equation:BoldPTwist}), along with the definition of $\Pbold$,
that
$$\begin{array}{rcl}
\left(\psi_x^{(\tau^{-N})}\right)^{q^\nu}&=&\theta_{x(t)}
(\Frob^\nu(\Pbold^{(\tau^{-N-1})}))\\
&=&\theta_{x(t)}\left((\Frob^\nu\Pbold) \prod_{i=0}^{\nu-1}[t-1\otimes T^{q^i}]\right)\\
&\equiv&
\theta_{x(t)}\left([t-T^{q^\nu}\otimes 1]^{-1}\prod_{i=0}^{\nu}
[t-1\otimes T^{q^i}]\right)
\bmod{\hat{\Delta}_0^{q^\nu+1}},
\end{array}$$
and this noted, (\ref{equation:SimClincher}) follows
by comparing Propositions~\ref{Proposition:DualJacobiTrudiRedux}
and \ref{Proposition:PeerIntoTheHole}
to formulas (\ref{equation:NegativeThetaCalc})
and (\ref{equation:NegativeCatalanCalc}). 
\qed

\begin{Lemma}\label{Lemma:UsefulColemanUnits}
Fix $x\in k^\times$ and an integer $n\geq 0$.
Put 
$$\Phi=\one_{x\tau^{-n+1}\OO}-\one_{\tau\OO}
-(q^{n}-1)
\one_{1+\tau\OO}.$$
Then there exists a Coleman unit $\varphi\in \kk_0$ such that
$$\lead_\Delta \varphi^{(\tau^{-N})}=\left(\begin{array}{c}
\tau^{-N}\\
\Phi
\end{array}\right)
$$ for all $N\gg 0$.
\end{Lemma}
\proof 
By example (\ref{equation:TrivialSolitonCalculation}) we may
assume without loss of generality that $x=1$. Now write $\Phi=\Phi_n$ to keep track of dependence on $n$.
We have
$$\Phi_0=0,\;\;\;\Phi_1=\one_{\OO}-\one_{\tau\OO}-(q-1)\one_{1+\tau\OO}=
\sum_{c\in \FF_q^\times}(\one_{1+\tau\OO}^{(c)}-\one_{1+\tau\OO}).$$
So the lemma holds trivially for $n=0$, and holds for $n=1$ by another
application of (\ref{equation:TrivialSolitonCalculation}).
For
all
$n\geq
0$
we
have
$$
\Phi_{n+2}-\Phi_{n+1}-q(\Phi_{n+1}-\Phi_n)=\Upsilon^{(\tau^{-n-1})},
$$
where
$$\Upsilon=\one_{\OO}-\one_{\tau\OO}-q(\one_{\tau\OO}-
\one_{\tau^{2}\OO}).$$
By induction on $n$ (twice) it suffices to find a Coleman unit
$\upsilon\in \kk_0$ such that
$$\lead_\Delta \upsilon^{(\tau^{-N})}=\left(\begin{array}{c}
\tau^{-N}\\
\Upsilon
\end{array}\right)$$
$$=\left( \frac{\Moore(T^N,\dots,1)}{
\Moore(T^{N-1},\dots,1)}\bigg/
\frac{\Moore(T^{N-1},\dots,1)^{q}}{\Moore(T^{N-2},\dots,1)^{q}}\right)^{q-1}=(T^{q^N}-T)^{q-1}$$
for $N\gg  0$. Clearly, 
$$\upsilon=(1\otimes T-T\otimes 1)^{q-1}=\varpi^{q-1}$$
has the desired properties.
\qed

\subsection{End of the proof}
Fix $x\in k^\times$. By Lemma~\ref{Lemma:StraightTheta} it remains only to prove that $\psi_x=(\theta_{x(t)}(P))^{(\tau^{-1})}$ is a Coleman unit.
 Fix monic $f\in \FF_q[T]$ 
and an integer $n\geq \deg f$. Note that $f^{-1}\tau^{-n+1}\OO\supset\tau\OO$ and hence $f^{-1}\tau^{-n+1}\OO$ is a finite union
of cosets of $\tau\OO$.
Consider the product
$$\varphi=\prod_{\xi\in k^\times\cap f^{-1}\tau^{-n+1}\OO}
\psi_\xi.$$
Since all the factors in the product by construction are nonzero elements of $\KK$,
and for suitable $f$ and $n$ the function $\psi_x$ is found among the factors
on the right, it suffices simply to show that $\varphi$ is a Coleman unit.
Consider the function
$$\Phi=\sum_{\xi\in k^\times\cap f^{-1}\tau^{-n+1}\OO}
\Psi_\xi=\one_{f^{-1}\tau^{-n+1}\OO}-\one_{\tau\OO}
-(q^{n}-1)
\one_{1+\tau\OO}.$$
By Lemma~\ref{Lemma:StraightTheta} we have $\varphi\sim_a\CCC[\Phi]$ for all $a\in \adeles^\times$. 
It follows by Lemma~\ref{Lemma:UsefulColemanUnits}
and the uniqueness principle (\ref{equation:UniquenessPrinciple})
that $\varphi$ is a Coleman unit. The proof of
Conjecture~\ref{TheConjecture} in the case
$g=0$ is complete.
\qed
\subsection{Remarks}
\subsubsection{}
The {\em solitons} defined in \cite[Thm.\ 2]{AndersonStickelberger},
and the {\em Coleman functions} applied in \cite{ABP}, \cite{Sinha} and \cite{SinhaDR} are all characterized by ``interpolation
formulas'' with righthand sides that can be put into
the form of the righthand side of equation (\ref{equation:PositiveCatalanCalc}) above. A similar remark applies to the formula (\ref{equation:ThakurInterpolation}) stated in the introduction. For example, by substituting $1+a/f$ for $x$ in the righthand side
of (\ref{equation:PositiveCatalanCalc}), one recovers the righthand
side of the interpolation formula stated in \cite[Thm.\ 2]{AndersonStickelberger}.
For another example, by substituting $x=\alpha/T+T+\beta$ in the righthand side of
(\ref{equation:PositiveCatalanCalc}), one recovers the righthand side of 
formula (\ref{equation:ThakurInterpolation}).
It follows via the uniqueness principle (\ref{equation:UniquenessPrinciple}) that  
solitons and Coleman functions belong to the class of two-variable algebraic functions 
engendered by Conjecture~\ref{TheConjecture} in the genus zero case. 
\subsubsection{}
What about the higher genus case?  We have four general methods
at our disposal for attacking the problem, namely the methods of ``rational Fourier analysis'' discussed in this paper, along with the methods discussed in the author's papers \cite{AndersonAHarmonic}, \cite{AndersonSpecialPoints},
and \cite{AndersonFock}. We are currently working with all these methods
in an effort to find a  proof of the conjecture.

\section{Horizontal specialization of Coleman units}
\label{section:HorizontalSpecialization}
We return to our studies at the full level of generality 
of Conjecture~\ref{TheConjecture}.
We study
embeddings of certain subrings of $\KK$ into certain power series
rings. Then we draw some general
conclusions concerning the operation of ``setting the second variable
equal to a constant'' in a Coleman unit, expressing the results
in the form of an easy-to-use calculus. Only at the very end
of this section do we draw any conclusions
conditional on Conjecture~\ref{TheConjecture}.

\subsection{Power series constructions}
Fix a field $\FF$.
\subsubsection{}
Let $t$, $t_1$ and $t_2$ be independent variables.  
Let $\FF((t))$ be the field obtained from the one-variable power series ring
$\FF[[t]]$ by inverting $t$. Let
$\FF((t_1,t_2))$ be the principal ideal domain obtained from the
two-variable power series ring
$\FF[[t_1,t_2]]$ by inverting $t_1$ and $t_2$. 
Given $F=F(t_1,t_2)\in \FF((t_1,t_2))$, put 
$$F^*=F^*(t_1,t_2)=F(t_2,t_1)\in \FF((t_1,t_2)),$$
and let the involutive automorphism $F\mapsto F^*$ of $\FF((t_1,t_2))$ be extended in the 
unique possible way to the fraction field of $\FF((t_1,t_2))$.
\subsubsection{}
Let $0\neq F=F(t_1,t_2)\in \FF((t_1,t_2))$ be given.
There is a unique factorization
$$F(t_1,t_2)=t_1^{\ell_1}t_2^{\ell_2}G(t_1,t_2)$$ where
$\ell_1,\ell_2\in \ZZ$ and
$G(t_1,t_2)\in \FF[[t_1,t_2]]$ is divisible by  neither $t_1$ nor $t_2$.
In turn, there is by the Weierstrass preparation theorem 
a unique factorization
$$G(t_1,t_2)=U(t_1,t_2)(t_1^m+t_2H(t_1,t_2))$$
where $U(t_1,t_2)\in \FF[[t_1,t_2]]^\times$, $m$ is a nonnegative integer
and $H(t_1,t_2)\in \FF[[t_1,t_2]]$ is a polynomial in $t_1$ of degree $<m$.
We now define 
$$\weight\, F=\ell_2\in \ZZ,\;\;\;\epsilon[F]=\epsilon[F](t)=t^{\ell_1+m}U(t,0)\in \FF((t))^\times.$$
Also put $\weight\, 0=+\infty$.
We have
\begin{equation}\label{equation:NumericalCriterion}
t^{-\weight\, F^*}\epsilon[F]\in \FF[[t]]^\times\Leftrightarrow m=0\Leftrightarrow
F\in \FF((t_1,t_2))^\times.
\end{equation}
The function $\weight$ is an additive valuation of $\FF((t_1,t_2))$. We extend $\weight$ in the unique possible way to a normalized
additive valuation of the fraction field of $\FF((t_1,t_2))$.
The $\FF((t))^\times$-valued function $\epsilon$
defined on $\FF((t_1,t_2))\setminus\{0\}$ is multiplicative. We extend $\epsilon$
in the unique possible way 
to a homomorphism from the multiplicative group of
the fraction field of $\FF((t_1,t_2))$ to $\FF((t))^\times$.

\subsection{Embeddings, valuations, and leading coefficients}
\label{subsection:Embeddings}
Fix a finite subextension $K/k$ of $k^\ab_\perf/k$.
Let $\FF$ be the field of constants of $K$. Put $K^{[2]}=K\otimes_{\FF}K$,
 which is a Dedekind domain. 

\subsubsection{}
Given a place $w$ of $K$
with residue field equal to $\FF$, a uniformizer $\pi$ in the completion $K_w$, and
$x\in K$, let
$x_{K,w,\pi}(t)\in
\FF((t))$ be the unique Laurent series such that the equality
$x_{K,w,\pi}(\pi)=x$ holds in $K_w$. 
Put $$\iota_{K,w,\pi}=(x\otimes y\mapsto x\cdot
y_{K,w,\pi}(t)):K^{[2]}\rightarrow K((t)),$$ 
thus defining an embedding of the ring $K^{[2]}$ into the 
field
$K((t))$.  We extend
$\iota_{K,w,\pi}$ in the unique   possible way to an embedding of
the fraction field of
$K^{[2]}$ into $K((t))$.

\subsubsection{}
With $w$ and $\pi$ as above, and given  nonzero $\varphi$ in the fraction field of
$K^{[2]}$, write
\begin{equation}\label{equation:PreWrite}
\iota_{K,w,\pi}\varphi=
\sum_{i=i_0}^\infty
e_it^i\;\;\;(i_0\in \ZZ,\;\;\;e_i\in K,\;\;\;e_{i_0}\neq 0),
\end{equation}
and put
$$\weight_{K,w,\pi}\varphi=i_0,\;\;\;
\epsilon_{K,w,\pi}\varphi=e_{i_0}.$$
 Also put $\weight_{K,w,\pi}0=+\infty$. 
The function 
$\weight_{K,w,\pi}$ is independent of $\pi$.
We write $\weight_{K,w}$ hereafter. 
The
function
$\weight_{K,w}$ is a normalized additive valuation of the fraction
field of
$K^{[2]}$. The function $\epsilon_{K,w,\pi}$ is a homomorphism from
the multiplicative group of the fraction field of $K^{[2]}$ to $K^\times$.
While $\epsilon_{K,w,\pi}$ does depend on $\pi$, the dependence
is rather mild: given another uniformizer
$\pi'\in K$ at
$w$, the ratio
$\epsilon_{K,w,\pi}/\epsilon_{K,w,\pi'}$ is a homomorphism
from the fraction field of $K^{[2]}$ to
$\FF^\times$ of the form $F\mapsto c^{\weight_{K,w}\,F}$
for some $c\in \FF^\times$.
\subsubsection{}
Let $w$, $\pi$ and $\varphi$ be as above,
but assume now that  $\varphi$ belongs to the fraction field of $K^{[2]}\cap\KK$.
Fix $a\in \adeles^\times$ such that $\norm{a}\geq 1$.
Then $\varphi^{(a)}$ is defined and belongs again to
the fraction field of $K^{[2]}\cap \KK$. We claim that
\begin{equation}\label{equation:TwistFunctoriality}
\weight_{K,w}\varphi^{(a)}=\norm{a}\weight_{K,w}\varphi,\;\;
\epsilon_{K,w,\pi}(\varphi^{(a)})=(\epsilon_{K,w,\pi}\varphi)^{\rho(a)}.
\end{equation}
In any case, notation as in (\ref{equation:PreWrite}), we must have
$$\iota_{K,w,\pi}\varphi^{(a)}=
\sum_{i=i_0}^\infty e_i^{\rho(a)}t^{i\norm{a}}$$
on account of the rule $x\otimes y\mapsto x^{\rho(a)}\otimes y^{\norm{a}}$
by which $\varphi^{(a)}$ is defined. The claim follows.

\subsubsection{}
 Given for $i=1,2$ a place $w_i$ of $K$ with residue field $\FF$
and a uniformizer $\pi_i$ in $K_{w_i}$, put
$$\iota_{K,w_1,w_2,\pi_1,\pi_2}=
(x\otimes y\mapsto x_{K,w_1,\pi_1}(t_1)
y_{K,w_2,\pi_2}(t_2)):K^{[2]}\rightarrow
\FF((t_1,t_2)),$$
thus defining an embedding of the ring $K^{[2]}$ into $\FF((t_1,t_2))$.
We extend 
$\iota_{K,w_1,w_2,\pi_1,\pi_2}$ in the unique
 possible way to an embedding of the fraction field of $K^{[2]}$
into the fraction field of $\FF((t_1,t_2))$.
Clearly, we have
\begin{equation}\label{equation:PowerSeriesInvariants0}
\weight_{K,w_2}\varphi=\weight\, \iota_{K,w_1,w_2,\pi_1,\pi_2}\varphi.
\end{equation}
It  is not difficult to verify that the equality
\begin{equation}\label{equation:PowerSeriesInvariants1}
\epsilon_{K,w_2,\pi_2}\varphi=\epsilon[\iota_{K,w_1,w_2,\pi_1,\pi_2} \varphi](\pi_1)
\end{equation}
holds in $K_{w_1}$. If $\varphi$ belongs to the fraction field of $\KK\cap K^{[2]}$,
then $\varphi^\dagger$ is defined and belongs again to the fraction field
of $\KK\cap K^{[2]}$, and
we have
\begin{equation}\label{equation:PowerSeriesInvariants2}
(\iota_{K,w_1,w_2,\pi_1,\pi_2}\varphi)^*=\iota_{K,w_2,w_1,\pi_2,\pi_1}(\varphi^\dagger),
\end{equation}
on account of the rule $x\otimes y\mapsto y\otimes x$ defining the ``dagger'' operation.

\begin{Lemma} Let $K$, $\FF$
and $K^{[2]}$ be as in \S\ref{subsection:Embeddings}. For $i=1,2$ let $w_i$ be a place
of $K$ with residue field equal to $\FF$, let $\pi_i$ be a uniformizer in $K_{w_i}$,
and let $v_i$ be the unique place of $k$ below $w_i$. Assume that $v_1\neq v_2$.
Fix $a\in \adeles^\times$ and put 
$$I(a)=\ker\left((x\otimes y\mapsto
x^{\rho(a)}y^{\norm{a}}):K^{[2]}\rightarrow K_\perf\right).$$ 
Then there exists $\psi\in I(a)$ such that
$$\iota_{K,w_1,w_2,\pi_1,\pi_2}\psi\in
\FF[[t_1,t_2]]^\times.$$
\end{Lemma}
\proof
For any $x\in K$ we have
$$\delta_a x=(x^{\rho(a)}\otimes 1-1\otimes
x^{\norm{a}})^{\max(1,\norm{a}^{-1})}
\in I(a).$$
Since the additive valuations
$$x\mapsto \ord_{w_1} x^{\rho(a)},\;\;\;
x\mapsto \ord_{w_2}x$$
are inequivalent (they already have inequivalent
restrictions to $k$ by hypothesis), we can find $x\in K$ by the
Artin-Whaples approximation theorem such that
$$\ord_{w_1}x^{\rho(a)}>0,\;\;\;\ord_{w_2}x=0.$$
Then $\psi=\delta_ax$ has the desired property.
\qed

\begin{Proposition}\label{Proposition:OffDiagonal} Let $K$, $\FF$
and $K^{[2]}$ be as in \S\ref{subsection:Embeddings}.
For $i=1,2$ let $w_i$ be a place
of $K$ with residue field equal to $\FF$, let $\pi_i$ be a uniformizer in $K_{w_i}$,
and let $v_i$ be the unique place of $k$ below $w_i$.  Assume that $v_1\neq v_2$.
For all Coleman
units $\varphi$ belonging to the fraction field of $K^{[2]}$,
we have $\ord_{w_1}\epsilon_{K,w_2,\pi_2}\varphi=\weight_{K,w_1}
\varphi^\dagger$.
\end{Proposition}
\proof By the preceding lemma and the definition of a Coleman
unit, we can find
$\psi\in K^{[2]}$ such that 
$$\psi\varphi^{\pm 1}\in
K^{[2]},\;\;\iota_{K,w_1,w_2,\pi_1,\pi_2}(\psi)\in
\FF[[t_1,t_2]]^\times,$$
in which case it follows that
$$\iota_{K,w_1,w_2,\pi_1,\pi_2}(\varphi)\in \FF((t_1,t_2))^\times.$$
We then have
$$ \begin{array}{rcl}
\ord_{w_1}\epsilon_{K,w_2,\pi_2}&=&\ord_{w_1}\epsilon[\iota_{K,w_2,\pi_2}\varphi](\pi_1)\\
&=&\weight\,\left(\iota_{K,w_1,w_2,\pi_1,\pi_2}\varphi\right)^*\\
&=&\weight\, \iota_{K,w_2,w_1,\pi_2,\pi_1}(\varphi^\dagger)\\
&=&\weight_{K,w_1,\pi_1}\varphi^\dagger
\end{array}
$$
at the first equality by (\ref{equation:PowerSeriesInvariants1}), at the second by 
(\ref{equation:NumericalCriterion}),
at the third by (\ref{equation:PowerSeriesInvariants2}),
and at the fourth by (\ref{equation:PowerSeriesInvariants0}).
 \qed

\begin{Proposition}\label{Proposition:GrossKoblitzRedux} Let $K$, $\FF$
and $K^{[2]}$ be as in \S\ref{subsection:Embeddings}.  Let $w$ be a place of $K$ with residue field equal to $\FF$ and let $\pi$ be a uniformizer in $K_w$. Let $v$ be the place of $k$ below $w$ and let $\tau$ be a uniformizer in $k_v$. For every nonzero $\varphi$ belonging to the fraction field of
$K^{[2]}\cap\KK$, we
have
$$\lim_{n\rightarrow\infty}\ord_{w}\left(
\epsilon_{K,w,\pi}\varphi-\pi^{-q_v^{nd}\weight_{K,w}\varphi}\lead_\Delta
\varphi^{(i_v^\times(\tau^{-nd}))}\right)=+\infty,
$$
where $d=[K:k]$, and $\ord_w$ is the normalized additive valuation of $K$
associated to $w$.
\end{Proposition}
\noindent
We are indebted to D.\ Thakur for teaching us about this sort of
convergence phenomenon. (See, for example, \cite[Thm.\ 4.8.1]{ThakurBook}.)

\proof In general, $\varphi$ is a quotient of elements of $K^{[2]}\cap\KK$.
But since we get the result if we can prove it for numerator and
denominator of $\varphi$ separately, we may assume without loss of generality
that
$\varphi\in
K^{[2]}\cap\KK$. Put 
$$\iota_{K,w,w,\pi,\pi}(\varphi)=F(t_1,t_2)\in \FF((t_1,t_2)).$$
Write
$$F(t_1,t_2)=t_1^{\ell_1}t_2^{\ell_2}G(t_1,t_2),$$ where
$G(t_1,t_2)\in \FF[[t_1,t_2]]$ is divisible neither by $t_1$ nor $t_2$.
By the Weierstrass preparation theorem, write $$G(t_1,t_2)=(t_1^m+t_2H(t_1,t_2))U(t_1,t_2),$$
where $m$ is a nonnegative integer,
$H(t_1,t_2)\in \FF[[t_1,t_2]]$ is a polynomial in $t_1$ of degree $<m$, and $U(t_1,t_2)\in \FF[[t_1,t_2]]^\times$.
We have
$$\weight_{K,w}\varphi=\ell_2,\;\;\;\epsilon_{K,w,\pi}\varphi=\pi^{\ell_1}U(\pi,0)$$
by (\ref{equation:PowerSeriesInvariants0}),
(\ref{equation:PowerSeriesInvariants1}) and the definitions.
Write
$$\varphi=\sum_i x_i\otimes y_i\;\;\;(x_i,y_i\in K).$$
Note that for $n\geq 0$ we have
$$\varphi^{(i^\times_v(\tau^{-nd}))}=\sum_i
x_i^{\rho(i^\times_v(\tau^{-n}))^d}\otimes y_i^{q_v^{nd}}=\sum_i
x_i\otimes y_i^{q_v^{nd}},$$
also we have $\card \FF\vert q_v^{d}$, and hence
equality
$$\lead_\Delta \varphi^{(i^\times_v(\tau^{-nd}))}=
\sum_i x_iy_i^{q_v^{nd}}=
F(\pi,\pi^{q_v^{nd}})$$
$$=\pi^{\ell_1+q_v^{nd}\ell_2}(\pi^m+\pi^{q_v^{nd}}H(\pi,\pi^{q_v^{nd}}))
U(\pi,\pi^{q_v^{nd}})$$
holds in the completion $K_w$. 
The result follows.
\qed

\subsection{The functions $\weight_{v}$ and $\epsilon_{\bar{v}}$}
\label{subsection:WeightEpsilonCalculus}
We boil the preceding somewhat complicated considerations
down  to a few easy-to-apply rules.

\subsubsection{} Let $v$ be a place of $k$,
let $\bar{v}$ be a place of $k^\ab_\perf$ above $v$
and let $\ord_{\bar{v}}$ be the unique additive valuation of $k^\ab_\perf$
belonging to the place $\bar{v}$ which extends
the normalized additive valuation $\ord_v$ of $k$. For every nonzero
$\varphi$ in the fraction field of $\KK$ we define
$$\begin{array}{rcl}
\weight_{v}\varphi&=&
\ord_{\bar{v}}\pi\cdot\weight_{K,w}\varphi,\\\\
\epsilon_{\bar{v}}\varphi&=&
\left\{\begin{array}{ll}
\epsilon_{K,w,\pi}\varphi&\mbox{if $\weight_{v}\varphi=0$,}\\
\epsilon_{K,w,\pi}\varphi\bmod{(\FF^\ab_q)^\times}&
\mbox{if $\weight_{v}\varphi\neq 0$,}
\end{array}\right.
\end{array}
$$
where:
\begin{itemize}
\item $K/k$ is any finite subextension of $k^\ab_\perf/k$ such that
\begin{itemize}
\item $\varphi$ belongs to the fraction field of $K^{[2]}$, 
\item $w$ is the place of $K$ below $\bar{v}$,
\item $\pi$ is a uniformizer at $w$, and
\item the constant field $\FF$ of $K$ equals
the residue field of $w$.
\end{itemize}
\end{itemize}
We are obliged to check that $\weight_v$ and $\epsilon_{\bar{v}}$ are well-defined.
In any case, for every $\varphi$ we can find suitable $K$ as above.
Moreover, for fixed
$\varphi$, the expression
on the right side of the definition of $\epsilon_{\bar{v}}$ is easily
verified to depend only on $\bar{v}$, and
the same is true of the expression on the right side of the
definition of $\weight_v
\varphi$. Finally, the latter depends only
on
$v$ by Proposition~\ref{Proposition:GrossKoblitzRedux}.
Thus
$\weight_v$ and $\epsilon_{\bar{v}}$ are indeed well-defined.
For convenience put $\weight_{v}0=+\infty$. Then the function
$\weight_v$ is an additive valuation of the fraction field
of $\KK$. The function $\epsilon_{\bar{v}}$ is a homomorphism
from the multiplicative group of the fraction field of $\KK$
to $(k^\ab_\perf)^\times/(\FF_q^\ab)^\times$ 
which on the subgroup $\{\weight_{\bar{v}}=0\}$
is refined to a homomorphism to
$(k^\ab_\perf)^\times$. 

\subsubsection{}
We claim that the following relations hold:
\begin{equation}\label{equation:TwistTransformationBis}
\weight_{v}\varphi^{(a)}=\norm{a}\weight_{v}\varphi,\;\;\;\;
\epsilon_{\bar{v}}\varphi^{(a)}
=(\epsilon_{\bar{v}}\varphi)^{\rho(a)}
\end{equation}
\begin{equation}\label{equation:GrossKoblitzReduxBis}
\ord_v\lead_\Delta \varphi^{(i_v^\times(\tau^{-n!}))}
=\ord_{\bar{v}}\epsilon_{\bar{v}}\varphi+
q_v^{n!}\weight_v\varphi\;\;\;(n\gg 0)
\end{equation}
\begin{equation}\label{equation:GrossKoblitzReduxTer}
\weight_v\varphi=0\Rightarrow
\lim_{n\rightarrow\infty}
\ord_{\bar{v}}\left(\epsilon_{\bar{v}}\varphi-\lead_\Delta
\varphi^{(i_v^\times(\tau^{-n!}))}\right)=+\infty
\end{equation}
\begin{equation}
\label{equation:OffDiagonal}
\ord_{\bar{v}_1}\epsilon_{\bar{v}_2}\varphi
=\weight_{v_1}\varphi^\dagger\;\;\;(v_1\neq v_2)
\end{equation}
Here $\bar{v}$ (resp., $\bar{v}_1$, $\bar{v}_2$)
are places of $k^\ab_\perf$ above places $v$ (resp., $v_1$, $v_2$) of
$k$,
$\varphi$ is a Coleman unit, $a\in \adeles^\times$, and $\tau\in k$
is a uniformizer at $v$. Relation
(\ref{equation:TwistFunctoriality}) justifies
(\ref{equation:TwistTransformationBis}).
Proposition~\ref{Proposition:GrossKoblitzRedux} justifies relations
(\ref{equation:GrossKoblitzReduxBis}) and
(\ref{equation:GrossKoblitzReduxTer}).
Proposition~\ref{Proposition:OffDiagonal} justifies
(\ref{equation:OffDiagonal}).

\subsubsection{} Keeping the notation of the preceding
paragraph, suppose further that $\varphi=\CCC[\Phi]$ for some
$\Phi\in
\sch_{00}(\adeles)$. (Notice that not until now have we invoked
Conjecture~\ref{TheConjecture}.) We claim
that
\begin{equation}\label{equation:WeightIntegral}
\weight_v\varphi=\int \tilde{\Phi}(i_v(t))d\mu_v^\times(t),
\end{equation}
\begin{equation}\label{equation:DaggerWeightIntegral}
\weight_v\varphi^\dagger=\int\Phi(i_v(t))d\mu_v^\times(t),
\end{equation}
\begin{equation}\label{equation:barv2}
\ord_{\bar{v}}\epsilon_{\bar{v}}\varphi=
\int\Phi(i_v(t))d\mu_v^\times(t)+\int\Theta(i_v^\times(t),\Phi)d\mu_v^\times(t).
\end{equation}
Equations (\ref{equation:WeightIntegral})
and (\ref{equation:barv2}) are proved by comparing
(\ref{equation:GrossKoblitzReduxBis})
to the simplified version (\ref{equation:AdelicStirlingBis})
of the adelic Stirling formula. Equation
(\ref{equation:DaggerWeightIntegral}) follows from
(\ref{equation:WeightIntegral}) since, as remarked in
\S\ref{subsection:CCalculus}, we
have
$\CCC[\tilde{\Phi}]=\CCC[\Phi]^\dagger$. 

\section{A conditional recipe for the Stark unit}
\label{section:StarkRecovery}
We continue in the setting of Conjecture~\ref{TheConjecture}.
Let $K/k$ be a finite subextension of $k^\ab/k$.
Let $S$ be a finite set of places of $k$. Assume the following:
\begin{itemize}
\item Some place $\infty\in S$ splits completely in $K$.
\item   The set $S_0=S\setminus\{\infty\}$ is nonempty. 
\item All places of $k$ ramified in $K$ belong
to
$S_0$.
\end{itemize}
We are going to show that Tate's formulation 
$\St(K/k,S)$ \cite[p.\ 89, Conj.\
2.2]{Tate} of the Stark conjecture  is a consequence of
Conjecture~\ref{TheConjecture}. After introducing suitable notation and making a convenient reduction, we recall the statement
of $\St(K/k,S)$ in detail below.
We must of course
remark that since we are working in the function field situation,
$\St(K/k,S)$ is already a theorem due to Deligne 
\cite{Tate} and (independently) to Hayes \cite{Hayes}.
 The point of deriving
$\St(K/k,S)$ from our conjecture is to establish that
the latter does in fact refine the former.

\subsection{Notation and a reduction}
\label{subsection:StarkNotation}
\subsubsection{}
For each place $v$ of $k$ we fix the following objects:
\begin{itemize}
\item
Let
$\tau_v\in \adeles^\times$ be  the image of a fixed choice of uniformizer of $k_v$
under the map
$i_v^\times:k_v^\times\rightarrow\adeles^\times$.
\item Put $\langle
\tau_v\rangle=\{\tau_v^n\vert n\in \ZZ\}\subset \adeles^\times$.
\item Fix a place $\bar{v}$ of $k^\ab_\perf$ above $v$.
\item Let $\ord_{\bar{v}}$ be the unique additive valuation of
$k^\ab_\perf$ belonging to $\bar{v}$ and extending the normalized
additive valuation $\ord_v$ of
$k$.
\end{itemize}
\subsubsection{}
According to Tate \cite[Prop.\ 3.5, p.\ 92]{Tate},
$\St(K/k,S)$ implies
$\St(K'/k,S)$ for any subextension
$K'/k$ of $K/k$.
Accordingly, after choosing a suitable finite subextension $\tilde{K}/K$
of $k^\ab/K$ and replacing $(K/k,S)$ by $(\tilde{K}/k,S)$, we may
assume without loss of generality that the data $(K/k,S)$ satisfy the
following further condition:
\begin{itemize}
\item $\ker\rho_{K/k}=k^\times\UUU\langle
\tau_\infty\rangle\subset\adeles^\times$, where $\UUU\subset \OO^\times$ is
an open subgroup with the following properties:
\begin{itemize}
\item $\UUU\supset i_v^\times (\OO_v^\times)$ for all places $v$ of $k$ not belonging to $S_0$.
\item $\UUU\cap k^\times=\{1\}$.
\end{itemize}
\end{itemize}
Note that under this further condition the constant field
$\FF$ of $K$ has cardinality $q_\infty$,
and hence there are exactly $q_\infty-1$ roots of unity in $K$.
\subsubsection{}
For $T=S_0$ or $T=S$, consider the Euler product
$$\theta_{T}(s)=\prod_{v\not\in T}
(1-F_vq_v^{-s})^{-1}=\sum_{\sigma\in G}\zeta_T(s,\sigma)\sigma^{-1}\in
\CC[G]\;\;\;(\Re(s)>1)$$
extended over places $v$ of $k$ not in $T$, where
$F_v=\rho_{K/k}(\tau_v)\in G$ is the geometric Frobenius at
$v$ (cf.\ \cite[p.\ 86, Prop.\ 1.6]{Tate}). 
If $T=S_0$, we drop the subscript and write
simply $\theta(s)$ and $\zeta(s,\sigma)$. It is
well-known that
$\theta_T(s)$ continues meromorphically to the entire $s$-plane with
no singularity other than a pole at
$s=1$. Note that
we have 
$$\theta_S(0)=0,\;\;\;\theta'_{S}(0)=
\log q_\infty\cdot \theta(0)$$
(cf.\ \cite[p.\ 86, Cor.\ 1.7]{Tate}) 
since $\infty$ splits completely in $K/k$.
\subsubsection{}
Let $U\subset K^\times$ be the $G$-submodule
consisting of $x$ satisfying the following condition:
\begin{itemize}
\item $\ord_{\bar{v}}x^\sigma=
\left\{\begin{array}{rl}
\ord_{\bar{v}}x&\mbox{if $S_0=\{v\}$,}\\
0&\mbox{otherwise,}
\end{array}\right.$ for all places $v$ of $k$ distinct from $\infty$ and
$\sigma\in G$.
\end{itemize}
Let
$U^\ab\subset U$ be the $G$-submodule
consisting of
$x$ satisfying the following further condition:
\begin{itemize}
\item $K(x^{1/(q_\infty-1)})/k$ is an abelian extension.
\end{itemize}
It is convenient to define a group homomorphism
$$\ordbold=\left(x\mapsto
\sum_{\sigma\in
G}(\ord_{\overline{\infty}}x^\sigma)\sigma^{-1}\right):U\rightarrow
\CC[G].$$
Note that $\ordbold$ is $G$-equivariant
and $\ker\ordbold=\FF^\times$.
\subsubsection{}\label{subsubsection:StarkUnitSpecs}
The
{\em Stark unit}
$$\epsilon(K/k,S)\in U^\ab$$ predicted to
exist by conjecture 
$$\St(K/k,S)\;\;\;(\mbox{\cite[p.\ 89, Conj.\
2.2]{Tate}}),$$
given the choice $\overline{\infty}\vert_{K}$ of a place of $K$
above
$\infty$,
is uniquely determined up to a factor in
$\FF^\times$  by the formulas
$$(q_\infty-1)\zeta(0,\sigma)=\ord_{\overline{\infty}}\,\epsilon(K/k,S)^\sigma\;\;
\mbox{for all $\sigma\in G$},$$
or, equivalently, the formula
\begin{equation}\label{equation:StarkUnit}
(q_\infty-1)\theta(0)=\ordbold\,\epsilon(K/k,S).
\end{equation}
The rest of our work in \S\ref{section:StarkRecovery} is devoted
to working out a ``recipe'' for the Stark unit in terms of the
transformation
$\CCC$ defined by Conjecture~\ref{TheConjecture}.

\subsection{``Ingredient list''}
\label{subsection:Recipe}

\subsubsection{}
Let $\overline{S}$ be the set of places of $k$ not belonging to $S$.
Let $\Gamma$ be the subgroup of
$\adeles^\times$ generated by $\{\tau_v\vert v\in \overline{S}\}$.  Let $\ZZ[\Gamma]$ be the group ring of $\Gamma$ over the integers. Let $\JJ$ be the kernel of the ring homomorphism
$\ZZ[\Gamma]\rightarrow \ZZ[1/q]$ induced by the restriction
of the idele norm function $\norm{\cdot}$ to $\Gamma$.
Since $\Gamma$ is a free abelian group with basis $\{\tau_v\vert
v\in \overline{S}\}$, the ring $\ZZ[\Gamma]$ may be viewed as the ring of
Laurent polynomials with integral coefficients in independent variables
$\tau_v$ indexed by $v\in \overline{S}$.
From the latter point of view it is clear that
the set 
$$\{1-q_v\cdot \tau_v\vert v\in \overline{S}\}\subset \ZZ[\Gamma]$$
generates $\JJ$ as an ideal of
$\ZZ[\Gamma]$.  Here, and in similar contexts below,
the expression $1-q_v\cdot
\tau_v$ is to be viewed as a formal
$\ZZ$-linear combination of elements of $\Gamma$.

\subsubsection{}
Put
$$\UUU(\infty)=\{a=[a_v]\in \UUU\vert \norm{a_\infty-1}_\infty<1\}.$$
Let $K(\infty)/k$ be the unique subextension of $k^\ab/k$ such that
$$\ker\rho_{K(\infty)/k}=k^\times \UUU(\infty)\langle \tau_\infty\rangle$$
and put
$$G(\infty)=\Gal(K(\infty)/k).$$
Crucially, the constant field $\FF$ of $K$ is also the constant
field of $K(\infty)$.

\subsubsection{}
Put
$$\begin{array}{rclcl}
\Pi&=&\Gamma\cap
k^\times\UUU\langle
\tau_\infty\rangle,\\
\Pi(\infty)&=&\Gamma\cap
k^\times\UUU(\infty)\langle
\tau_\infty\rangle,\\\Pi_1(\infty)&=&\Gamma\cap
k^\times\UUU(\infty)=\{a\in
\Pi(\infty)\vert\;
\norm{a}=1\}.
\end{array}$$
Extend the group homomorphism
$\Gamma\xrightarrow{\rho_{K/k}}G$
to a ring homomorphism $\ZZ[\Gamma]\xrightarrow{\rho_{K/k}}
\ZZ[G]$, and define
$\ZZ[\Gamma]\xrightarrow{\rho_{K(\infty)/k}}\ZZ[G(\infty)]$ analogously.
Let $J(\infty)\subset\ZZ[G(\infty)]$ be the annihilator of $\FF^\times$.
For any subgroup $\Gamma'\subset\Gamma$ let
$\II(\Gamma')\subset\ZZ[\Gamma]$ be the ideal generated by differences
of elements of $\Gamma'$.  From the exactness of the sequence
$$1\rightarrow \Pi\rightarrow
\Gamma\xrightarrow{\rho_{K/k}}G\rightarrow 1,$$
which is well-known, we deduce an exact sequence
\begin{equation}\label{equation:CurlyJExact0}
0\rightarrow \II(\Pi)\subset
\ZZ[\Gamma]\xrightarrow{\rho_{K/k}}\ZZ[G]\rightarrow 0.
\end{equation}
Analogously we have an exact sequence
\begin{equation}\label{equation:CurlyJExact1}
0\rightarrow \II(\Pi(\infty))\subset
\ZZ[\Gamma]\xrightarrow{\rho_{K(\infty)/k}}\ZZ[G(\infty)]\rightarrow 0.
\end{equation}
We claim that we have an exact sequence
\begin{equation}\label{equation:CurlyJExact2}
0\rightarrow \II(\Pi_1(\infty))+\II(\Pi(\infty))
\cdot \JJ\subset \JJ\xrightarrow{\rho_{K(\infty)/k}}
J(\infty)\rightarrow 0.
\end{equation}
Exactness at $\JJ$ follows from Lemma~\ref{Lemma:IdealTheoryInGroupRings}.
The set
$$\{v\in \overline{S}\vert
1-q_v\rho_{K(\infty)/k}(\tau_v)\}$$
by \cite[p.\ 82, Lemme 1.1]{Tate} 
generates $J(\infty)$ as an ideal of $\ZZ[G(\infty)]$,
 whence exactness at $J(\infty)$.
The claim is proved. 

\subsubsection{}
From exactness of (\ref{equation:CurlyJExact1})
and (\ref{equation:CurlyJExact2}) it follows that there exists
$\abold_\infty\in \JJ$ such that
$\rho_{K(\infty)/k}(\abold_\infty)=q_\infty-1$. 

\subsubsection{}
Put
$$\VVV=\left\{a=[a_v]\in \OO\left|v\in S_0\Rightarrow a_v\neq 0,\;\;\;\prod_{v\in
S_0}i_v^{\times}(a_v)\in
\UUU\right.\right\}\subset \OO,$$
which is an open compact subset of $\adeles$. Note that
$\VVV$ is stable under the action of the group $\UUU$
by multiplication.
Extend the map
$$(a\mapsto \one_\VVV^{(a)}):\Gamma\rightarrow
\sch_{0}(\adeles)$$ 
$\ZZ$-linearly to a homomorphism
$$(\abold\mapsto \one_\VVV^{(\abold)}):\ZZ[\Gamma]\rightarrow
\sch_0(\adeles)$$
of abelian groups,
noting that
\begin{equation}\label{equation:PreliminaryAbold}
\one_\VVV^{(a\abold)}=(\one_\VVV^{(\abold)})^{(a)},\;\;\;\;\abold\in
\JJ\Leftrightarrow  \one_\VVV^{(\abold)}\in
\sch_{00}(\adeles)
\end{equation}
for all $\abold\in \ZZ[\Gamma]$ and $a\in \Gamma$.

\subsubsection{} We claim that for all
$\abold\in\JJ$ we have
\begin{equation}\label{equation:WeightInfinity}
\int
\FFF_0[\one_\VVV^{(\abold)}](i_\infty(t))d\mu_\infty^\times(t)=0.
\end{equation}
In any case, the set of $\abold\in \JJ$ satisfying the equation above
forms an ideal of $\ZZ[\Gamma]$, and so to prove the claim we may assume
without loss of generality that
$\abold=1-q_v\cdot \tau_v$ for some $v\in \overline{S}$, in which case
already by the scaling rule (\ref{equation:AdelicRationalScaling})
the integrand above vanishes identically as a function
of $t\in k_\infty^\times$.
The claim is proved.

\subsubsection{}
So far all our constructions and definitions make sense unconditionally.
But from this point onward we must assume that
Conjecture~\ref{TheConjecture} holds so that the transformation
$\CCC$ is defined.
From (\ref{equation:WeightIntegral}) and (\ref{equation:WeightInfinity}) 
 it follows that for every $\abold\in \JJ$ we have
$$\weight_\infty \CCC[\one_\VVV^{(\abold)}]=0$$
and hence
$$\epsilon(\abold)=\epsilon_{\overline{\infty}}\CCC[\one_\VVV^{(\abold)}]
\in (k^\ab_\perf)^\times$$
is well-defined. 
Note that by
(\ref{equation:TwistTransformationBis}) and
(\ref{equation:PreliminaryAbold}) the homomorphism
$$\epsilon:\JJ\rightarrow (k^\ab_\perf)^\times$$
is $\Gamma$-equivariant, where we view the target in the natural way
as a
$\Gamma$-module via the reciprocity law $\rho$.

\begin{Theorem}[``Recipe'']\label{Theorem:Recipe} Hypotheses and
notation as above, we have 
$\epsilon(K/k,S)=
\epsilon(\abold_\infty)$.
\end{Theorem}
\noindent We stress that this result is conditional on
Conjecture~\ref{TheConjecture}. The proof consists of an analysis of 
the
$\Gamma$-equivariant homomorphism
$\epsilon$ proceeding by way of several lemmas.

\begin{Lemma}
$\epsilon(\JJ)\subset 
\left(\mbox{$\bar{\infty}$-closure of $k$ in $k^\ab_\perf$}\right)^\times
\subset (k^\ab)^\times$.
\end{Lemma}
\proof For all
$\abold\in
\JJ$ we have
$$\lim_{n\rightarrow\infty}
\ord_{\bar{\infty}}\left(\epsilon(\abold)-
\left(\begin{array}{c}
\tau_\infty^{-n!}\\\\
\one_\VVV^{(\abold)}
\end{array}\right)\right)=+\infty
$$
by limit formula (\ref{equation:GrossKoblitzReduxTer}) and the
definitions. Moreover, the Catalan symbol in question
takes values in $k^\times$ for all $n\gg 0$
by Proposition~\ref{Proposition:GenericCatalanBehavior}. The result follows.
\qed

\begin{Lemma}
$\ker\epsilon\supset\II(\Pi)\cdot \JJ+
\II(\Pi_1(\infty))$.
\end{Lemma}
\proof Fix $a\in \Pi$ arbitrarily and write
$a=ux\tau_\infty^N$ with $u\in \UUU$, $x\in
k^\times$ and $N\in \ZZ$, noting that since $\UUU\cap k^\times=\{1\}$,
this factorization is unique.
 For
all
$\abold\in
\JJ$ we have
$$\begin{array}{rcl}
\epsilon(\abold
a)=\epsilon_{\overline{\infty}}\,\CCC[\one_\VVV^{(\abold
a)}]&=&\epsilon_{\overline{\infty}}\,\CCC[\one_\VVV^{(\abold
x\tau_\infty^n)}]\\
&=&\epsilon_{\overline{\infty}}\,\CCC[(\one_\VVV^{(\abold
)})^{(x\tau_\infty^n)}]\\
&=&\epsilon_{\overline{\infty}}\,(\CCC[\one_\VVV^{(\abold
)}]^{(\tau_\infty^n)})=\epsilon(\abold)^{\rho(\tau_\infty^n)}=\epsilon(\abold)
\end{array}$$
where the fourth equality is justified by the
$\adeles^\times$-equivariance of $\CCC$ noted in
\S\ref{subsection:CCalculus}, the fifth  by (\ref{equation:TwistTransformationBis}), and the last by the preceding lemma. Therefore we have $\ker\epsilon\supset\II(\Pi)\cdot
\JJ$.  Now suppose that $a\in \Pi_1(\infty)$. Then necessarily
$N=0$ and $xu_\infty=1$, hence $\norm{x-1}_\infty<1$, and hence
$$\epsilon(\one_\VVV^{(a)}-\one_\VVV)=
\epsilon_{\overline{\infty}}\,\CCC[
\one_\VVV^{(x)}-\one_\VVV]=\epsilon_{\overline{\infty}}(1\otimes
x^{\mu(\VVV)})=1,$$
where the middle equality is justified by example
(\ref{equation:TrivialSolitonCalculation}). Therefore we have
$\ker\epsilon\supset\II(\Pi_1(\infty))$. 
\qed
\begin{Lemma}
 $\epsilon(\JJ)\subset U^\ab$.
\end{Lemma}
\proof By exactness of the sequence (\ref{equation:CurlyJExact0})
and the
previous two lemmas,
$\epsilon$ takes values in
$K^\times$. We claim that
\begin{equation}\label{equation:WeightOfSpecialEpsilon}
\int
\one_\VVV^{(\abold a)}(i_v(t))d\mu_v^\times(t)
=\left\{\begin{array}{rl}
\displaystyle\int
\one_\VVV^{(\abold)}(i_v(t))d\mu_v^\times(t)
&\mbox{if
$S_0=\{v\}$}\\ 0&\mbox{if
$S_0\neq\{v\}$}
\end{array}\right.
\end{equation}
for all $\abold\in \JJ$, $a\in \Gamma$ and places $v$ of $k$.
If $S_0\neq \{v\}$, then we have
$b\VVV\cap i_v(k_v)=\emptyset$ for all $b\in \adeles^\times$
and {\em a fortiori} (\ref{equation:WeightOfSpecialEpsilon}) holds;
otherwise, if
$S_0=\{v\}$, we get (\ref{equation:WeightOfSpecialEpsilon}) by an evident
manipulation of integrals. The claim is proved. It follows by formula
(\ref{equation:DaggerWeightIntegral}) of
\S\ref{subsection:WeightEpsilonCalculus} that 
$$\weight_v
\CCC[\one_\VVV^{(\abold a)}]^\dagger=\left\{\begin{array}{rl}
\displaystyle\int
\one_\VVV^{(\abold)}(i_v(t))d\mu_v^\times(t)
&\mbox{if
$S_0=\{v\}$}\\ 0&\mbox{if
$S_0\neq\{v\}$}
\end{array}\right.$$
for all $\abold\in \JJ$, $a\in \Gamma$ and places $v$ of $k$.
In turn it follows by formula (\ref{equation:OffDiagonal}) of
\S\ref{subsection:WeightEpsilonCalculus} and the $\Gamma$-equivariance
of $\epsilon$ that
$\epsilon(\abold)\in U$ for all $\abold\in \JJ$. Finally, since by the
preceding lemma and exactness of sequences (\ref{equation:CurlyJExact1})
and (\ref{equation:CurlyJExact2}) we may view $\epsilon$
as a $G(\infty)$-equivariant homomorphism $J(\infty)\rightarrow U$,
in fact $\epsilon$ takes values in $U^\ab$ by
Lemma~\ref{Lemma:AbelianRootRedux}. \qed
\begin{Lemma}
$\rho_{K/k}(\abold)\theta(0)=\ordbold \,\epsilon(\abold)$
for all $\abold\in \JJ$.
\end{Lemma}
\proof For all $\abold\in \JJ$, by (\ref{equation:barv2}) and
(\ref{equation:WeightOfSpecialEpsilon}) in the case $v=\infty$, we have
$$
\ord_{\overline{\infty}}\epsilon(\abold)=
\int\Theta(i_\infty^\times(t),\one_{\VVV}^{(\abold)})d\mu_\infty^\times(t)=\sum_{n\in\ZZ}\Theta(    \tau_\infty^n
,\one_\VVV^{(\abold)}),
$$
 and so with $\abold$ as above, the proof boils down to verifying
the  analytic identity
\begin{equation}\label{equation:ThetaRecipeForTateTheta}
\rho_{K/k}(\abold)\theta(0)=\sum_{a\in
\frac{\adeles^\times}{k^\times
\UUU\langle
\tau_\infty\rangle}}\left(
\sum_{n\in\ZZ}\Theta(a^{-1}\tau_\infty^n
,\one_\VVV^{(\abold)})\right)\rho_{K/k}(a)
\end{equation}
relating partial zeta values to the theta symbol.
 Since the set of $\abold\in \JJ$ satisfying the identity in
question forms an ideal of $\ZZ[\Gamma]$, we may assume without loss of
generality that
$\abold=1-q_{v_0}\cdot \tau_{v_0}$ for some $v_0\in \overline{S}$.
Let
$\chi:G\rightarrow\CC^\times$ be any character, and extend $\chi$ in
$\CC$-linear fashion to a ring homomorphism
$\chi:\CC[G]\rightarrow\CC$. It is enough to prove
the $\chi$-version of (\ref{equation:ThetaRecipeForTateTheta}),
i.~e., the relation obtained by applying
$\chi$ to both sides of
(\ref{equation:ThetaRecipeForTateTheta}).
Now we adapt to the present case
the method presented in Tate's thesis for meromorphically continuing abelian $L$-functions.
At least for $\Re(s)>1$, when we have absolute convergence
and can freely exchange limit processes, we have
$$
\begin{array}{cl}
&(1-q_{v_0}^{1-s}\chi(F_{v_0}))\cdot\mu^\times
\UUU\cdot\chi(\theta(s))\\\\ =&\displaystyle
(1-q_{v_0}^{1-s}\chi(F_{v_0}))\cdot\mu^\times\UUU\cdot\prod_{v\not\in S_0}
(1-\chi(F_v)q_v^{-s})^{-1}\\\\ =&\displaystyle
(1-q_{v_0}^{1-s}\chi(F_{v_0}))\cdot\int
\one_\VVV(a)\chi(\rho_{K/k}(a))
\norm{a}^sd\mu^\times(a)\\\\
=&\displaystyle
\int
\one_\VVV(a)(\chi(\rho_{K/k}(x))\norm{a}^s-q_{v_0}\chi(\rho_{K/k}(\tau_{v_0}a))
\norm{\tau_{v_0}a}^s)d\mu^\times(a)\\\\
=&\displaystyle
\int
(\one_\VVV-q_{v_0}\one_{\tau_{v_0}\VVV})(a)\cdot\chi(\rho_{K/k}(a))
\norm{a}^sd\mu^\times(a)\\\\
=&\displaystyle
\int_{\adeles^\times/k^\times}
\Theta(a^{-1},\one_{\VVV}-q_{v_0}\one_{\tau_{v_0}\VVV})\chi(\rho_{K/k}(a))
\norm{a}^sd\mu^\times(a).
\end{array}
$$
The last integrand is constant on
cosets of $(k^\times \UUU)/k^\times$ and by Proposition~\ref{Proposition:ThetaAsymptotics} has compact support in $\adeles^\times/k^\times$.
So
the last integral defines an entire function of $s$.
Now since $\UUU\cap k^\times=\{1\}$ we have
$\mu(\UUU)=\mu((k^\times \UUU)/k^\times)$.
So by plugging in $s=0$ at beginning and end of the calculation 
above and breaking the last integral down as a sum of integrals over
cosets of $(k^\times \UUU)/k^\times$ in $\adeles^\times/k^\times$, we
recover the
$\chi$-version of (\ref{equation:ThetaRecipeForTateTheta})
in the case $\abold=1-q_{v_0}\cdot \tau_{v_0}$, as desired.
\qed

\subsection{End
of the
proof}\label{subsection:EndOfRecipeProof} By the preceding two lemmas
$\epsilon(\abold_\infty)$ indeed has the properties
specified in \S\ref{subsubsection:StarkUnitSpecs}
characterizing $\epsilon(K/k,S)$.
\qed


\begin{thebibliography}{999999}


\bibitem[An92]{AndersonStickelberger}
G.\ W.\ Anderson, {\em A two-dimensional analogue of Stickelberger's theorem,}
in: The arithmetic of function fields (Columbus, OH, 1991), 51--73,
Ohio State Univ.\ Math.\ Res.\ Inst.\ Publ., 2, de Gruyter, Berlin, 1992.

\bibitem[An94]{AndersonAHarmonic}
G.\ W.\ Anderson, {\em Rank one elliptic $A$-modules and $A$-harmonic series,}
Duke Math.\ J.\ \textbf{73}(1994), 491--542.

\bibitem[An96]{AndersonSpecialPoints}
G.\ W.\ Anderson, {\em Log-algebraicity of twisted $A$-harmonic series 
and special values of $L$-series in characteristic $p$,}
J.\ Number Theory \textbf{60}(1996),  165--209.

\bibitem[ABP04]{ABP}
G.\ W.\ Anderson, W.\ D.\ Brownawell, M.\ A.\ Papanikolas, 
{\em Determination of the algebraic relations among special 
$\Gamma$-values in positive characteristic,}
Ann.\ of Math.\ (2) \textbf{160}(2004),  237--313.

\bibitem[An04]{AndersonFock}
G.\ W.\ Anderson, 
{\em A generalization of Weyl's identity for $D_n$,}
Adv.\ in Appl.\ Math.\ \textbf{33}(2004), 573--614.

\bibitem[AT]{ArtinTate}
E.\ Artin, J.\ Tate, 
{\em Class field theory.}
Second edition. Advanced Book Classics.
Addison-Wesley Publishing Company, Advanced Book Program, Redwood City, CA, 1990. 

\bibitem[Co88]{Coleman}
R.\ F.\ Coleman, {\em On the Frobenius endomorphisms of Fermat and Artin-Schreier curves,}
Proc.\ Amer.\ Math.\ Soc.\ \textbf{102}(1988),  463--466.


\bibitem[Goss]{Goss} 
D.\ Goss, 
{\em Basic structures of function field arithmetic.}
Ergebnisse der Mathematik und ihrer Grenzgebiete (3) [Results in Mathematics and
Related Areas (3)], \textbf{35}, Springer-Verlag, Berlin, 1996. 

\bibitem[Ha85]{Hayes} D.\ R.\ Hayes, 
{\em Stickelberger elements in function fields,} 
Compositio Math.\ \textbf{55}(1985), 209--239.

\bibitem[Mac]{Macdonald}
I.\ G.\ Macdonald, {\em Symmetric functions and Hall
polynomials.} Second edition, Oxford Mathematical Monographs, Oxford
University Press 1995.

\bibitem[RV]{Ramakrishnan} D.\ Ramakrishnan, R.\ J.\ Valenza,
{\em Fourier Analysis on Number Fields,}
Graduate Texts in Mathematics \textbf{186}, Springer, New York 1999.

\bibitem[Serre]{Serre} 
J.-P.\ Serre, {\em Algebraic groups and class fields,} Graduate Texts in Math.\ \textbf{117},
Springer, New York 1988.

\bibitem[Si97a]{Sinha}
S.\ Sinha, {\em Periods of $t$-motives and transcendence}, Duke
Math.\ J.\ {\bf 88} (1997), 465--535.

\bibitem[Si97b]{SinhaDR}
S.\ Sinha, {\em Deligne's reciprocity for function fields,} J.\ Number Theory {\bf 63}(1997),
 65--88.

\bibitem[St80]{Stark} H.\ Stark, {\em $L$-functions at $s=1$.
First derivatives at $s=0$.} Adv.\ in Math.\ \textbf{35}(1980), 197-235.
\bibitem[Ta79]{TateBackground}
J.\ Tate, {\em Number-theoretic background},
Proc.\ of Symp.\ in Pure Math.\ \textbf{33}(1979), part 2, pp.\ 3-26.
\bibitem[Tate]{Tate}
J. Tate, {\em Les Conjectures de Stark sur les
Fonctions $L$ d'Artin en $s=0$,}
Prog.\ in Math.\ {\bf 47}, Birkh\"{a}user Boston 1984.

\bibitem[Th91]{Thakur}
D.\ S.\ Thakur, {\em Gamma functions for function fields and
Drinfeld modules,} Ann.\ of Math.\ {\bf 134} (1991), 25--64.

\bibitem[Thak]{ThakurBook} D.\ S.\ Thakur, {\em Function Field Arithmetic,}
World Scientific Pub., 2004.

\end{thebibliography}
\end{document}